\documentclass{article}

\usepackage[utf8]{inputenc}
\usepackage[english]{babel}
\usepackage[T1]{fontenc}
\usepackage{graphicx}
\usepackage{tikz}
\usepackage{quiver}
\usepackage{float}
\usepackage{url} 
\usepackage{enumerate} 
\usepackage{enumitem}
\usepackage{color}
\usepackage{hyperref}
\usepackage{lmodern} 
\usepackage{amssymb}
\usepackage{amsmath}
\usepackage{amsfonts}
\usepackage{mathrsfs}
\usepackage{amsthm}
\usepackage{mathtools}
\usepackage{titlesec}
\usepackage[top=4cm, bottom=4cm, left=3cm, right=3cm]{geometry} 
\usepackage{quiver}
\usepackage{xparse}
\usepackage{algpseudocode}

\newtheoremstyle{problem}
{\topsep}
{\topsep}
{\tt}
{0pt}
{\sc}
{---}
{ }
{\textsc{\thmnote{ #3 }}}

\newtheorem{The}{Theorem} 
\newtheorem{Pro}{Proposition}

\theoremstyle{definition}

\newtheorem{Def}{Definition}
\newtheorem{Le}{Lemma}
\newtheorem*{Ex}{Example}
\newtheorem*{Rem}{Remark}
\newtheorem*{Cor}{Corollary} 
\newtheorem*{Fac}{Fact} 
 
\newtheorem*{Not}{Notation} 
\newtheorem*{Rec}{Recall} 

\theoremstyle{problem}

\makeatletter
\newcommand*{\delsmo}{\mathrel{\mathpalette\reflectmathsymbol\models}}
\newcommand*{\reflectmathsymbol}[2]{\reflectbox{$\m@th#1#2$}}

\title{Polyteam logics for quantum contextuality} 
\author{Loan Le Maguet \\ \small ENS de Lyon \\ \small \texttt{loan.le\_maguet@ens-lyon.fr}}
\date{}

\usepackage[backend=biber,url=false,doi=false,isbn=false,date=year,style=numeric]{biblatex}
\begin{document}
\newcommand{\dom}{\textnormal{dom}}
\newcommand{\ran}{\textnormal{ran}}
\newcommand{\supp}{\textnormal{Supp}}
\newcommand{\iffdef}{\overset{\textnormal{def}}{\iff}}
\newcommand{\niff}{\iff \hspace{-17pt} / \hspace{12pt}}
\newcommand{\texists}{\tilde \exists}
\newcommand{\tforall}{\tilde \forall}
\newcommand{\verteq}{\rotatebox{90}{$\,=$}}
\newcommand{\leftlongmapsto}{\rotatebox{180}{$\longmapsto$}}
\newcommand\indep{\protect\mathpalette{\protect\independenT}{\perp}}
\def\independenT#1#2{\mathrel{\rlap{$#1#2$}\mkern2mu{#1#2}}}
\NewDocumentCommand{\independent}{ O{} O{} }{\indep^{\!\!\! #1}_{\!\! #2}\ }
\NewDocumentCommand{\dep}{ O{} }{=^{\! #1} \!}
\maketitle

\vspace{-20pt}
\begin{abstract}
We introduce logics of dependence and independence over polyteams (\textit{i.e.} databases with multiple relation schemes), in a similar way as in \cite{hannula_polyteam_2020}. On this basis, we develop a logical framework in which we can express (in formulas) relevant properties about the statistics of a repeated quantum experiment. In this framework, we formalize Bell's theorem, and the Kochen-Specker theorem as failure of logical entailments, as it was suggested in \cite{abramsky_team_2026}. The key feature of polyteam logics in action here is their power to express both local and global properties, with respect to the set of possible measurements. The point of this is at least twofold. Firstly, as in \cite{abramsky_sheaf-theoretic_2011}, it gives a high-level (information-theoretic) understanding of these physical phenomena, that is not intrinsic to the physical objects involved. Secondly, since logical entailment is often (partially) axiomatizable, this framework incites us to tackle the question of contextuality computationally. We also define semiring-weighted polyteams along with logics on them, and use them to recover the probabilistic versions of the contextuality theorems mentioned above. 
\end{abstract}

\section*{Introduction}

\paragraph*{}
Quantum physics sometimes prevents us from measuring a set $M$ of attributes simultaneously, \textit{e.g.} Heisenberg's inequality for position and momentum. Thus, we are only left with local information : for each compatible $C \subseteq M$, we may determine through experiment the statistical behavior of the attributes $C$ on a given system. This yields a set of tables, \textit{i.e.} essentially a so called polyteam. Even though it is inaccessible to experiment, it is natural to think this local information is derived (by repeated projections) from one big table : the one describing the statistical behavior of the whole set of attributes $M$. Quantum contextuality is the realization that this natural assumption turns out to be false, and the study of counter-examples. 

\paragraph*{}
This work tackles quantum contextuality through the lens of polyteam logics. It is the purpose of Section \ref{sec:definitions} to define teams, polyteams (Section \ref{sec:def_polyteams}), polyteam logics (Sections \ref{sec:polyteam_atoms} and \ref{sec:def_polyteam_logics}), and a polyteam-based framework for representing quantum experiments (Section \ref{sec:polyteam_representation_experiments}). A logical analysis of quantum contextuality is interesting for several reasons : it is a higher level (information-theoretic) understanding of the phenomena (hence, more general and less intrinsic), and the restricted syntax prompts us into searching for axiomatizations of logical entailment. 

\paragraph*{}
(Poly)team logics feature atomic formulas such as $\dep(x,y)$ "$x$ determines $y$", or $x \perp y$ "$x$ is independent from $y$". These (in)dependence notions are very common in experimental sciences, hence we may expect polyteam formulas to reflect natural scientific discourse. It is the purpose of Section \ref{sec:theorems} to show most of the relevant physical properties of quantum experiments are expressible in polyteam formulas (Sections \ref{sec:physical_properties} and \ref{sec:global_properties_formulas}), and that Bell's theorem and the Kochen-Specker theorem are formalizable as failure of entailments in polyteam logics (Section \ref{sec:formalization}). There was already work in this area, and specifically about providing \emph{team}-based accounts of contextuality. Section \ref{sec:previous_work} discusses the relation between this contribution and previous literature. In short, we believe our polyteam account is at the same time efficient (formulas are somewhat expressive), consistent (no need to add alien atoms to the natural polyteam logics) and computation-friendly (polyteams are like databases). 

\paragraph*{}
In Section \ref{sec:semirings}, we introduce semiring-weighted polyteams and logics of (in)dependence on them. Here again, physical properties are expressible by formulas (Section \ref{sec:R_properties_meaning}). They allow us to formalize contextuality theorems in their original probabilistic version, in Section \ref{sec:R_entailments}. 

\paragraph*{}
We leave computational questions for future work.

\newpage
\section{Polyteam logics, and representation of quantum experiments} \label{sec:definitions}

\subsection{Structures, teams, and polyteams} \label{sec:def_polyteams}
Let us fix a first-order relational language $\tau$, possibly many-sorted. It contains the following. 
\begin{itemize}
\item Symbols for sorts : $i,j,k, \ldots$
\item Countably many variables for each sort : $x_0^i, x_1^i,x_2^i, \cdots$, $x_0^j, x_1^j, x_2^j, \cdots$, $x_0^k, x_1^k, x_2^k, \cdots$, $\ldots$
\item Symbols for relations : $R_0, R_1, R_2, \ldots$ 
\item For each relation symbol $R$, its arity and the tuple of sorts it takes as input. 
\item Sorted constant symbols, if needed. 
\end{itemize}

\begin{Def}
$\tau'$ is called an \emph{extension of $\tau$ by the sort $\ell$} whenever it has the additional sort $\ell$, the associated variables, and that all the rest is exactly the same as $\tau$. 
\end{Def}

\begin{Def}
A $\tau$-structure $\mathcal A$ is given by the following information : 
\begin{itemize}
\item For each sort $i$, an associated domain $A_i \neq \emptyset$
\item For each relation symbol $R$ with sort-domain $(i_0, \cdots, i_{k-1})$, a set $R^{\mathcal A} \subseteq A_{i_0} \times \cdots \times A_{i_{k-1}}$
\item For each constant symbol $c^i$, an element $(c^i)^{\mathcal A} \in A_i$
\end{itemize}
A $\tau'$-structure $\mathcal B$ is called an \emph{expansion of $\mathcal A$ by the sort $\ell$} whenever its domains for $i \neq \ell$, and its interpretations for relations and constants are the same as that of $\mathcal A$.
\end{Def}

\begin{Not}
Let us denote $A$ for the collection of elements that appear in at least one of the $A_i$'s.
\end{Not}

\begin{Rem}
Unless otherwise specified, let us assume the language $\tau$ has no constant symbols (they are harmless, but also useless).
\end{Rem}

\begin{Def}
Let $\mathcal A$ be a $\tau$-structure, and $\bar x$ be a set of variables (with arbitrary sorts). An \emph{assignment} on $\mathcal A$ with domain $\bar x$ is a map $s : \bar x \to A$ such that each $x^i \in \bar x$ is sent in $A_i$. A \emph{team} on $\mathcal A$ with domain $\bar x$ is a set of such assignments. 
\end{Def}

\begin{Def}
A \emph{labeling} of $\bar x$ is a surjective map $\chi : \bar x \to \mathcal C$ where $\mathcal C$ is a set. A \emph{polyteam} on $\mathcal A$ with domain $\bar x$ labeled as $\chi : \bar x \to \mathcal C$ is a set of teams $\{L_C, C \in \mathcal C\}$ on $\mathcal A$ such that $\dom(L_C) = \{v \in \bar x : \chi(v) = C\}$. 
\end{Def}

The role of the labeling map is to ensure the $\dom(L_C)$'s are disjoint, which is crucial for the upcoming semantics to be unambiguous. \cite{hannula_polyteam_2020} reaches the same result by forcing the $\dom(L_C)$'s to contain variables of different sorts, \textit{i.e.} the authors use sorts as labels. In this paper, we allow multiple sorts to inhabit the same $\dom(L_C)$, in order to later use the notion of sorts in a more intrinsic way (in the semantics of $\texists, \tforall$). 

\begin{Rem}
A (poly)team may be infinite, but we will focus on the finite case. Up to bijection, its set of labels $\mathcal C$ is unique. It neither has ordering, occurrence, nor probability (yet). We use tuples or sets of variables interchangeably (for convenience), since it has no impact on the content of our study.
\end{Rem}

\begin{Not}
In the following, teams will typically be denoted by $X, Y$ or $Z$, and polyteams by $\mathcal L$ or $\mathcal L'$ (as they represent \emph{local} information). We write $x^C$ to indicate $x$ is labeled $C$. To avoid confusion with sorts, say by convention that labels are capital letters whereas sorts are lowercase. 
\end{Not}

\subsection{Dependence and independence atoms for polyteams} \label{sec:polyteam_atoms}

Polyteam semantics were introduced in \cite{hannula_polyteam_2020}. We stick to it in this section, with labels instead of sorts. 

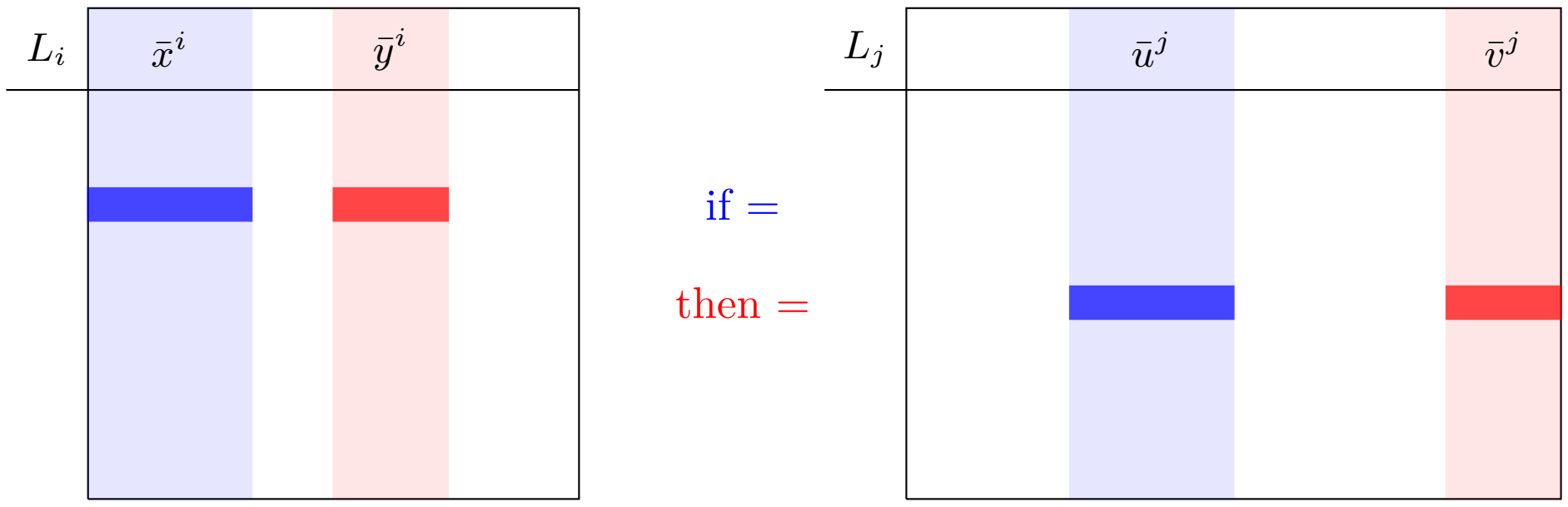
\begin{figure}
\centering
\begin{tikzpicture}[scale = 1.3]
\draw (0,0) -- (3,0) -- (3,3) -- (0,3) -- (0,0);
\draw (-0.5,2.5) -- (3,2.5);
\node at (-0.25,2.75) {\small $L_i$};
\filldraw[color=blue,opacity=0.1] (0,0) rectangle (1,3);
\node at (0.5,2.75) {\small $\bar x^i$}; 
\filldraw[color=red,opacity=0.1] (1.5,0) rectangle (2.2,3);
\node at (1.85,2.75) {\small $\bar y^i$}; 

\draw (5,0) -- (9,0) -- (9,3) -- (5,3) -- (5,0); 
\draw (4.5,2.5) -- (9,2.5); 
\node at (4.75,2.75) {\small $L_j$};
\filldraw[color=blue,opacity=0.1] (6,0) rectangle (7,3);
\node at (6.5,2.75) {\small $\bar u^j$}; 
\filldraw[color=red,opacity=0.1] (8.3,0) rectangle (9,3);
\node at (8.65,2.75) {\small $\bar v^j$};

\filldraw[color=blue,opacity=0.7] (0,1.7) rectangle (1,1.9);
\filldraw[color=blue,opacity=0.7] (6,1.1) rectangle (7,1.3);
\filldraw[color=red,opacity=0.7] (1.5,1.7) rectangle (2.2,1.9);
\filldraw[color=red,opacity=0.7] (8.3,1.1) rectangle (9,1.3);
\node[color=blue] at (4,1.8) {if $=$};
\node[color=red] at (4,1.2) {then $=$};
\end{tikzpicture}
\caption{Semantics of polydependence atoms}
\label{fig:polydependence_atom}
\end{figure}

\begin{Def} 
Let $\bar x^C, \bar y^C, \bar u^{C'}, \bar v^{C'}$ be tuples of labeled variables from $\tau$ such that $|\bar x| = |\bar u|$ and $|\bar y| = |\bar v|$. The syntactical expression $\dep(\bar x^C, \bar y^C / \bar u^{C'}, \bar v^{C'})$ is called a \emph{polydependence atom} (for polyteams). It has the following semantics (see Figure \ref{fig:polydependence_atom}). 
$$\mathcal A \models_{\mathcal L}\ \dep(\bar x^C, \bar y^C / \bar u^{C'}, \bar v^{C'}) \iffdef \forall s \in L_C, \forall s' \in L_{C'}, s(\bar x^C) = s'(\bar u^{C'}) \to s(\bar y^C) = s'(\bar v^{C'})$$
\end{Def}

\begin{Not}
We keep this notation in this article to remain consistent with previous literature on (poly)team semantics. It seems that the notation $\bar x^C / \bar u^{C'} \to \bar y^C / \bar v^{C'}$ sticks better to its meaning. 
\end{Not}

\begin{Rem}
$\dep(\bar x^C, \bar y^C / \bar x^C, \bar y^C)$ is true with respect to $\mathcal L \iff \dep(\bar x^C, \bar y^C)$ is true with respect to $L_C$, according to \emph{team} semantics. $\dep(\emptyset^C, \bar y^C / \emptyset^{C'}, \bar v^{C'})$ is true with respect to $\mathcal L$ if and only if $\bar y^C$ and $\bar v^{C'}$ are constant, equal to the same value, in $L_C$ and $L_{C'}$ respectively. 
\end{Rem}

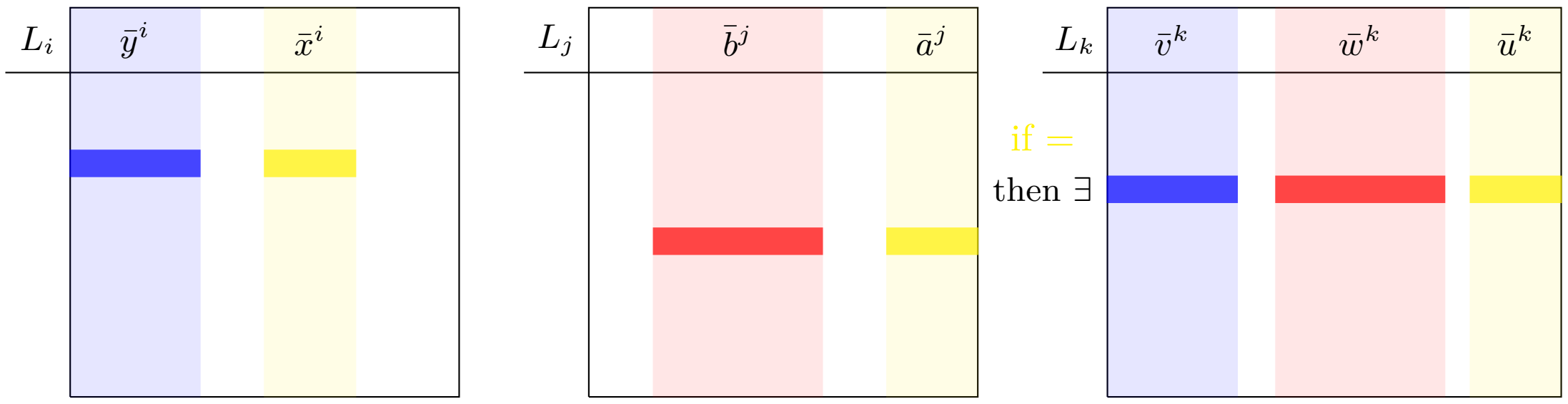
\begin{figure}
\centering
\begin{tikzpicture}[scale = 1.2]
\draw (0,0) -- (3,0) -- (3,3) -- (0,3) -- (0,0);
\draw (-0.5,2.5) -- (3,2.5);
\node at (-0.25,2.75) {\small $L_i$};
\filldraw[color=blue,opacity=0.1] (0,0) rectangle (1,3);
\node at (0.5,2.75) {\small $\bar y^i$}; 
\filldraw[color=yellow,opacity=0.1] (1.5,0) rectangle (2.2,3);
\node at (1.85,2.75) {\small $\bar x^i$}; 

\draw (4,0) -- (7,0) -- (7,3) -- (4,3) -- (4,0); 
\draw (3.5,2.5) -- (7,2.5); 
\node at (3.75,2.75) {\small $L_j$};
\filldraw[color=red,opacity=0.1] (4.5,0) rectangle (5.8,3);
\node at (5.15,2.75) {\small $\bar b^j$}; 
\filldraw[color=yellow,opacity=0.1] (6.3,0) rectangle (7,3);
\node at (6.65,2.75) {\small $\bar a^j$};

\draw (8,0) -- (11.5,0) -- (11.5,3) -- (8,3) -- (8,0); 
\draw (7.5,2.5) -- (11.5,2.5); 
\node at (7.75,2.75) {\small $L_k$};
\filldraw[color=blue,opacity=0.1] (8,0) rectangle (9,3);
\node at (8.5,2.75) {\small $\bar v^k$}; 
\filldraw[color=red,opacity=0.1] (9.3,0) rectangle (10.6,3);
\node at (9.95,2.75) {\small $\bar w^k$}; 
\filldraw[color=yellow,opacity=0.1] (10.8,0) rectangle (11.5,3);
\node at (11.15,2.75) {\small $\bar u^k$};

\filldraw[color=blue,opacity=0.7] (0,1.7) rectangle (1,1.9);
\filldraw[color=red,opacity=0.7] (4.5,1.1) rectangle (5.8,1.3);
\filldraw[color=yellow,opacity=0.7] (1.5,1.7) rectangle (2.2,1.9);
\filldraw[color=yellow,opacity=0.7] (6.3,1.1) rectangle (7,1.3);

\filldraw[color=blue,opacity=0.7] (8,1.5) rectangle (9,1.7);
\filldraw[color=red,opacity=0.7] (9.3,1.5) rectangle (10.6,1.7);
\filldraw[color=yellow,opacity=0.7] (10.8,1.5) rectangle (11.5,1.7);

\node[color=yellow] at (7.5,2) {if $=$}; 
\node at (7.5,1.6) {\small{then $\exists$}}; 
\end{tikzpicture}
\caption{Semantics of polyindependence atoms}
\label{fig:polyindependence_atom}
\end{figure}

\begin{Def}
Consider tuples of labeled variables from $\tau$ such that $|\bar x^C| = |\bar a^{C'}| = |\bar u^{C''}|$, $|\bar y^C| = |\bar v^{C''}|$ and $|\bar b^{C'}| = |\bar w^{C''}|$. The syntactical expression $\bar y^C / \bar v^{C''} \perp_{\bar x^C, \bar a^{C'} / \bar u^{C''}} \bar b^{C'} / \bar w^{C''}$ is called a \emph{polyindependence atom} (for polyteams). It has the following semantics (see Figure \ref{fig:polyindependence_atom}).
$$\hspace{-35pt}\mathcal A \models_{\mathcal L} \bar y^C / \bar v^{C''} \perp_{\bar x^C, \bar a^{C'} / \bar u^{C''}} \bar b^{C'} / \bar w^{C''} \iffdef \forall s \in L_C, \forall s' \in L_{C'}, s(\bar x^C) = s'(\bar a^{C'}) \to \exists s'' \in L_{C''}\ \begin{array}[t]{|l} s''(\bar u^{C''}) = s'(\bar a^{C'}) = s(\bar x^C),\\ s''(\bar v^{C''}) = s(\bar y^C), \text{ and} \\ s''(\bar w^{C''}) = s'(\bar b^{C'}) \end{array}$$
\end{Def}

\begin{Fac}
\begin{itemize}
\item Polyindependence atoms are symmetric, and polydependence atoms are symmetric with respect to the slash. 
\item Polydependence atoms are transitive, \textit{i.e.} if $\mathcal L \models\ \dep(\bar x^C, \bar y^C / \bar u^{C'}, \bar v^{C'})$ and $\mathcal L \models\ \dep(\bar y^C, \bar z^C / \bar v^{C'}, \bar w^{C'})$, then $\mathcal L \models\ \dep(\bar x^C, \bar z^C / \bar u^{C'}, \bar w^{C'})$. 
\end{itemize}
\end{Fac}

\begin{proof}
Follows easily from the definitions. 
\end{proof}

\begin{Not}
$\bar y^C / \bar v^{C''} \perp_{\emptyset^C, \emptyset^{C'} / \emptyset^{C''}} \bar b^{C'} / \bar w^{C''}$ is simply denoted $\bar y^C / \bar v^{C''} \perp \bar b^{C'} / \bar w^{C''}$ and called a \emph{pure polyindependence atom}. 
\end{Not}

\begin{Le} \label{le:polydep_atoms_as_polyindep_atoms}
Any polydependence atom is equivalent to a polyindependence atom : 
$$\dep(\bar x^C, \bar y^C / \bar u^{C'}, \bar v^{C'}) \equiv \bar y^C / \bar y^C \perp_{\bar x^C, \bar u^{C'} / \bar x^C} \bar v^{C'} / \bar y^C$$
\end{Le}

\begin{proof}
\begin{itemize}
\item[$\models$.] Assume $\mathcal L \models\ \dep(\bar x^C, \bar y^C / \bar u^{C'}, \bar v^{C'})$, and that $s \in L_C$, $s' \in L_{C'}$ are so that $s(\bar x) = s'(\bar u)$. Then $s(\bar y) = s'(\bar v)$, and it follows that $s$ witnesses $\mathcal L \models \bar y^C / \bar y^C \perp_{\bar x^C, \bar u^{C'} / \bar x^C} \bar v^{C'} / \bar y^C$ for $s,s'$. 
\item[$\delsmo$.] Assume $\mathcal L \models \bar y^C / \bar y^C \perp_{\bar x^C, \bar u^{C'} / \bar x^C} \bar v^{C'} / \bar y^C$, and that $s \in L_C$, $s' \in L_{C'}$ are so that $s(\bar x) = s'(\bar u)$. There is $s'' \in L_C$ such that $s''(\bar x) = s'(\bar u) = s(\bar x)$, $s''(\bar y) = s(\bar y)$ and $s''(\bar y) = s'(\bar v)$. In particular $s(\bar y) = s'(\bar v)$. 
\end{itemize}
\end{proof}

\begin{Pro}
There is a complete proof system for polydependence atoms. 
\end{Pro}

\begin{proof}
See Definition 3 and Theorem 4 in \cite{hannula_polyteam_2020}. 
\end{proof}

One may also define \emph{polyinclusion atoms} as $\mathcal L \models \bar x^C \subseteq \bar y^{C'} \iff \forall s \in L_C, \exists s' \in L_{C'}, s(\bar x^C) = s'(\bar y^{C'})$ and \emph{polyexclusion atoms} as $\mathcal L \models \bar x^C\ |\ \bar y^{C'} \iff \forall s \in L_C, \forall s' \in L_{C'}, s(\bar x^C) \neq s'(\bar y^{C'})$. They are not necessary to understand what follows, although they may allow slight simplifications.  

\subsection{Dependence and independence logics for polyteams} \label{sec:def_polyteam_logics}

\begin{Def}[Syntax]
The formulas of polyteam logics are built according to the following syntax. 
\begin{center}
$\varphi := $ Single-labeled first-order literal | Non-first order atom | $\psi \wedge \theta$ | $\psi \vee \theta$ | $\exists x\ \psi$ | $\forall x\ \psi$
\end{center}
By "single-labeled first-order literal", we mean an atom involving variables from only one context, which is either an equality atom, a relational atom, or the negation of any of these. By non-first order atom, we refer to any of the atoms defined in previous section. The set of formulas thus obtained is denoted accordingly to the set of polyatoms allowed. In particular FO($=$) stands for \emph{polydependence logic}, and FO($\perp$) stands for \emph{polyindependence logic}. These logics must not be mistaken for their \emph{team}-logic analogues : in ambiguous contexts, we index them with "poly". 
\end{Def}

\begin{Rem}
Unless otherwise specified, polyteam logics have no negation operator.
\end{Rem}

\begin{Def}
Let $X$ be a team on $\mathcal A$, $x \notin \dom(X)$ with sort $i$, and $F : X \to \mathcal P^+(A_i)$ a function mapping each assignment in $X$ to a non-empty set of $A_i$-values. The \emph{supplement team of $X$ over $x$ with $F$} is the following. 
$$X(F/x) := \bigcup_{s \in X} \{s(a/x), a \in F(s)\},\ \text{ where } 
s(a/x) :\begin{array}[t]{rcl}
\dom(X) \cup \{x\} & \to & A \\ 
y & \mapsto & \begin{cases} a & \text{if } y \text{ is } x \\ s(y) & \text{otherwise} \end{cases}
\end{array}$$
When $F$ is the constant $A_i$, we denote $X(A_i/x)$ or $X(A/x)$ and we call it the \emph{duplicate team of $X$ over $x$}. 
\end{Def}

\begin{Def}
Let $\mathcal L$ be a polyteam on $\mathcal A$, $x^{C,i}$ a variable from $\tau$, and $F : L_C \to \mathcal P^+(A_i)$. The \emph{supplement polyteam of $\mathcal L$ over $x$ with $F$} is the polyteam $\mathcal L(F/x)$ obtained from $\mathcal L$ by replacing $L_C$ with $L_C(F/x)$. It is called \emph{duplicate polyteam} $\mathcal L(A/x)$ whenever $F$ is the constant $A_i$. 
\end{Def}

\begin{Ex}
Consider the polyteam $\mathcal L$ explicited on the left, on a 3-sorted structure that has domains $A_i = \{0,1,2\}$, $A_j = \{0,1\}$ and $A_k = \{0,1,2\}$. The polyteam on the right is $\mathcal L(F/z^{\chi(x),i})$, where $F$ maps the assignments of $L_{\chi(x)}$ to $\{0\}$, $\{0,2\}$, $\{0,1,2\}$ in that order. 

\begin{center}
$\begin{array}{cc}
x^i & y^i \\ \hline
0 & 0 \\ 
0 & 1 \\ 
2 & 0 
\end{array}$
\hspace{20pt}
$\begin{array}{cc}
a^j & b^k \\ \hline
0 & 0 \\ 
0 & 1 \\ 
1 & 1 
\end{array}$
\hspace{20pt}
$\begin{array}{c}
u^k \\ \hline
1 \\ 
2
\end{array}$
\hspace{20pt}
$\ast$
\hspace{20pt}
$\begin{array}{ccc}
x^i & y^i & z^i \\ \hline
0 & 0 & 0 \\ 
0 & 1 & 0 \\ 
0 & 1 & 2 \\ 
2 & 0 & 0 \\ 
2 & 0 & 1 \\ 
2 & 0 & 2
\end{array}$
\hspace{20pt}
$\begin{array}{cc}
a^j & b^k \\ \hline
0 & 0 \\ 
0 & 1 \\ 
1 & 1 
\end{array}$
\hspace{20pt}
$\begin{array}{c}
u^k \\ \hline
1 \\ 
2
\end{array}$
\end{center}
\end{Ex}

\begin{Rem}
We also define $X(F/x)$ and $X(A/x)$ when $x \in \dom(X)$, and by the same token $\mathcal L(F/x)$ and $\mathcal L(A/x)$, by first forgetting the $x$-information brought by $X$ and then supplementing.
\end{Rem}

\begin{Def}
We say $(\mathcal L_1, \mathcal L_2)$ is a \emph{split} of the polyteam $\mathcal L$, denoted $\mathcal L_1 \cup \mathcal L_2 = \mathcal L$, whenever $L_{1,C} \cup L_{2,C} = L_C$ for all $C$, \textit{i.e.} whenever it is a teamwise split. 
\end{Def}

\begin{Def}[Semantics] \label{def:polyteam_semantics}
Satisfaction of polyteam logics formulas in a $\tau$-structure $\mathcal A$ is so defined.
\begin{itemize}
\item[(First order literal)] $\mathcal A \models_{\mathcal L}\ \alpha \iffdef \forall s \in L_C, \mathcal A \models_s\ \alpha$, where $C$ is the uniform label of FV$(\alpha)$
\item[(Non-first order atom)] As detailed in section \ref{sec:polyteam_atoms}.
\item[($\wedge$)] $\mathcal A \models_{\mathcal L}\ \varphi \wedge \psi \iffdef \mathcal A \models_{\mathcal L}\ \varphi$ and $\mathcal A \models_{\mathcal L}\ \psi$
\item[($\vee$)] $\mathcal A \models_{\mathcal L}\ \varphi \vee \psi \iffdef \exists \mathcal L_1, \mathcal L_2 \subseteq \mathcal L,\ \mathcal L_1 \cup \mathcal L_2 = \mathcal L$, $\mathcal A \models_{\mathcal L_1}\ \varphi$ and $\mathcal A \models_{\mathcal L_2}\ \psi$
\item[($\exists$)] $\mathcal A \models_{\mathcal L}\ \exists x^{C,i}\ \varphi \iffdef \exists F : L_C \to \mathcal P^+(A_i),\ \mathcal A \models_{\mathcal L(F/x)}\ \varphi$
\item[($\forall$)] $\mathcal A \models_{\mathcal L}\ \forall x^{C,i}\ \varphi \iffdef \mathcal A \models_{\mathcal L(A/x)}\ \varphi$ 
\end{itemize}
In section \ref{sec:theorems}, we will moreover use the quantifiers $\texists, \tforall$ with the following semantics. Let $\tau'$ extend $\tau$ by a sort $h$, and $z^{C,h}$ a variable of that sort in $\tau'$. 
\begin{itemize}
\item[($\texists$)] $\mathcal A \models_{\mathcal L}\ \texists z^{C,h}\ \varphi \iffdef \mathcal B \models_{\mathcal L}\ \exists z^{C,h}\ \varphi$ for some $\mathcal B$ expanding $\mathcal A$ by the sort $h$
\item[($\tforall$)] $\mathcal A \models_{\mathcal L}\ \tforall z^{C,h}\ \varphi \iffdef \mathcal B \models_{\mathcal L}\ \forall z^{C,h}\ \varphi$ for all $\mathcal B$ expanding $\mathcal A$ by the sort $h$
\end{itemize}
\end{Def}

\begin{Rem}
This is called the \emph{global} splitjunction semantics for $\vee$. It makes $\vee$ associative and commutative. One can replace $\vee$ with \emph{local} $\vee^C$'s (one per label), whose semantics only allow splits that leave the $L_{C'}$'s ($C' \neq C$) untouched. Notably, it happens that it would yield the same expressive power (see \cite{hannula_polyteam_2020}). We will not use the latter. 
\end{Rem}

\begin{Not}
In FO$(\perp)$, we allow ourselves the use of $\dep(\bar x^C, \bar y^C / \bar u^{C'}, \bar v^{C'})$ for $\bar y^C / \bar y^C \perp_{\bar x^C, \bar u^{C'} / \bar x^C} \bar v^{C'} / \bar y^C$, in line with Lemma \ref{le:polydep_atoms_as_polyindep_atoms}..\end{Not}

\begin{figure}
\centering
\begin{tikzpicture}
\node at (0,0) {$\boxed{\text{FO}}$};
\node at (-1,1) {$\boxed{\text{FO}(=), \text{FO}(|)}$};
\node at (1.3,1) {$\boxed{\text{FO}(\subseteq)}$};
\node at (0,2) {$\boxed{\text{FO}(\perp), \text{FO}(\perp_\emptyset), \text{FO}(=,\subseteq), \text{FO}(|,\subseteq)}$};
\draw (-1,1.27) -- (-1,1.73);
\draw (1.3,1.27) -- (1.3,1.73); 
\draw (-0.3,0.25) -- (-0.3,0.73);
\draw (0.3,0.25) -- (1.3,0.73);`
\node at (-4.5,2) {ESO$(R_1, \cdots, R_n) \equiv$};
\end{tikzpicture}
\caption{Expressive power of $n$-sorted polyteam logics}
\label{fig:expressive_power_polyteams}
\end{figure}
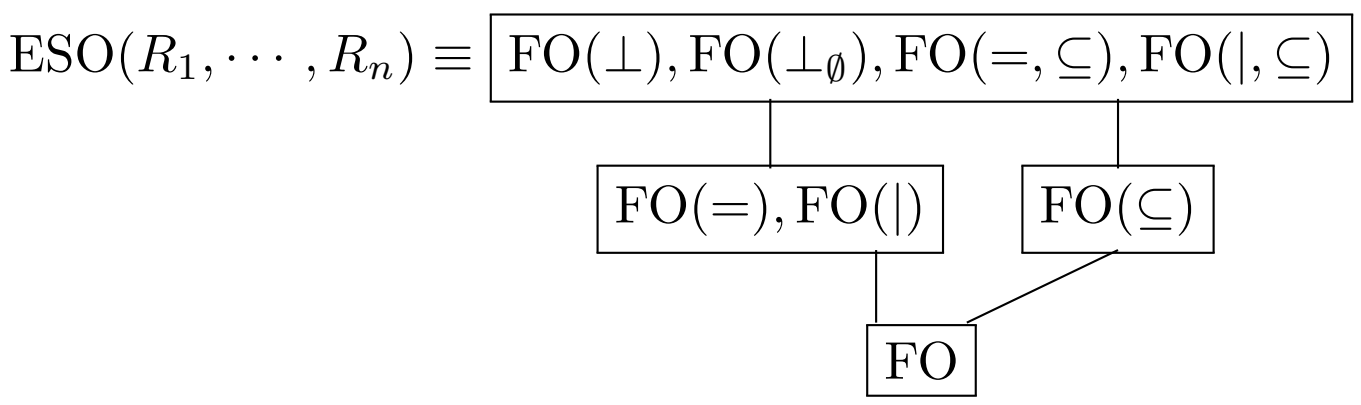

The expressive powers of various polyteam logics are detailed in Figure \ref{fig:expressive_power_polyteams}. Read \cite{hannula_polyteam_2020} for proofs. 

Before giving examples of satisfaction, let us prove the \emph{empty polyteam property} and \emph{domain locality} for polyteams.  

\begin{Pro}[Empty polyteam property]
$\mathcal A \models_{\{\emptyset, \cdots, \emptyset\}} \varphi$ holds, whatever $\mathcal A$ and $\varphi$. 
\end{Pro}

\begin{proof}
The proof goes by induction on $\varphi$. It vacuously holds for atoms. Assume it holds for $\varphi$ and $\psi$.
\begin{itemize}
\item It holds for $\varphi \wedge \psi$ by the very definition of $\wedge$. 
\item It holds for $\varphi \vee \psi$ through the split $\{\emptyset, \cdots, \emptyset\} \cup \{\emptyset, \cdots, \emptyset\} = \{\emptyset, \cdots, \emptyset\}$. 
\item It holds for $\exists x\ \varphi$ and $\texists x\ \varphi$, since it suffices to supplement via the empty map. 
\item It holds for $\forall x\ \varphi$ and $\tforall x\ \varphi$, since the duplicate polyteam remains of the form $\{\emptyset, \cdots, \emptyset\}$. 
\end{itemize} 
\end{proof}

\begin{Not}
$\mathcal L \restriction \bar v := \{L_C \restriction \bar v,\ C \in \text{labels}(\mathcal L)\}$
\end{Not}

\begin{Le} \label{le:locality}
Let $\varphi$ be a formula, and $\mathcal L$ a polyteam on some structure, such that $\dom(\mathcal L) \supseteq$ FV($\varphi$). Then : 
$$\mathcal L \models \varphi \iff \mathcal L \restriction \text{FV}(\varphi) \models \varphi$$
\end{Le}

\begin{proof}
The proof goes by induction on $\varphi$. Atomic cases are easy. Assume (H) it holds for $\varphi$ and $\psi$. 
\begin{itemize}
\item It is easy to show the result for $\varphi \wedge \psi$ if one restricts by FV$(\varphi \wedge \psi)$ instead of FV$(\varphi)$, FV$(\psi)$. 
\item Similarly for $\vee$, and further notice that restricting both sides of a split yields a split of the restriction.
\item $\forall$ case. $\mathcal L \models \forall x\ \varphi \iff \mathcal L(A/x) \models \varphi \overset{\text{H}_\varphi}{\iff} \mathcal L(A/x) \restriction \text{FV}(\varphi) \models \varphi \overset{\text{equal}}{\iff} (\mathcal L \restriction \text{FV}(\varphi))(A/x) \models \varphi\\ \iff \mathcal L \restriction \text{FV}(\varphi) \models \forall x\ \varphi \iff \mathcal L \restriction \text{FV}(\forall x \varphi) \models \forall x\ \varphi$
\item The case $\theta = \exists x^{C,i}\ \varphi$ is trickier. Consider the following claim.  
$$\exists F : L_C \to \mathcal P^+(A_i), \mathcal L(F/x) \restriction \text{FV}(\varphi) \models \varphi \iff \exists G : (\mathcal L \restriction \text{FV}(\varphi))_C \to \mathcal P^+(A_i), (\mathcal L \restriction \text{FV}(\varphi))(G/x) \models \varphi$$ 
If it is true, the proof ends as follows : 
\begin{multline*}
\mathcal L \models \theta \iff \exists F\ \mathcal L(F/x) \models \varphi \overset{\text{H}_\varphi}{\iff} \exists F\ \mathcal L(\mathcal F/x) \restriction \text{FV}(\varphi) \models \varphi \\ \overset{\text{claim}}{\iff} \exists G\ (\mathcal L \restriction \text{FV}(\varphi))(G/x) \models \varphi \iff \mathcal L \restriction \text{FV}(\varphi) \models \exists x\ \varphi \iff \mathcal L \restriction \text{FV}(\theta) \models \theta
\end{multline*}

Remains to prove the claim. For $\Leftarrow$, just consider $F(s) := G(s \restriction \text{FV}(\varphi))$. For $\Rightarrow$, it is clearly sufficient to show that we get $\mathcal L(F/x) \restriction \text{FV}(\varphi) = (\mathcal L \restriction \text{FV}(\varphi))(G/x)$ when $G$ is so defined :  
$$G(s) := \bigcup_{\underset{\text{extending } s}{\hat s \in L_C}} F(\hat s)$$
But indeed, the $C$-labeled team of $\mathcal L(F/x) \restriction \text{FV}(\varphi)$ exactly contains those assignments obtained by augmenting some assignment $s \in L_C$ with an $x$-value taken from $F(s)$, and forgetting everything but FV$(\varphi)$-values. One would get the same list by first selecting $\hat s \in L_C$, forgetting everything but FV($\varphi$)-values to get $s$, and augmenting with an $x$-value that is $F$-admissible to some $L_C$-extension of $s$. Whence the claim, the existential case, and the lemma. 
\end{itemize}
\end{proof}

\begin{Cor}[Domain locality] \label{cor:domain_locality}
Let $\mathcal A$ be a structure, $\mathcal L, \mathcal L'$ polyteams on $\mathcal A$, $\varphi$ a formula, and $\bar v$ such that FV$(\varphi) \subseteq \bar v \subseteq \dom(\mathcal L) \cap \dom(\mathcal L')$. If $\mathcal L \restriction \bar v = \mathcal L' \restriction \bar v$, then $\mathcal A \models_{\mathcal L}\ \varphi \iff \mathcal A \models_{\mathcal L'}\ \varphi$. In words, truth only depends on relevant information. 
\end{Cor}

\begin{Ex}
Reuse the set up of previous example, and consider the following polyteam formulas.

\begin{minipage}{0.35\textwidth}
\centering
$\begin{array}{cc}
x^i & y^i \\ \hline
0 & 0 \\ 
0 & 1 \\ 
2 & 0 
\end{array}$
\hspace{10pt}
$\begin{array}{cc}
a^j & b^k \\ \hline
0 & 0 \\ 
0 & 1 \\ 
1 & 1 
\end{array}$
\hspace{10pt}
$\begin{array}{c}
u^k \\ \hline
1 \\ 
2
\end{array}$
\end{minipage}
\begin{minipage}{0.63\textwidth}
\centering 
\begin{align*}
\alpha & :\ \neg R(a^j,b^k) \\ 
\beta & :\ \dep(x^i,y^i / a^j, b^k) \\ 
\varphi & :\ \alpha \vee \beta \\ 
\psi & :\ \exists v^{\chi(u),k}\ \forall w^{\chi(u),k}\ y / v \perp_{x, a / u} b / w
\end{align*}
\end{minipage}
\begin{itemize}
\item $\alpha$ is false because of the assignment $a^jb^k \mapsto (1,1)$. 
\item $\beta$ is false (for example) because of the assignments $x^iy^i \mapsto (0,0)$ and $a^jb^k \mapsto (0,1)$. 
\item $\varphi$ is true, for instance splitting as $\{\emptyset, \{a^jb^k \mapsto (0,0), a^jb^k \mapsto (0,1)\},\emptyset\}$ and its complement. 
\item $\psi$ is false, because $x^i$ and $a^j$ share the common value 0 that $u^k$ does not have. 
\end{itemize}
\end{Ex}

\begin{Fac} \label{fac:general_polyteam_properties}
\begin{itemize} 
\item Any team $X$ may be seen as the single-labeled polyteam $\{X\}$. Truth of team logic formulas with respect to $X$ can be reduced to truth of polyteam logic formulas with respect to $\{X\}$ : a dependence atom $\dep(\bar x, \bar y)$ is turned into $\dep(\bar x, \bar y / \bar x, \bar y)$, an independence atom $\bar x \perp_{\bar z} \bar y$ is turned into $\bar x / \bar x \perp_{\bar z, \bar z / \bar z} \bar y / \bar y$, and the rest is left unchanged. 
\item Distributivity of $\wedge$ and $\vee$ on one another fails in polyteams, as it does for teams. 
\end{itemize}
\end{Fac}

We say $\mathcal L'$ is a subpolyteam of $\mathcal L$, denoted $\mathcal L' \subseteq \mathcal L$, whenever $L'_C \subseteq L_C$ for each $C$. Besides, the union of $\mathcal L$ and $\mathcal L'$ is defined as $(\mathcal L \cup \mathcal L')_C := L_C \cup L'_C$. 

\begin{Def}
Let $\varphi$ be a formula in any polyteam logic. $\varphi$ is said to be \emph{closed downwards} ($\downarrow$-closed) if whenever $\mathcal A \models_{\mathcal L} \varphi$ and $\mathcal L' \subseteq \mathcal L$ we have $\mathcal A \models_{\mathcal L'} \varphi$. $\varphi$ is said to be \emph{closed under unions} ($\cup$-closed) if each time we have $I \neq \emptyset$, $(\mathcal L_i)_{i \in I}$ and $\forall i \in I, \mathcal A \models_{\mathcal L_i} \varphi$, then we also have $\mathcal A \models_{\cup_i \mathcal L_i} \varphi$.
\end{Def}

\begin{Pro} \label{pro:polyteam_logics_closures}
$\textnormal{FO}(=)$ is $\downarrow$-closed, $\textnormal{FO}(\subseteq)$ is $\cup$-closed, and $\textnormal{FO}(\perp)$ is neither.
\end{Pro}

\begin{proof}
Easy inductions, and counter examples. 
\end{proof}

\emph{Strict polyteam semantics} is what comes out of Definition \ref{def:polyteam_semantics} when we further impose that splits are disjoint (for $\vee$) and that supplementing functions map each assignment to a singleton (for $\exists$). It is important to note that, in general, strict semantics do not coincide with regular semantics, in particular domain locality fails in strict semantics. Nevertheless, it is easily seen that strict and regular semantics coincide on downward-closed formulas. In this case, it is often more convenient to think "in strict mode". 

\begin{Le}
Quantification on different labels commute. Moreover : 
$$\exists x^C\ \forall y^C\ \varphi(x,y,\bar v) \ \equiv\ \forall y^C\ \exists x^C\ \dep(\bar u, x / \bar u, x) \wedge \varphi(x,y,\bar v)$$
where $\varphi$ is any formula in FO($=$) with free variables $xy\bar v$, and $\bar u := \{v_0 \in \bar v : \chi(v_0) = C\}$.
\end{Le}

\begin{proof}
\begin{itemize}
\item[$\models$.] Assume the left hand side formula is true in some $\mathcal A$ and $\mathcal L$. Since polydependence logic is $\downarrow$-closed, the existential $\exists x^C$ is witnessed through a singleton-only supplementing function $F$, \textit{i.e.} such that $\mathcal L(F/x) \models \dep(\bar u, x / \bar u, x)$ holds. It is easy to check that $\mathcal L(A/y)(F'/x)$ witnesses the right hand side, where $F'(s) = F(s \restriction \bar u)$. 
\item[$\delsmo$.] Assume the right hand side holds, witnessed by $\mathcal L(A/y)(F/x)$. For each assignment $s'$ on domain $\bar u$, let $F'(s') := F(s)$ for any extension $s$ of $s'$, which is well defined thanks to the polydependence atom. It is easy to check that $\mathcal L(F'/x)(A/y)$ witnesses the left hand side. 
\end{itemize}
\end{proof}

\begin{Pro}
Every formula $\varphi$ in \textnormal{FO}(=) is equivalent to one of the form $\forall \bar x\ \exists \bar y\ \theta$ where $\theta$ is quantifier free.
\end{Pro}

\begin{proof}
First, any $\varphi$ can be transformed into an equivalent formula in quantifier-prenex normal form, in a similar way as in \cite{hannula_logics_2026} Proposition 2.4.16. Then, quantifiers are easily re-ordered inductively, using previous Lemma. 
\end{proof}

\subsection{Polyteam representation of physical experiments, and their sheaf-theoretic structure} \label{sec:polyteam_representation_experiments}

Let us fix a finite set of measurements $M := \{m_1, \cdots, m_n\}$ and a set of outcomes $O$. It does not matter what the elements of $M$ and $O$ actually are. We just need to agree on a correspondence between the $m_i$'s and which practical measurement they refer to, and between the $o$'s and which instrument-outputs they refer to. 

\begin{Def}
Any set $\mathcal M \subseteq \mathcal P(M)$ satisfying the following is called a \emph{maximal context set}. 
\begin{enumerate}
\item[(cover)] $\cup \mathcal M = M$
\item[(antichain)] There are no $C \subsetneq C'$ in $\mathcal M$.
\end{enumerate}
$\mathcal M$ is intended to describe the set of measurements from $M$ that can be performed together. The covering condition says that all measurements in $M$ are possible, and the antichain condition roughly says we always measure as many things as possible. 
\end{Def}

When one actually has a physical system in mind, there is one "real $\mathcal M$", and we may compute it. Indeed, the system state can be viewed as a density operator $\rho$ in a Hilbert space $H$, and each $m_i$ as a bounded self-adjoint operator $M_i$ on $H$. Then compatibility is exactly commutativity, and : 
\begin{center}
$\mathcal M := \{C \subseteq \{M_1, \cdots, M_n\} :$ the elements of $C$ pairwise commute and $C$ is maximal such$\}$
\end{center}
For this reason, from now on, let us consider $\mathcal M$ as fixed. 

\begin{Def} \label{def:sigma}
Let us define the labeled language $\sigma$ as follows. 
\begin{itemize}
\item $\sigma$ has one sort : the "measurement" sort. 
\item Each $C \in \mathcal M$ is a label in $\sigma$, and there is no other label. 
\item For each label $C$, there is one \emph{measurement variable} $m^C$ per $m \in C$.
\item Add whatever constants or relations that are relevant to the studied system (typically none). 
\end{itemize}
\end{Def}

\begin{Def} \label{def:tau}
Let us define the labeled language $\tau$ as follows. 
\begin{itemize}
\item $\tau$ has two sorts : the "measurement" sort and the "hidden" sort. 
\item Each $C \in \mathcal M$ is a label in $\tau$, and there is no other label. 
\item For each label $C$, there is one \emph{measurement variable} $m^C$ per $m \in C$, and one \emph{hidden variable} $z^C$.
\item Add whatever constants or relations that are relevant to the studied system (typically none). 
\end{itemize}
\end{Def}

\begin{Not}
If $V$ is a consistent set of measurements $\subseteq C$, we systematically denote by $\bar v^C$ (or simply $\bar v$) for $\{m^C,\ m \in V\}$. In particular $\bar c$ unambiguously means $\bar c^C$, thanks to the antichain condition. We always omit sorts by the following convention : a variable has the hidden sort if and only if it uses the letter $z$.
\end{Not}

\begin{Def}
Let $U \subseteq M$. An \emph{event} on $U$ is a map $s : U \to O$. Their collection is denoted by $\mathcal E(U)$. When considering a fixed $C \supseteq U$, $\mathcal E(U)$ may be identified with the set of assignments from $\bar u^C$ to $O$. 
\end{Def}

\begin{Pro} \label{pro:sheaf_events}
$\mathcal E$ is a \emph{sheaf}, \textit{i.e.} the following is true. 
\begin{itemize}
\item \emph{Restriction maps.} There are maps $\textnormal{res}^{U'}_U : \mathcal E(U') \to \mathcal E(U)$ ($U \subseteq U' \subseteq M$) such that : 
$$\forall U \subseteq U' \subseteq U'' \subseteq M,\ \textnormal{res}^{U'}_U \circ \textnormal{res}^{U''}_{U'} = \textnormal{res}^{U''}_U$$
\item \emph{Gluing property.} Assume $\bigcup_{n \in I} U_n = U \subseteq M$ and we are given $s_n \in \mathcal E(U_n)$ for each $n \in I$ such that : 
$$\forall n, p \in I,\ \textnormal{res}^{U_n}_{U_n \cap U_p}(s_n) = \textnormal{res}^{U_p}_{U_n \cap U_p}(s_p)$$
Then there is $s \in \mathcal E(U)$ that restricts to $s_n$ via $\textnormal{res}^U_{U_n}$ for each $n \in I$. 
\end{itemize}
In words, we have a way to test the consistency of some local information (restriction maps), and as soon as the test is passed, we can merge it into global information (gluing). 
\end{Pro}

\begin{figure}
\centering
\begin{tikzcd}
	{\mathcal E(U'')} && {\mathcal E(U')} && {\mathcal E(U)}
	\arrow["{\text{res}^{U''}_{U'}}", from=1-1, to=1-3]
	\arrow["{\text{res}^{U''}_U}"', curve={height=12pt}, from=1-1, to=1-5]
	\arrow["{\text{res}^{U'}_U}", from=1-3, to=1-5]
\end{tikzcd}
\hspace{20pt}
\raisebox{-35pt}{
\begin{tikzpicture}[scale=1]
\draw (0,0) -- (1,0); 
\node at (0.5,-0.2) {$s_1$}; 
\draw (0.6,0.5) -- (2,0.5); 
\node at (1.3,0.7) {$s_2$}; 
\draw[dotted] (0.6,0) rectangle (1,0.5);
\node at (0.8,0.25) {$\verteq$}; 
\draw (1.5,0) -- (2.2,0); 
\node at (1.85,-0.2) {$s_3$};
\draw[dotted] (1.5,0) rectangle (2,0.5); 
\node at (1.75,0.25) {$\verteq$};
\draw (2.2,0.5) -- (2.5,0.5); 
\node at (2.35,0.7) {$s_4$}; 
\draw (0,-1) -- (2.5,-1);  
\node at (1.25,-1.2) {$s$}; 
\draw[dotted] (0,0) -- (0,-1); 
\draw[dotted] (2.5,0.5) -- (2.5,-1); 
\draw[->] (2.6,0.5) -- (2.6,-1);
\end{tikzpicture}}
\end{figure}

\begin{proof}
Restriction maps are obviously given by domain-restriction of maps : $res^{U'}_U = \cdot \restriction U$. The hypothesis of the gluing property exactly says that $s := \cup_n s_n$ is a well defined map, and that map is of course suitable.  
\end{proof}

In general sheaf theoretic terms, the $s_n$'s are called (local) \emph{sections}, and $s$ is called a \emph{global section}. 

\begin{Not}
For further-coming reasons, use the notation $\mathcal D_{\mathbb B}\mathcal E(U)$ for the set $\mathcal P^+(\mathcal E(U))$ of non-empty subsets of events on $U$. When considering a fixed $C \supseteq U$, elements of $\mathcal D_{\mathbb B}\mathcal E(U)$ may be identified with non-empty teams on $O$ with domain $\bar u^C$. 
\end{Not}

\begin{Def}
An \emph{empirical model} (or \emph{empirical polyteam}) is a polyteam $\mathcal L = \{L_C, C \in \mathcal M\}$ such that $L_C \in \mathcal D_{\mathbb B}\mathcal E(C)$ for each $C$. Intuitively, each component team $L_C$ lists all possible outputs one may get while measuring $C$ (there is at least one). 
\end{Def}

\begin{Fac}
$\mathcal D_{\mathbb B}\mathcal E$ is a presheaf, \textit{i.e.} it has restriction maps. 
\end{Fac}

\begin{proof}
Define res$^{U'}_U(X) := X \restriction U$ for each $\emptyset \subsetneq X \subseteq \mathcal E(U')$. Now, given $Y \in \mathcal D_{\mathbb B}\mathcal E(U'')$ and $s \in \mathcal E(U)$ : 
\begin{multline*}
s \in \text{res}^{U'}_U \circ \text{res}^{U''}_{U'}(Y)
\iff \exists s' \in \text{res}^{U''}_{U'}(Y), s = s' \restriction U
\iff \exists s'' \in Y, s = (s'' \restriction U') \restriction U \\
\iff \exists s'' \in Y, s = s'' \restriction U
\iff s \in \text{res}^{U''}_U(Y)
\end{multline*}
\end{proof}

\begin{Rem}
We have already used this result multiple times, for instance in the $\wedge \vee$ cases of Lemma \ref{le:locality}. 
\end{Rem}

\begin{Pro} \label{pro:not_sheaf}
$\mathcal D_{\mathbb B}\mathcal E$ is not a sheaf, as witnessed by the Hardy model (Figure \ref{fig:hardy}).
\end{Pro}

\begin{proof}
Assume by contradiction the Hardy model $\mathcal L$ has a global section $G$ (team with domain $aba'b'$). Since $ab \mapsto (0,0)$ occurs in $L_{\{a,b\}}$, and $G \restriction \{a,b\} = L_{\{a,b\}}$ one of the following global assignments must witness it in $G$. 
\begin{align*}
s_1 : aba'b' & \mapsto (0,0,0,0) \\ 
s_2 : aba'b' & \mapsto (0,0,0,1) \\ 
s_3 : aba'b' & \mapsto (0,0,1,0) \\ 
s_4 : aba'b' & \mapsto (0,0,1,1)
\end{align*}
However, $s_1$ cannot be in $G$ since otherwise $(s_1 \restriction \{a',b\}) : a'b \mapsto (0,0)$ would be in $L_{\{a',b\}}$. Similarly $s_2 \notin G$ (see $L_{\{a',b\}}$), $s_3 \notin G$ (see $L_{\{a,b'\}}$), and $s_4 \notin G$ (see $L_{\{a',b'\}}$). A contradiction. 
\end{proof}

\begin{figure}
\centering
$\begin{array}{cc}
a & b \\ \hline
0 & 0 \\ 
0 & 1 \\ 
1 & 0 \\ 
1 & 1
\end{array}$
\hspace{20pt}
$\begin{array}{cc}
a' & b \\ \hline
0 & 1 \\ 
1 & 0 \\ 
1 & 1
\end{array}$
\hspace{20pt}
$\begin{array}{cc}
a & b' \\ \hline
0 & 1 \\ 
1 & 0 \\ 
1 & 1
\end{array}$
\hspace{20pt}
$\begin{array}{cc}
a' & b' \\ \hline
0 & 0 \\ 
0 & 1 \\
1 & 0 
\end{array}$
\caption{The Hardy model}
\label{fig:hardy}
\end{figure}

\begin{Not}
We have omitted labels in Figure \ref{fig:hardy}, because they may be unambiguously deduced. For example, the first $a$ should be labeled $\{a,b\}$ whereas the second $a$ should be labeled $\{a,b'\}$. 
\end{Not}

\begin{Rem} \label{rem:hardy_models}
In other words, Proposition \ref{pro:not_sheaf} says that the polyteam defining the Hardy model is some local information that undoubtedly does \emph{not} come from global information on the variables $\{a,a',b,b'\}$, \textit{i.e.} that there is no list of possible values for the tuple $aa'bb'$ that yields the four mentioned lists respectively when truncating its domain. \emph{There is nothing to be amazed of}, this is just combinatorial play. Indeed, the proof of Proposition \ref{pro:not_sheaf} would apply to any polyteam satisfying the following.
\begin{center}
\begin{tabular}{l}
$\exists C \supseteq \{a,b\}, ab \mapsto (0,0)$ is (the restriction of some assignment) in $L_C$ \\ 
$\exists C \supseteq \{a',b\}, a'b \mapsto (0,0)$ is not (the restriction of some assignment) in $L_C$ \\
$\exists C \supseteq \{a,b'\}, ab' \mapsto (0,0)$ is not (the restriction of some assignment) in $L_C$ \\
$\exists C \supseteq \{a',b'\}, ab \mapsto (1,1)$ is not (the restriction of some assignment) in $L_C$ \\
\end{tabular}
\end{center}
What would be interesting though, is to know that this counter example polyteam describes a real quantum system. The latter system would then exemplify contextuality, in the form known as \emph{non-locality}. And indeed, the Hardy model is such, as we shall see in Section \ref{sec:theorems}. 
\end{Rem}

Since the foundational work \cite{einstein_can_1935} of Einstein, Podolsky and Rosen in 1935, the possibility that the probabilistic nature of quantum measurements is in fact an illusion has been studied. The idea consist in saying that maybe quantum systems have attributes that we do not (or cannot) know about, which actually determine the values of the measurements which seem probabilistic to the human level. For example, maybe all electrons are either of type $\Delta$, $\triangleleft$ or $\nabla$, all those that have type $\Delta$ or $\triangleleft$ have spin $\uparrow$, all the remaining have spin $\downarrow$, and actually none of them are in a so called \emph{superposed state}. Such a description of a system is called a \emph{hidden-variable model}. These hidden attributes are exactly what our hidden variables $z^C$ are intended to account for. They are given a different sort specifically because we do not know if they can be measured, and what kind of values they may take. The quantifier $\texists$ has been introduced to guess that set of values, which $\exists$ cannot do. 

\begin{Not}
Let $z$ represent a hidden property (whatever mathematical object it really is), and $\Lambda$ be the set of values $z$ can take. 
\end{Not}

\begin{Def}
Let $U \subseteq M$. A \emph{hidden event} on $U$ is a map $s : U \cup \{z\} \to O \cup \Lambda$ such that $s[U] \subseteq O$ and $s(z) \in \Lambda$. Their collection is denoted $\mathcal E^h(U)$. When considering a fixed $C \supseteq U$, $\mathcal E^h(U)$ may be identified with the set of assignments from $\bar u^Cz^C$ to $O^{|\bar u|} \times \Lambda$. 
\end{Def}

\begin{Pro}
$\mathcal E^h$ is a sheaf. 
\end{Pro}

\begin{proof}
It holds for the same reason as for $\mathcal E$ (see Proposition \ref{pro:sheaf_events}). 
\end{proof}

\begin{Not}
For further-coming reasons, use the notation $\mathcal D_{\mathbb B}\mathcal E^h(U)$ for the set $\mathcal P^+(\mathcal E^h(U))$ of non-empty subsets of hidden events on $U$. When considering a fixed $C \supseteq U$, elements of $\mathcal D_{\mathbb B}\mathcal E^h(U)$ may be identified with non-empty teams with domain $\bar u^Cz^C$ on a two-sorted structure with domains $O$ and $\Lambda$. 
\end{Not}

\begin{Def}
A \emph{hidden variable model} (or \emph{hidden variable polyteam}) is a polyteam $\mathcal H = \{H_C, C \in \mathcal M\}$ such that $H_C \in \mathcal D_{\mathbb B}\mathcal E^h(C)$ for each $C$. We may think of each $H_C$ as being the succession of all possible outputs one may get while measuring $C$, for each value of the hidden variable. 
\end{Def}

\begin{Fac}
$\mathcal D_{\mathbb B}\mathcal E^h$ is a presheaf. 
\end{Fac}

\begin{proof}
Define $\text{res}^{U'}_U(Y) := Y \restriction U \cup \{z\}$ for each $\emptyset \subsetneq Y \subset \mathcal E^h(U)$, and conclude as for $\mathcal D_{\mathbb B}\mathcal E$. 
\end{proof}

\begin{Def}
Let $\mathcal L$ be an empirical model, and $\mathcal H$ be a hidden variable model. We say that $\mathcal H$ is a \emph{realization} of $\mathcal L$ whenever restricting the domain of $\mathcal H$ to measurement variables (wiping out the $z^C$'s) yields $\mathcal L$. 
\end{Def}

\begin{Le} 
Let $\varphi$ be any polyteam $\tau$-formula, $\mathcal L$ an empirical polyteam, and enumerate $\mathcal M$ as $\{C_1, \cdots, C_n\}$.
\begin{center}
$\mathcal L$ has a hidden variable realization $\mathcal H$ satisfying $\varphi \iff \mathcal L \models \texists z^{C_1} \exists z^{C_2} \cdots \exists z^{C_n} \varphi$
\end{center}
\end{Le}

\begin{Not}
We allow ourselves to write $\texists z$ as a shorthand for $\texists z^{C_1} \exists z^{C_2} \cdots \exists z^{C_n}$. 
\end{Not}

\begin{Rem} \label{rem:single_hv}
If we had given ourselves countably many $z$-variables ($z_1, z_2, \cdots$) for each $C$, it would have been possible to prove that for any $\varphi(z_1)$ : 
$$\texists z_1 \varphi \equiv \overline{\bigvee_{\ell \in \mathbb N^*}} \texists z_1 \exists z_2 \cdots \exists z_\ell\ \varphi(z_1\cdots z_\ell/z_1)$$
where $\varphi(z_1\cdots z_\ell/z_1)$ is the formula obtained from $\varphi$ by replacing each occurrence of $z_1$ by the tuple $z_1 \cdots z_\ell$, and $\bar \vee$ stands for the intuitionistic disjunction ($\mathcal L \models \psi\ \bar \vee\ \theta \iffdef \mathcal L \models \psi$ or $\mathcal L \models \theta$). 
Indeed, $\models$ is obvious. Conversely, if we have $\mathcal L \models \texists z_1 \exists z_2 \cdots \exists z_\ell\ \varphi(z_1 \cdots z_\ell/z_1)$ for some $\ell$, get supplementing functions $F_1, \cdots, F_\ell$ and transform $\mathcal L(F_1/z_1) \cdots (F_\ell/z_\ell)$ into $\mathcal L(F/z_1)$ as follows in each component team $L_C$. 

\begin{center}
\begin{tabular}{ccccc}
$m_1^C$ & $m_2^C$ & $\cdots$ & $z_1^C$ & $z_2^C$ \\ \hline
0 & 0 & $\cdots$ & $-1$ & 1 \\ 
0 & $\pi/2$ & $\cdots$ & $-1$ & 0 \\ 
$-\pi/2$ & 0 & $\cdots$ & 0 & 1\\
$-\pi/2$ & $\pi/2$ & $\cdots$ & 1 & 1
\end{tabular}
\hspace{20pt}
$\longmapsto$
\hspace{20pt}
\begin{tabular}{cccc}
$m_1^C$ & $m_2^C$ & $\cdots$ & $z_1^C$ \\ \hline
0 & 0 & $\cdots$ & $(-1,1)$ \\ 
0 & $\pi/2$ & $\cdots$ & $(-1,0)$ \\ 
$-\pi/2$ & 0 & $\cdots$ & $(0,1)$ \\
$-\pi/2$ & $\pi/2$ & $\cdots$ & $(1,1)$
\end{tabular}
\end{center}
As a result, it is legitimate to think of $\texists z$ as meaning "There is a certain number of hidden \emph{elementary} physical attributes such that ..." : elementariness can be achieved by splitting $z_1^C$ into more understandable bits. 
\end{Rem}

\begin{Def}
Recall that the maximal context set $\mathcal M$ is fixed, and let $G \in \mathcal D_{\mathbb B}\mathcal E(M)$ be an empirical global section. The \emph{contextualization} of $G$ is the empirical polyteam $\{\text{res}^M_C(G), C \in \mathcal M\}$ obtained through projecting $G$ on each $C \in \mathcal M$. Similarly for hidden variable global sections $H \in \mathcal D_{\mathbb B}\mathcal E^h(M)$. 
\end{Def}

\begin{Rem}
Proposition \ref{pro:not_sheaf} can be rephrased as "the contextualization map is not surjective".
\end{Rem}

\begin{Le} \label{le:polyteams_commutative_diagram}
The following diagram commutes. 

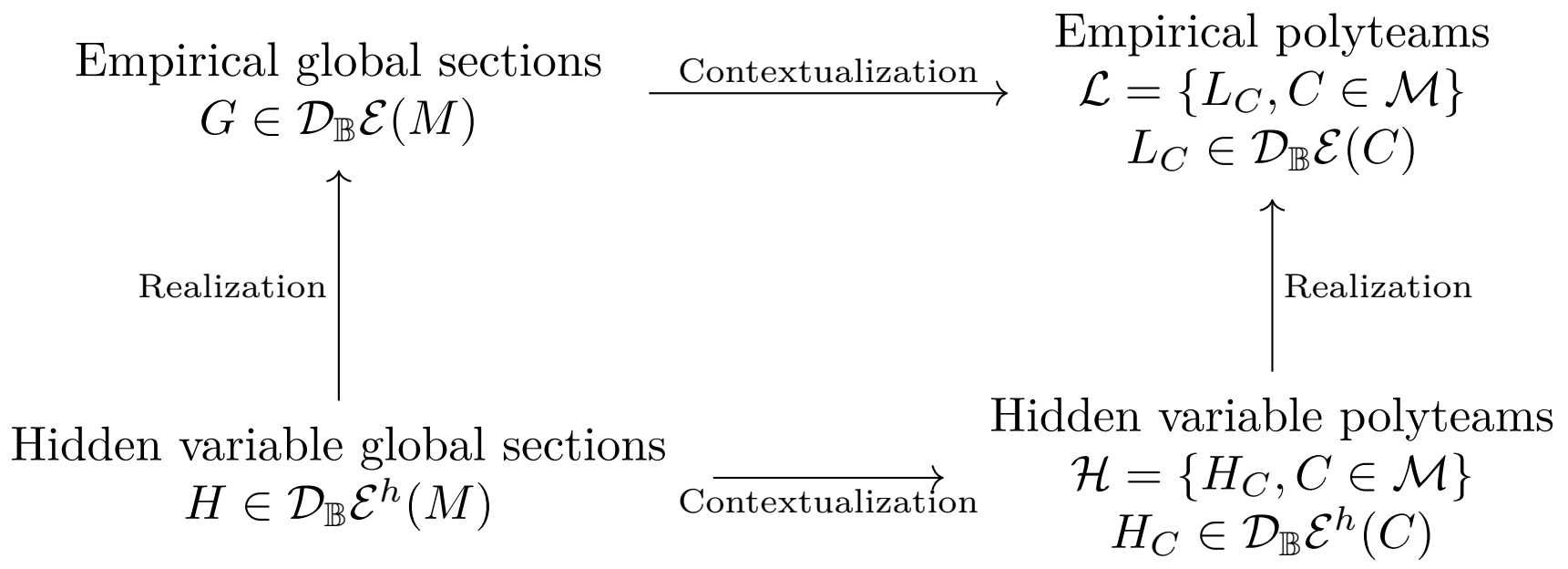
\[\begin{tikzcd}
	\begin{array}{c} \text{Empirical global sections} \\ G \in \mathcal D_{\mathbb B}\mathcal E(M) \end{array} && \begin{array}{c} \text{Empirical polyteams} \\ \mathcal L = \{L_C, C \in \mathcal M\} \\ L_C \in \mathcal D_{\mathbb B}\mathcal E(C) \end{array} \\
	\\
	\begin{array}{c} \text{Hidden variable global sections} \\ H \in \mathcal D_{\mathbb B}\mathcal E^h(M) \end{array} && \begin{array}{c} \text{Hidden variable polyteams} \\ \mathcal H = \{H_C, C \in \mathcal M\} \\ H_C \in \mathcal D_{\mathbb B}\mathcal E^h(C) \end{array}
	\arrow["{\text{Contextualization}}", from=1-1, to=1-3]
	\arrow["{\text{Realization}}", from=3-1, to=1-1]
	\arrow["{\text{Contextualization}}"', from=3-1, to=3-3]
	\arrow["{\text{Realization}}"', from=3-3, to=1-3]
\end{tikzcd}\]
\end{Le}

\begin{proof}
Easy verifications.
\end{proof}

\begin{Pro}
$\mathcal D_{\mathbb B}\mathcal E^h$ is not a sheaf. 
\end{Pro}

\begin{proof}
Assume by contradiction it is, and let $\mathcal L$ denote the Hardy model. Let $\mathcal H$ be an arbitrary hidden variable polyteam $\mathcal H$ that realizes $\mathcal L$ : it exists since it suffices to throw any values under each $z^C$ column. By assumption, there is a hidden variable global section for $\mathcal H$, call it $H$. Let $G$ be the unique empirical global section realized by $H$. By Lemma \ref{le:polyteams_commutative_diagram}, $G$ contextualizes as $\mathcal L$, thus contradicting Proposition \ref{pro:not_sheaf}. 
\end{proof}

The previous sheaf theoretic framework for contextuality is inspired from, and isomorphic to that of \cite{abramsky_sheaf-theoretic_2011}. More precisely, \cite{abramsky_sheaf-theoretic_2011} has its own $\mathcal E$ and $\mathcal E^h$, its own commuting square, there are explicit bijections connecting the two squares, and the resulting cube commutes. The essential difference between these two frameworks is that ours has been intentionally fine-tuned to support logical analysis. 

To end this section, let's exemplify the new notions with a model we shall deal with later. 

\begin{Ex}[Possibilistic collapse of Bell's model]
Consider the empirical polyteam $\mathcal L$ from Figure \ref{fig:bell_model_possibilistic}. It is build from the measurement set $M = \{a,a',b,b'\}$, the outcome set $O = \{0,1\}$, and the maximal context set $\mathcal M = \{\{a,b\}, \{a',b\}, \{a,b'\}, \{a',b'\}\}$. This empirical model corresponds to the repetition of a physical experiment, similar to that performed by Bell in 1964, as reported in \cite{bell_einstein_1964}. 

\begin{quote}
Consider a pair of spin one-half particles formed somehow in the singlet spin state and moving freely in opposite directions. Measurements can be made, say by Stern-Gerlach magnets, on selected components of the spins $\vec \sigma_1$ and $\vec \sigma_2$. If measurement of the component $\vec \sigma_1 \cdot \vec a$ where $\vec a$ is some unit vector, yields the value $+1$ then, according to quantum mechanics, measurement of $\vec \sigma_2 \cdot \vec a$ must yield the value $-1$ and vice versa. (\textit{cf.} \cite{bell_einstein_1964} p.195)
\end{quote}
Now assume Alice and Bob are spatially separated (unable to communicate). Alice can choose to perform either measurement $\vec \sigma_1 \cdot \vec a$ (variable $a$) or measurement $\vec \sigma_1 \cdot \vec a'$ (variable $a'$) on the first particle. Bob may perform the same measurements on the second particle (variables $b$ and $b'$). Then, if one builds the table of all possible joint observations, one will get $\mathcal L$ (using 0 and 1 instead of $\pm 1$). These observations have a global explanation, \textit{i.e.} the empirical global section $G$ does contextualize to $\mathcal L$. 

\begin{figure}
\centering
$\mathcal L =\hspace{15pt} \begin{array}{cc}
a & b \\ \hline
0 & 0 \\ 
1 & 1
\end{array}$
\hspace{20pt}
$\begin{array}{cc}
a' & b \\ \hline
0 & 0 \\
0 & 1 \\ 
1 & 0 \\ 
1 & 1
\end{array}$
\hspace{20pt}
$\begin{array}{cc}
a & b' \\ \hline
0 & 0 \\ 
0 & 1 \\ 
1 & 0 \\ 
1 & 1
\end{array}$
\hspace{20pt}
$\begin{array}{cc}
a' & b' \\ \hline
0 & 0 \\ 
0 & 1 \\ 
1 & 0 \\ 
1 & 1
\end{array}$
\caption{The possibilistic collapse of Bell's model}
\label{fig:bell_model_possibilistic}
\end{figure}

\begin{center}
$G =\hspace{15pt}$
\begin{tabular}{cccc}
$a$ & $a'$ & $b$ & $b'$ \\ \hline
0 & 0 & 0 & 0 \\ 
0 & 0 & 0 & 1 \\ 
0 & 1 & 0 & 0 \\ 
0 & 1 & 0 & 1 \\ 
1 & 0 & 1 & 0 \\ 
1 & 0 & 1 & 1 \\ 
1 & 1 & 1 & 0 \\ 
1 & 1 & 1 & 1
\end{tabular}
\end{center}
It would be even better, in line with what Einstein, Podolsky and Rosen expect in \cite{einstein_can_1935}, to have a deterministic hidden variable global model of the system, such as the following.

\begin{center}
$H =\hspace{15pt}$
\begin{tabular}{ccccc}
$a$ & $a'$ & $b$ & $b'$ & $z$ \\ \hline
0 & 0 & 0 & 0 & $\triangleleft \uparrow$ \\ 
0 & 0 & 0 & 1 & $\triangleleft \rightarrow$ \\ 
0 & 1 & 0 & 0 & $\triangleleft \downarrow$ \\ 
0 & 1 & 0 & 1 & $\triangleleft \leftarrow$ \\ 
1 & 0 & 1 & 0 & $\triangleright \uparrow$ \\ 
1 & 0 & 1 & 1 & $\triangleright \rightarrow$ \\ 
1 & 1 & 1 & 0 & $\triangleright \downarrow$ \\ 
1 & 1 & 1 & 1 & $\triangleright \leftarrow$
\end{tabular}
\end{center}
$H$ is deterministic in the following sense : $H \models\ \dep(z,aa'bb')$ (details in Section \ref{sec:theorems}). $H$ contextualizes as $\mathcal H$, and the following square commutes. 
\begin{figure}[H] \label{fig:ex1}
\centering
$\hspace{-50pt}\mathcal H =\hspace{10pt}$ \begin{tabular}{ccc}
$a$ & $b$ & $z^{ab}$ \\ \hline
0 & 0 & $\triangleleft \uparrow$ \\ 
0 & 0 & $\triangleleft \rightarrow$ \\ 
0 & 0 & $\triangleleft \downarrow$ \\ 
0 & 0 & $\triangleleft \leftarrow$ \\ 
1 & 1 & $\triangleright \uparrow$ \\ 
1 & 1 & $\triangleright \rightarrow$ \\ 
1 & 1 & $\triangleright \downarrow$ \\ 
1 & 1 & $\triangleright \leftarrow$
\end{tabular}
\hspace{15pt}
\begin{tabular}{ccc}
$a'$ & $b$ & $z^{a'b}$ \\ \hline
0 & 0 & $\triangleleft \uparrow$ \\ 
0 & 0 & $\triangleleft \rightarrow$ \\ 
1 & 0 & $\triangleleft \downarrow$ \\ 
1 & 0 & $\triangleleft \leftarrow$ \\ 
0 & 1 & $\triangleright \uparrow$ \\ 
0 & 1 & $\triangleright \rightarrow$ \\ 
1 & 1 & $\triangleright \downarrow$ \\ 
1 & 1 & $\triangleright \leftarrow$
\end{tabular}
\hspace{15pt}
\begin{tabular}{ccc}
$a$ & $b'$ & $z^{ab'}$ \\ \hline
0 & 0 & $\triangleleft \uparrow$ \\ 
0 & 1 & $\triangleleft \rightarrow$ \\ 
0 & 0 & $\triangleleft \downarrow$ \\ 
0 & 1 & $\triangleleft \leftarrow$ \\ 
1 & 0 & $\triangleright \uparrow$ \\ 
1 & 1 & $\triangleright \rightarrow$ \\ 
1 & 0 & $\triangleright \downarrow$ \\ 
1 & 1 & $\triangleright \leftarrow$
\end{tabular}
\hspace{15pt}
\begin{tabular}{ccc}
$a'$ & $b'$ & $z^{a'b'}$ \\ \hline
0 & 0 & $\triangleleft \uparrow$ \\ 
0 & 1 & $\triangleleft \rightarrow$ \\ 
1 & 0 & $\triangleleft \downarrow$ \\ 
1 & 1 & $\triangleleft \leftarrow$ \\ 
0 & 0 & $\triangleright \uparrow$ \\ 
0 & 1 & $\triangleright \rightarrow$ \\ 
1 & 0 & $\triangleright \downarrow$ \\ 
1 & 1 & $\triangleright \leftarrow$
\end{tabular}
\hspace{15pt}
$\begin{tikzcd}
	G && \mathcal L \\
	\\
	H && \mathcal H
	\arrow["{\text{Cont.}}", from=1-1, to=1-3]
	\arrow["{\text{Real.}}", from=3-1, to=1-1]
	\arrow["{\text{Cont.}}"', from=3-1, to=3-3]
	\arrow["{\text{Real.}}"', from=3-3, to=1-3]
\end{tikzcd}$
\end{figure}
\end{Ex}

\newpage

\section{Logical formalization of physical properties and contextuality theorems} \label{sec:theorems}

\subsection{A first idea of contextuality theorems}

The original theorem of Bell (see \cite{bell_einstein_1964}) uses probabilistic models. For this reason, we deal with its formalization only in Section \ref{sec:semirings}. However, we may already try to formalize a strengthening of Bell's theorem without probabilities. The \emph{principle of locality} is the view according to which an object can only be influenced by its immediate surroundings. \emph{Local realism} is a slight strengthening of locality, which we shall define precisely in Section \ref{sec:locality_local_realism}.

\begin{The}[Bell, possibilistic version] \label{th:bell}
There is a quantum system that violates local realism. 
\end{The}

Let us now turn to the Kochen-Specker theorem, which states an even wildest form of contextuality. A quantum system is called \emph{strongly contextual} whenever for any global assignment $s$ of the relevant attributes, there is at least one measurement context $C$ for which it is impossible to observe $s(\bar c)$ together. Put differently, thinking that the quantum system might have all its attributes well-defined at any moment simply makes no sense. 

\begin{The}[Kochen-Specker] \label{th:kochen_specker}
There is a strongly contextual quantum system. 
\end{The}

\subsection{Logical formalization of physical properties} \label{sec:physical_properties}
There are many relevant properties of physical models (whether they are empirical or using hidden variables). Abramsky gives in \cite{abramsky_relational_2013} (2013) a precise definition of some of these properties, supported by what he calls \emph{relational models} (which are essentially databases). From now on, let us refer to this set of properties as $\mathcal P$. The purpose of this section is to express within our logics the properties from $\mathcal P$, supported by our empirical and hidden variable models. As we said before, the interest in doing so is (at least) twofold : being able to recognize contextuality phenomena in other topics that are suitably-represented by polyteams, and trying to axiomatize truth in order to benefit from computers when searching for contextuality. The results are grouped in Table \ref{tab:physical_properties}. 

\subsubsection{Weak determinism}
\begin{Def}
Let $\mathcal L$ be an empirical model. $\mathcal L$ is said to support \emph{Weak Determinism} whenever $L_C$ is a singleton for each $C \in \mathcal M$. In words, as soon as one chooses what to measure, one knows what will come out : the measurement context determines the output. 
\end{Def}

\begin{Le}
$\mathcal L$ supports Weak Determinism $\iff \mathcal L \models WD$, where $WD$ is the following formula. 
$$WD : \bigwedge_{m \in C \in \mathcal M} \dep(\emptyset^C, m^C / \emptyset^C, m^C) \ \overset{\text{not}}{\equiv} \bigwedge_{m \in C \in \mathcal M} \dep(m^C)\ \equiv \bigwedge_{C \in \mathcal M} \dep(\bar c)$$
\end{Le}

\begin{Rem}
The satisfaction $\mathcal L \models WD$ should be understood as $\mathcal A \models_{\mathcal L} WD$, where $\mathcal A$ is any $\sigma$-structure on domain $O$, and $\sigma$ is what remains from the definition of $\tau$ (Definition \ref{def:tau}) when the hidden sort and the hidden variables are wiped out. 
\end{Rem}

\begin{Ex}
Figure \ref{fig:weak_determinism} shows an empirical polyteam that supports Weak Determinism. 
\end{Ex}

\begin{figure}
\centering
$\begin{array}{cc}  
a & b \\ \hline
0 & -1
\end{array}$
\hspace{20pt}
$\begin{array}{ccc} 
b & c & d \\ \hline
1 & 0 & 2
\end{array}$
\hspace{20pt}
$\begin{array}{cc} 
d & e \\ \hline
1 & -1
\end{array}$
\caption{An empirical polyteam supporting Weak Determinism}
\label{fig:weak_determinism}
\end{figure}

\subsubsection{Strong determinism}
\begin{Def}
Let $\mathcal L$ be an empirical model. $\mathcal L$ is said to support \emph{Strong Determinism} whenever each $m \in M$ always get the same value, whatever the context $C$. In words, each physical attribute $m$ has a well-defined deterministic value $o$ that is observed in all the contexts $C \ni m$. 
\end{Def}

\begin{Le}
$\mathcal L$ supports Strong Determinism $\iff \mathcal L \models SD$, where $SD$ is the following formula. 
$$SD : \bigwedge_{m \in C,C' \in \mathcal M} \dep(\emptyset^C, m^C / \emptyset^{C'}, m^{C'}) \ \overset{\text{not}}{\equiv} \bigwedge_{m \in C,C' \in \mathcal M} \dep(m^C/m^{C'})$$
\end{Le}

\begin{Ex}
Figure \ref{fig:strong_determinism} shows an empirical polyteam that supports Strong Determinism. 
\end{Ex}

\begin{figure}
\centering
$\begin{array}{cc}  
a & b \\ \hline
0 & -1
\end{array}$
\hspace{20pt}
$\begin{array}{ccc} 
b & c & d \\ \hline
-1 & 0 & 2
\end{array}$
\hspace{20pt}
$\begin{array}{cc} 
d & e \\ \hline
2 & -1
\end{array}$
\caption{An empirical polyteam supporting Strong Determinism}
\label{fig:strong_determinism}
\end{figure}

\subsubsection{No-signalling}
This property is easier to understand in the special case of the so called Bell-type models.
Assume (as in \cite{bell_einstein_1964}) that the system is spatially separated, and that each experimenter is to perform a measurement on its respective part of the system. Saying the system is \emph{no-signalling} roughly means that Alice has no way of knowing which measurements the other experimenters have decided to perform for themselves. 

\begin{Ex}
Let $\mathcal L$ be the Hardy model, explicited in Figure \ref{fig:hardy}. Let us see that $\mathcal L$ supports no-signalling. Assume Bob repeats measurement $b$ and never chooses $b'$. Let us identify as Alice, and let us repeat measurement $a$ for some time. The measurement context is then $\{a,b\}$, and we (Alice) will observe the values 0 and 1 for our measurement $a$. Now let us repeat measurement $a'$ for some time. The measurement context is then $\{a',b\}$, and we (Alice) will observe the values 0 and 1 again. It is easy to check that precisely the same thing would have happened if Bob had chosen $b'$. In other words, Alice has no way of knowing what Bob does (and conversely) : there is no "signal". 
\end{Ex}

\begin{Fac}
\begin{tabular}[t]{rcl}
$\mathcal L \text{ supports No-Signalling}$ & $\iff$ & $\mathcal L \text{ is locally consistent}$ \\ 
& $\iff$ & $\mathcal L$ satisfies the premise of the gluing property \\
& $\iff$ & $\forall C, C' \in \mathcal M, L_C \restriction (C \cap C') = L_{C'} \restriction (C \cap C')$
\end{tabular}
\end{Fac}

\begin{Le} \label{le:ns_formula}
$\mathcal L$ supports No-Signalling $\iff \mathcal L \models NS$ where $NS$ is the following formula. 
$$NS : \bigwedge_{U \subseteq C, C' \in \mathcal M} \bar u^C / \bar u^{C'} \perp \emptyset^C / \emptyset^{C'}\ \equiv \bigwedge_{U \subseteq C,C' \in \mathcal M} \bar u^C \subseteq \bar u^{C'}$$
\end{Le}

\subsubsection{Hidden weak determinism}
\begin{Def}
Let $\mathcal H$ be a hidden variable model. $\mathcal H$ is said to support \emph{Hidden Weak Determinism} whenever in each context $C$, the value of the hidden variable determines the output. 
\end{Def}

\begin{Le}
$\mathcal H$ supports Hidden Weak Determinism $\iff \mathcal H \models HWD$ where $HWD$ is the following formula. 
$$HWD : \bigwedge_{C \in \mathcal M} \dep(z^C, \bar c / z^C, \bar c)\ \overset{\text{not}}{\equiv} \bigwedge_{C \in \mathcal M} \dep(z^C, \bar c) $$
\end{Le}

\begin{Rem}
The satisfaction $\mathcal H \models HWD$ should be understood as $\mathcal B \models_{\mathcal H} HWD$, where $\mathcal B$ is any $\tau$-structure with measurement-sort domain $O$ and hidden-sort domain $\Lambda$.
\end{Rem}

\subsubsection{Hidden strong determinism}
\begin{Def}
Let $\mathcal H$ be a hidden variable model. $\mathcal H$ is said to support \emph{Hidden Strong Determinism} whenever the output of a measurement $m \in M$ is fully determined by the value of the hidden variable. 
\end{Def}

\begin{Le}
$\mathcal H$ supports Hidden Strong Determinism $\iff \mathcal H \models HSD$ where $HSD$ is the following formula. 
$$HSD :\ \bigwedge_{m \in C, C' \in \mathcal M} \dep(z^C, m^C / z^{C'}, m^{C'})\ \equiv \bigwedge_{U \subseteq C, C' \mathcal M} \dep(z^C, \bar u^C / z^{C'}, \bar u^{C'})$$
\end{Le}

\subsubsection{Single-valuedness}
\begin{Def}
Let $\mathcal H$ be a hidden variable model. $\mathcal H$ is said to support \emph{Single-Valuedness} whenever the hidden variable is set to the same value in every assignment from any measurement context.  
\end{Def}

\begin{Le}
$\mathcal H$ supports Single-Valuedness $\iff \mathcal H \models SV$ where $SV$ is the following formula. 
$$SV :\ \bigwedge_{C,C' \in \mathcal M} \dep(\emptyset^C, z^C / \emptyset^{C'}, z^{C'}) \ \overset{\text{not}}{\equiv} \bigwedge_{C, C' \in \mathcal M} \dep(z^C/z^{C'})$$
\end{Le}

\subsubsection{Parameter Independence}
This property is the reformulation of no-signalling for hidden variable models. Assume again the studied system is Bell-type, \textit{i.e.} that multiple agents are assigned spatially separated bits of the system, and that each is asked to choose and perform a measurement from a previously fixed list. Saying this experiment supports \emph{Parameter Independence} means that Alice (for instance) cannot get any information on what the other experimenters are choosing to measure, \emph{for each value of the hidden variable}. Put differently, Alice may obtain information on the value of the hidden variable of each instance of the system (for instance, she may learn from her measurement that her particle has hidden attribute either $\triangleleft \uparrow$ or $\triangleleft \downarrow$), but that hidden value fully determines what she may possibly observe for each of her measurements, independently of what the others choose to measure, hence she may not empirically learn anything about the latter. In fact, Parameter Independence also says that, even if Alice and Bob were to gather to acquire information on what Charlie is measuring, they would fail. 

\begin{Le}
$\mathcal H$ supports Parameter Independence $\iff \mathcal H \models PI$ where $PI$ is the following formula. 
$$PI :\ \bigwedge_{U \subseteq C,C' \in \mathcal M} \bar u^C / \bar u^{C'} \perp_{z^C, z^{C'} / z^{C'}} \emptyset^{C'} / \emptyset^{C'} \ \overset{\text{not}}{\equiv} \bigwedge_{U \subseteq C, C' \in \mathcal M} \bar u^C \subseteq_{z^C, z^{C'}} \bar u^{C'}$$
\end{Le}

\begin{Ex}
Figure \ref{fig:parameter_independence} shows a hidden variable polyteam supporting Parameter Independence. Indeed, suffices to check it for $U = \{a,b\}$, since one may easily prove that it implies what is needed for $U = \{a\}$ and $U = \{b\}$. If we restrict to $z = \triangle$, both teams agree to say that the possible values for $ab$ are $(0,0)$, $(0,1)$, and $(1,0)$. If we restrict to $z = \triangledown$, they agree to say these values are now $(0,1)$, $(1,0)$ and $(1,1)$. There is nothing to check for $\triangleright$ because it is not a shared value of $z^{abc}$ and $z^{abd}$. 
\end{Ex}

\begin{Rem}
$PI$ is \emph{not} equivalent to $\displaystyle \bigwedge_{m \in C, C' \in \mathcal M} m^C / m^{C'} \perp_{z^C, z^{C'} / z^{C'}} \emptyset^{C'} / \emptyset^{C'}$. 
\end{Rem}

\begin{figure}
\centering
$\begin{array}{cccc}
a & b & c & z^{abc} \\ \hline 
0 & 0 & 1 & \triangle \\ 
0 & 1 & 0 & \triangle \\ 
1 & 0 & 0 & \triangle \\ 
1 & 0 & 1 & \triangle \\ 
0 & 1 & 1 & \triangledown \\ 
1 & 0 & 2 & \triangledown \\ 
1 & 1 & 2 & \triangledown
\end{array}$
\hspace{20pt}
$\begin{array}{cccc}
a & b & d & z^{abd} \\ \hline 
0 & 0 & 0 & \triangle \\ 
0 & 1 & 1 & \triangle \\ 
1 & 0 & 2 & \triangle \\ 
1 & 1 & 1 & \triangledown \\ 
1 & 0 & 2 & \triangledown \\ 
0 & 1 & 0 & \triangledown \\ 
0 & 0 & 1 & \triangleright
\end{array}$
\caption{A hidden variable polyteam supporting Parameter Independence}
\label{fig:parameter_independence}
\end{figure}

\subsubsection{Outcome Independence} \label{sec:outcome_independence}
\begin{Def}
Let $\mathcal H$ be a hidden variable model. $\mathcal H$ is said to support \emph{Outcome Independence} if for every $m \in C \in \mathcal M$ and every possible outcomes $o_m$ and $\bar o_{c \setminus m}$ that are both possible along with the hidden value $\lambda$, we have that all three may happen simultaneously. 
\end{Def}

\begin{Le}
$\mathcal H$ supports Outcome Independence $\iff \mathcal H \models OI$ where $OI$ is the following formula. 
$$OI :\ \bigwedge_{m \in C \in \mathcal M} m^C / m^C \perp_{z^C, z^C / z^C} \overline{c \setminus m}^C /\ \overline{c \setminus m}^C \ \overset{\text{not}}{\equiv} \bigwedge_{m \in C \in \mathcal M} m^C \perp_{z^C} \overline{c \setminus m}^C$$
\end{Le}

\begin{Fac} \label{fac:oi_equivalence}
$OI$ is equivalent to $\displaystyle \bigwedge_{U \subseteq C \in \mathcal M} \bar u^C / \bar u^C \perp_{z^C, z^C / z^C} \overline{c \setminus u}^C /\ \overline{c \setminus u}^C \ \overset{\text{not}}{\equiv} \bigwedge_{U \subseteq C \in \mathcal M} \bar u^C \perp_{z^C} \overline{c \setminus u}^C$. 
\end{Fac}

\begin{proof}
Let us denote by $\varphi$ this formula. The fact that $\varphi \models OI$ is trivial. Conversely, assume $\mathcal H \models OI$, let $U \subseteq C \in \mathcal M$ and $s, s' \in H_C$ be such that $s(z^C) = s'(z^C) =: \lambda$. Enumerate $U = \{m_1, \cdots, m_k\}$ and $C = \{m_1, \cdots, m_\ell\}$ and let us prove by induction that 
\begin{equation} \label{eqn:oi_equivalence_ih}
\exists s_n \in H_C, s_n(z^C) = \lambda \text{ and } \forall i \leq n, s_n(m_i) = \begin{cases} s(m_i) & \text{if } i \leq k \\ s'(m_i) & \text{if } i > k \end{cases}
\end{equation}
is true for all $1 \leq n \leq \ell$. First, $s_1 := s$ witnesses that (\ref{eqn:oi_equivalence_ih}) is true for $n = 1$. Then, let $n \leq \ell$ and assume (\ref{eqn:oi_equivalence_ih}) holds for $n-1$, witnessed by $s_{n-1}$. Let $t$ denote $s$ if $n \leq k$, and $s'$ if $n > k$, so that $t(m_n)$ is the desired value for $s_n(m_n)$. It holds that $s_{n-1}(z^C) = t(z^C) = \lambda$, hence by the $(m_n \in C)$-clause of $OI$ (which is assumed for $\mathcal H$), we may find $s_n \in H_C$ such that $s_n(z^C) = \lambda$, $s_n(m_n) = t(m_n)$ and $s_n$ agrees with $s_{n-1}$ elsewhere. $s_n$ witnesses (\ref{eqn:oi_equivalence_ih}) for $n$. By induction, (\ref{eqn:oi_equivalence_ih}) holds for $n = \ell$, and let $s''$ witness it. Then, $s''$ witnesses that the $(U \subseteq C)$-clause of $\varphi$ is true for $s, s'$. 
\end{proof}

\begin{Rem}
One may easily define an empirical version $EOI$ of Outcome Independence : it suffices to omit each occurrence of $z^C$ in $OI$. Doing so, $EOI$ is equivalent to $\texists z\ SV \wedge OI$.
\end{Rem}

\begin{Ex}
The polyteam from Figure \ref{fig:parameter_independence} does not support Outcome Independence, because (for example) in the first team, under $z = \triangle$, both $a$ and $b$ may be 1, but $ab$ can never be $(1,1)$. In the possibilistic collapse of Bell's model (shown in Figure \ref{fig:bell_model_possibilistic}), Empirical Outcome Independence does not hold either. Indeed, $a$ may be 0, $b$ may be 1, but $ab$ cannot be $(0,1)$. The Hardy model from Figure \ref{fig:hardy} also fails Empirical Outcome Independence. 
\end{Ex}

\subsubsection{$z$-Independence}
\begin{Def}
Let $\mathcal H$ be a hidden variable model. $\mathcal H$ is said to support \emph{$z$-Independence} whenever the possible hidden values for $z$ are the same whatever the measurement context $C$. 
\end{Def}

\begin{Le}
$\mathcal H$ supports $z$-Independence $\iff \mathcal H \models zI$ where $zI$ is the following formula. 
$$zI :\ \bigwedge_{C,C' \in \mathcal M} z^C / z^{C'} \perp \emptyset^{C'} / \emptyset^{C'}\ \overset{\text{not}}{\equiv} \bigwedge_{C, C' \in \mathcal M} z^C \subseteq z^{C'}$$
\end{Le}

\begin{Ex}
The hidden variable model from Figure \ref{fig:parameter_independence} does not support $z$-Independence, because the value $\triangleright$ appears for $z^{abd}$ but not for $z^{abc}$. On the contrary, the hidden variable model $\mathcal H$ from page \pageref{fig:ex1} does support $z$-Independence. 
\end{Ex}

\begin{Rem}
From a physical point of view, $z$-Independence is an extremely natural property to impose. Recall that we think of hidden variables as capturing the most fundamental attributes of a system, \textit{i.e.} what it is "in essence". It would then be very puzzling that the possible values for these attributes vary throughout the measurement contexts, basically by definition of what we mean by "in essence" : what is "fundamental" or "essential" to something does not depend on how exterior beings (us) interact with that thing. From a philosophical point of view, this is close to what is called \emph{realism} : reality is not (purely) a social or scientific construct, some objects have some attributes in themselves. 
\end{Rem}

\subsubsection{Locality and local realism} \label{sec:locality_local_realism}
As mentioned before, locality is intended to mean that the possible outcomes of a measurement can only be influenced by its immediate surroundings. In the simpler Bell-type case, this can be rephrased as "Alice's outcomes are independent from Bob's outcomes, for each value of the hidden variable." 

\begin{Def}
Let $\mathcal H$ be a hidden variable model. $\mathcal H$ is said to support \emph{Locality} whenever : 
\begin{multline*}
\forall C \in \mathcal M, \forall s \in H_C, \forall (s_m)_{m \in C} \text{ s.t. } s_m \in H_{C_m} \text{ where } m \in C_m \in \mathcal M, \\ \text{if } \left(\forall m \in C, s_m(z^{C_m}) = s(z^C)\right) \text{ then } \left(\exists s' \in H_C, \begin{array}[t]{l}
\forall m \in C, s'(m^C) = s_m(z^{C_m})\\ \text{and } s'(z^C) = s(z^C)
\end{array} \right)
\end{multline*}
If we allow ourselves to denote $H_C^\lambda[m^C]$ for the set of values $m^C$ may take along with $\lambda$ in $H_C$, and $\mathcal H^\lambda[m] := \displaystyle \bigcup_{C \ni m} H^\lambda_C[m^C]$, Locality reformulates as : 
$$\forall \lambda \in \Lambda, \forall C \in \mathcal M, \mathcal H^\lambda_C[\bar c] = \prod_{m \in C} \mathcal H^\lambda[m]$$
\end{Def}

\begin{Le}
$\mathcal H$ supports Locality $\iff \mathcal H \models Loc$ where $Loc := PI \wedge OI$.
\end{Le}

\begin{proof}
\begin{itemize}
\item[$\models$.] Assume $\mathcal H$ supports Locality. Let's first show $\mathcal H \models PI$. Let $U \subseteq C, C' \in \mathcal M$, and suppose $s \in H_C$, $s' \in H_{C'}$ are such that $s(z^C) = s'(z^{C'})$. By locality (take $s_m = s$ always) there is $s'' \in H_{C'}$ such that $s''(z^{C'}) = s(z^C)$ and $s'' \restriction U = s \restriction U$. Whence $\mathcal H \models \bar u^C \subseteq_{z^C, z^{C'}} \bar u^{C'}$, and $\mathcal H \models PI$. Secondly, let $m \in C \in \mathcal M$, and suppose $s,s' \in H_C$ are such that $s(z) = s'(z)$. By locality (take $s_m = s$ and $s_{m'} = s'$ if $m' \neq m$), there is $s'' \in H_C$ such that $s''(m) = s(m)$ and $s''$ agrees with $s'$ elsewhere. Whence $\mathcal H \models OI$ and $\mathcal H \models Loc$. 
\item[$\delsmo$.] Assume $\mathcal H \models PI \wedge OI$. Let $C \in \mathcal M$, $s \in H_C$, and $(s_m)_{m \in C}$ as in the Definition. The $(\{m\} \subseteq C_m,C)$-clause from $PI$ says there is $s_m^C \in H_C$ such that $s_m^C(m^C) = s_m(m^C)$ and $s_m^C(z^C) = s_m(z^{C_m})$. Now that all the $s_m^C$'s are in $H_C$, one can set up an induction argument that is very similar to that of the Fact at page \pageref{fac:oi_equivalence}. It follows that there is $s' \in H_C$ such that $s'(z) = s(z)$ and $s'$ regroups the measurement values of the $s_m^C$'s. All in all, $\mathcal H$ supports Locality. 
\end{itemize}
\end{proof}

\begin{Def}
$\mathcal H$ is said to support \emph{Local Realism} whenever $\mathcal H$ supports Locality and $z$-Independence. 
\end{Def}

\begin{Cor}
$\mathcal H$ supports Local Realism $\iff \mathcal H \models LocReal$ where $LocReal := PI \wedge OI \wedge zI$. 
\end{Cor}

\begin{Rem}
Again, we may easily define Empirical Locality as $ELoc := NS \wedge EOI$. As we have seen in Section \ref{sec:outcome_independence}, Empirical Outcome Independence (and by the same token Empirical Locality) is very likely to fail in the quantum systems we study. Thus, we are interested in the following question. Given an empirical model that does not support Empirical Locality, does it have a hidden variable realization that supports Local Realism? Put differently, does it satisfy $\texists z\ LocReal$ ?

Beware : in physics, "non-locality" most often means the negation of $\texists z\ LocReal$, as physicians see $z$-Independence as a triviality.  
\end{Rem}

\subsubsection{Local consistency}
It is a very notable fact that local consistency is expressible in our logics, as nothing has been done with that purpose in mind. As we have seen, the physical notion of no-signalling happens to coincide with local consistency in the empirical case. Since Parameter Independence is supposed to translate No-Signalling in the hidden variable framework, one expects it to be equivalent to local consistency of hidden variable models. This is almost the case. Indeed, Parameter Independence only speaks about the local consistency of component teams truncated to common values of $z$. One must further impose that those $z$-values are the same in any two component teams (\textit{i.e.} $z$-Independence). For instance the polyteam from Figure \ref{fig:parameter_independence} supports Parameter Independence without being locally consistent, because the $\{a,b,z^{abd}\}$-restriction of the second team yields an extra $(0,0,\triangleright)$. 

\begin{Le}
$\mathcal H$ is locally consistent $\iff \mathcal H \models PI \wedge zI$. 
\end{Le}

\begin{proof}
Easy verifications. 
\end{proof}

\begin{Fac}
$NS \equiv \texists z\ PI \wedge zI$
\end{Fac}

\begin{proof}
If $\mathcal L$ is locally consistent, then any hidden variable realization $\mathcal H$ of $\mathcal L$ supporting Single-Valuedness remains locally consistent. Conversely, if $\mathcal L$ has a locally consistent hidden variable realization $\mathcal H$, we have 
$$\text{res}^C_{C \cap C'}(L_C) = \text{res}^C_{C \cap C'}(H_C) \restriction (C \cap C') = \text{res}^{C'}_{C \cap C'}(H_{C'}) \restriction (C \cap C') = \text{res}^{C'}_{C \cap C'}(L_{C'})$$
\end{proof}

\subsubsection{Some remarks}
It is not clear what we are required to show for most of the lemmas in this section, because the physical property in discussion is sometimes vaguely defined. There are multiple ways of understanding these results. Either one consider them as definitions of the corresponding physical properties, in which case the only thing to care about is whether the formula accurately translates the informal property. Or, and we believe it is best, we understand them as follows. In \cite{abramsky_relational_2013}, Abramsky has defined the properties from $\mathcal P$ for relational models, in broad mathematical terms. It is much easier to be convinced of the faithfulness of these definitions in that framework, because broad mathematical language is much closer to physical language than our narrow syntax. Then, one may easily think of explicit reciprocate bijections to turn a relational model $e$ into a polyteam $\mathcal L(e)$, and to turn a polyteam $\mathcal L$ into a relational model $e(\mathcal L)$. Now, most of the lemmas from this section can be understood as follows : 
\begin{center}
$e$ supports $P \iff \mathcal L(e) \models \varphi_P$ where $\varphi_P$ is the formula for $P$ 
\end{center}
Most of these statements are easy to show. For example, let us prove Lemma \ref{le:ns_formula}. Beware : what we use is the "stronger version" of No-Signalling (\textit{cf.} \cite{abramsky_relational_2013} p.6).

\begin{proof} Let $\mathcal L$ be an empirical polyteam. 
\begin{align*}
\mathcal L \models NS & \iff \forall U \subseteq C,C' \in \mathcal M, \mathcal L \models \bar u^C \subseteq \bar u^{C'} \\ 
& \iff \forall U \subseteq C, C' \in \mathcal M, \forall s \in L_C, \exists s' \in L_{C'}, s'(\bar u^{C'}) = s(\bar u^C) \\ 
& \iff \forall C,C' \in \mathcal M, \forall U \subseteq C, \forall \bar o_u \subseteq O, (\bar o_u \in L_C[\bar u] \text{ and } U \subseteq C') \to \bar o_u \in L_{C'}[\bar u]\\ 
& \iff e(\mathcal L) \text{ supports No-signalling} 
\end{align*}
\end{proof}

\subsection{Logical entailments of physical properties} \label{sec:entailments}
Unsurprisingly, there are entailment relations between the physical properties from $\mathcal P$, other than those already shown. First, let us check that the team-based entailments from \cite{abramsky_team_2026} remain for polyteams. 

\begin{Le} \label{le:entailments}
$HWD \models OI$ and $HSD \equiv HWD \wedge PI$.
\end{Le}

\begin{proof}
Firstly, assume $\mathcal H \models HWD$, let $m \in C \in \mathcal M$ arbitrary. By Hidden Weak Determinism, any $s,s' \in H_C$ such that $s(z) = s'(z)$ must be equal, hence $\mathcal H \models m^C \perp_{z^C} \overline{c \setminus m}^C$ trivially holds. Thus $\mathcal H \models OI$. Secondly, assume $\mathcal H \models HSD$, and denote by $o_m$ the uniform output of $m \in M$ throughout $\mathcal H$. Note that $HSD$ is the conjunction of $HWD$ with some residue, whence $\mathcal H \models HWD$. Let $U \subseteq C,C' \in \mathcal M$. By Hidden Strong Determinism, any two $s \in H_C$, $s' \in H_{C'}$ such that $s(z^C) = s'(z^{C'})$ must be equal, \textit{i.e.} $\mathcal H \models \bar u^C \subseteq_{z^C, z^{C'}} \bar u^{C'}$. Thus $\mathcal H \models PI$. Thirdly, assume $\mathcal H \models HWD \wedge PI$, let $m \in C,C' \in \mathcal M$, and suppose $s \in H_C$, $s' \in H_{C'}$ are such that $s(z^C) = s(z^{C'})$. By $PI$, one may find $\hat s \in H_{C'}$ such that $\hat s(z^{C'}) = s(z^C)$ and $\hat s(m^{C'}) = s(m^C)$. Then $\hat s(z) = s'(z)$ and by Hidden Weak Determinism, it must be that $\hat s(m) = s'(m)$. In particular, $s(m^C) = s'(m^{C'})$. Thus $\mathcal H \models HSD$. 
\end{proof}

\begin{The}
An empirical model has a deterministic hidden variable realization (in the sense of Einstein-Podolsky-Rosen) if and only if it has a locally realistic hidden variable realization. In formulas : 
$$\texists z\ Det\ \overset{\text{not}}{\equiv}\ \texists z\ HSD \wedge zI \ \equiv \texists z\ PI \wedge OI \wedge zI\ \overset{\text{not}}{\equiv}\ \texists z\ LocReal$$
\end{The}

\begin{proof}
Let us adapt the proof of \cite{albert_unifying_2022} to polyteams. Lemma \ref{le:entailments} implies that $Det \models LocReal$, in a very straightforward way. The converse direction remains. Assume $\mathcal L$ is realized by $\mathcal H \models PI \wedge OI \wedge zI$. $\mathcal H$ satisfying Locality exactly means : 
$$\forall C \in \mathcal M, \forall \lambda \in \Lambda, H_C^{\lambda}[\bar c] = \prod_{m \in C} \mathcal H^\lambda[m]$$
All these sets are non-empty, since $zI$ forces the $H_C$'s to feature the same hidden values, and we recall that any $H_C$ is non-empty. Make the following definitions. 
\begin{align*}
\Lambda' & := \{(\lambda, t) \in \Lambda \times \mathcal E(M) : \forall m \in M, t(m) \in \mathcal H^\lambda[m]\}\\ 
H'_C & := \bigcup_{(\lambda,t) \in \Lambda'} \{s \in \mathcal E^h(C) : s(\bar c) = t(\bar c) \text{ and } s(z) = (\lambda,t)\}\\ 
\mathcal H' & := \{H_C, C \in \mathcal M\}
\end{align*}
First, if $m \in C,C' \in \mathcal M$ and $s \in H'_C, s' \in H'_{C'}$ are such that $s(z^C) = s'(z^{C'}) =: (\lambda,t)$, it must be that $s(m^C) = t(m) = s'(m^{C'})$. Hence $\mathcal H' \models HSD$. Secondly, let $C \in \mathcal M$ and suppose $(\lambda,t)$ occurs as a hidden value for $s' \in H'_{C'}$. By definition, it must hold that $\forall m \in M, t(m) \in \mathcal H^\lambda[m]$, so that $t(\bar c) \in \prod_{m \in C} \mathcal H^\lambda[m] = H_C^\lambda[\bar c]$ by Locality. This means that there is $s \in H_C$ such that $s(\bar c) = t(\bar c)$ and $s(z^C) = \lambda$. By definition, $s((\lambda,t)/z) \in H'_C$ as well. Thus $\mathcal H' \models zI$. Finally, it is easy to verify that $\mathcal H'$ realizes $\mathcal L$ from the assumption that $\mathcal H$ does. 
\end{proof}

\subsection{Formalization of global properties} \label{sec:global_properties_formulas}
Recall that our concern is to examine whether or not a given set of local information (a polyteam) coherently derives from global information about the system (a global section). It is thus strongly desirable, for our polyteam account to be somewhat comprehensive, that such global properties turn out to be expressible in our logics. This section deals with two notions of contextuality, namely non-locality and strong contextuality, in reverse order of strength, since the formulas involved are then of increasing complexity. 

\subsubsection{Strong contextuality}
\begin{Def}
Let $\mathcal L$ be an empirical polyteam. We say $s \in \mathcal E(M)$ is \emph{consistent with $\mathcal L$} whenever for all $C \in \mathcal M$, res$^M_C(s) \in L_C$. In words, $s$ is an empirically plausible global state of the system. $\mathcal L$ is called \emph{strongly contextual} if there is not even one such $s$. 
\end{Def}

\begin{Not}
Let $\top$ stand for any tautology, \textit{e.g.} $m^C = m^C$ for some $m \in C$. 
\end{Not}

\begin{Rem}
In the remaining of Section \ref{sec:theorems}, let us harmlessly assume that all the polyteams that may appear have non-empty components. Put differently, we temporarily rule out from the definition of polyteam those polyteams that have at least one empty component team. 
\end{Rem}

\begin{Le}
$\mathcal L$ is not strongly contextual $\iff \mathcal L \models NSC$ where $NSC$ is the following formula. 
$$NSC = (WD \wedge NS) \vee \top$$
\end{Le}

\begin{proof}
\begin{itemize}
\item[$\implies$.] Assume $\mathcal L$ is not strongly contextual, witnessed by $s \in \mathcal E(M)$. Then $\mathcal S:= \{\{s \restriction C\}, C \in \mathcal M\}$ is an empirical polyteam that conforms to previous Remark. The component teams of $\mathcal S$ are singletons, whence $\mathcal S \models WD$. Moreover, we have : 
$$\forall C, C' \in \mathcal M, S_C \restriction C \cap C' = S_{C'} \restriction C \cap C' = \{s \restriction C \cap C'\}$$
Hence, $\mathcal S \models NS$. Therefore $\mathcal L \models NSC$ with respect to the split $\mathcal L = \mathcal S \cup \mathcal L$. 
\item[$\impliedby$.] Assume $\mathcal L \models NSC$, witnessed by $\mathcal S \cup \mathcal L' = \mathcal L$. Since $\mathcal S \models WD$ and by previous Remark, $\mathcal S$ is of the form $\{\{s_C\}, C \in \mathcal M\}$. Now $\mathcal S \models NS$ exactly says that $\{s_C, C \in \mathcal M\}$ is a locally consistent set of events. Recall (Proposition \ref{pro:sheaf_events}) that $\mathcal E$ is a sheaf. Thus, there is $s \in \mathcal E(M)$ such that $\forall C \in \mathcal M, s \restriction C = s_C \in L_C$, and $\mathcal L$ is not strongly contextual.
\end{itemize}
\end{proof}

\begin{Rem}
The global meaning of $NSC$ stems from the global semantics of the splitjunction connective. 
\end{Rem}

\subsubsection{Global sections}
\begin{Pro}
$\mathcal L$ has a global section $\iff \mathcal L \models GS$ where $GS = \texists z\ Det \equiv \texists z\ LocReal$. 
\end{Pro}

\begin{proof}
\begin{itemize}
\item[$\implies$.] Assume $\mathcal L$ has a global section $G$. Consider the hidden variable global section $H = \{s \in \mathcal E^h(M) : s \restriction M \in G \text{ and } s(z) = s \restriction M\}$ obtained from $G$ by suffixing each assignment by itself. It is clear that $H$ realizes $G$. Let $\mathcal H$ be the contextualization of $H$, \textit{i.e.} $\mathcal H = \{\text{res}^M_C(H), C \in \mathcal M\}$. Then $H_C[z] = H[z] = G$ for each $C$, whence $\mathcal H \models zI$. Moreover, if $s \in H_C$ and $s' \in H_{C'}$ are such that $s(z^C) = s'(z^{C'}) = t \in G$, it must be that $s \restriction (C \cap C') = s' \restriction (C \cap C') = t \restriction (C \cap C')$, whence $\mathcal H \models HSD$. Now, Lemma \ref{le:polyteams_commutative_diagram} ensures that $\mathcal H$ realizes $\mathcal L$, \textit{i.e.} we have shown that $\mathcal L \models GS$. 
\item[$\impliedby$.] Assume $\mathcal L \models GS$, and let $\mathcal H$ be a realization of $\mathcal L$ such that $\mathcal H \models HSD \wedge zI$. For each $\lambda \in \mathcal H[z]$, and each $C \in \mathcal M$, let $s_C^\lambda \in H_C$ unique such that $s_C^\lambda(z^C) = \lambda$. Existence follows from $z$-Independence, and unicity from Hidden Strong Determinism. For each $\lambda$, the family $\{s_C^\lambda, C \in \mathcal M\}$ is locally consistent, by virtue of which it has a gluing $s^\lambda \in \mathcal E^h(M)$. Let $H = \{s^\lambda, \lambda \in \mathcal H[z]\}$. It is clear that $H$ is a hidden variable global section for $\mathcal H$. Let $G$ be the empirical global section realized by $H$. Lemma \ref{le:polyteams_commutative_diagram} ensures that $G$ contextualizes as $\mathcal L$, \textit{i.e.} $\mathcal L$ has a global section. 
\end{itemize}
\end{proof}

\begin{Fac}
$GS \models NSC$. In words, strong contextuality is a "wilder" form of contextuality than non-locality (the inexistence of global sections).
\end{Fac}

\begin{proof}
Any global section of $\mathcal L$ features assignments that are consistent with $\mathcal L$.
\end{proof}

\subsection{Logical formulations of contextuality theorems} \label{sec:formalization}

\subsubsection{The QM atom}
As mentioned before (Proposition \ref{pro:not_sheaf} and its remark), the inexistence of global sections is absolutely not surprising if we do not restrict the view to those polyteams that represent a quantum experiment. Indeed, any $\mathcal L \not\models NS$ exemplifies the wildest forms of contextuality one may think of. That is why, in order to provide a faithful logical account of contextuality theorems in quantum mechanics, we must introduce a formula $\varphi$ such that $\mathcal L \models \varphi \iff \mathcal L$ describes a quantum system. The only precise definition of the latter is through the Hilbert space formalism of quantum mechanics.

\begin{Def}
A \emph{complex Hilbert space} is a complex vector space $\mathbb H$ endowed with a hermitian product $\langle \cdot | \cdot \rangle$, \textit{i.e.} a bilinear form on $\mathbb H$ that is antilinear with respect to the first argument, linear with respect to the second argument, symmetric (when switching arguments conjugation appears), positive, and definite. An \emph{observable} is a bounded self-adjoint operator on $\mathbb H$. Two observables are called \emph{compatible} when they comute, in which case their product remains an observable. A \emph{positive operator-valued measure} on $\mathbb H$ is a set $\{A_i, i \in I\}$ of observables whose sum is $\text{id}_{\mathbb H}$.  
\end{Def}

\begin{Def} \label{def:finite_dim_quantum_system}
A \emph{finite-dimensional quantum system of type $(\mathcal M, O)$} is a tuple $\mathcal S = (\mathbb H, \{A_m^o,\ (m,o) \in M \times O\}, \rho)$, where : 
\begin{itemize}
\item $\mathbb H$ is a complex Hilbert space whose dimension is finite, 
\item For all $m \in M$, $\{A_m^o, o \in O\}$ is a positive operator-valued measure on $\mathbb H$. 
\item $\rho$ is a \emph{density operator} on $\mathbb H$, \textit{i.e.} 
$$\rho = \sum_{j < k} p_j \langle\ v_j\ |\ \cdot\ \rangle v_j \overset{\text{not}}{=} \sum_{j < k} p_j |v_j\rangle \langle v_j|$$
for some $k \in \mathbb N$, $v_j \in \mathbb H$, and $p_j$'s $\geq 0$ such that $\sum p_j = 1$. 
\end{itemize}
\end{Def}

\begin{Rem}
The second point implies that $ \displaystyle \left\{\prod_{m \in C} A_m^{o_m} \overset{\text{not}}{=} A_C^{\bar o},\ \bar o \in O^{|C|}\right\}$ is a p.o.v.m. on $\mathbb H$ as well. 
\end{Rem}

\begin{Fac}
According to quantum mechanics, given $\mathcal S$ a finite dimensional quantum system of type $(\mathcal M, O)$ and $C \in \mathcal M$, the probability of observing the outcome $\bar o$ while measuring $C$ is : 
$$p_C^{\mathcal S}(\bar o) := \text{Tr}(A_C^{\bar o} \cdot \rho) \hspace{30pt} \text{where Tr stands for the trace.}$$
\end{Fac}

The following definition anticipates on that of Section \ref{sec:semirings}, in the more specific case of possibilistic (as opposed to probabilistic) polyteams. 

\begin{Def}
Let $\mathcal L$ be an empirical polyteam. We say that $\mathcal L$ is a \emph{finite dimensional quantum mechanical polyteam} whenever there exists $\mathcal S$ a finite dimensional quantum system of type $(\mathcal M,O)$ such that $L_C = \{s \in \mathcal E(C) : p_C^{\mathcal S}(s(\bar c)) > 0\}$. Obviously, in the current state of our logics, we cannot hope to capture this notion with a formula. That is why we define the new atom $QM$ as follows : 
$$\mathcal L \models QM \iff \mathcal L \text{ is a finite dimensional quantum mechanical system of type } (\mathcal M,O)$$ 
\end{Def}

\begin{Le} \label{le:qm_ns}
$QM \models NS$. In words, all quantum systems are no-signalling. 
\end{Le}

\begin{proof}
Assume $\mathcal L$ is a finite dimensional quantum mechanical polyteam. Let $\mathcal S = (\mathbb H, \{A_m^o,\ (m,o) \in M \times O\}, \rho)$ witness it. Let $U \subseteq C, C' \in \mathcal M$, and suppose $\bar o_u \in O^{|U|}$ appears as a value for $\bar u^C$, \textit{i.e.} there is $\bar o_c \in O^{|C|}$ extending $\bar o_u$ such that Tr$(A_C^{\bar o_c}\rho) > 0$. If we denote by $E^C_{\bar o_u}$ the set of such $\bar o_c$'s, that is equivalent to : 
$$\text{Tr}(A_U^{\bar o_u}\rho) > 0 \hspace{20pt} \text{where } A_U^{\bar o_u} := \sum_{\bar o_c \in E_{\bar o_u}} A_C^{\bar o_c}$$
And this last claim does not depend on $C$, since : 
$$A_U^{\bar o_u} = \sum_{\bar o_c \in E^C_{\bar o_u}} \prod_{m \in C} A_m^{o_m} = \prod_{m \in C} \sum_{o_m \in O_m} A_m^{o_m} = \prod_{m \in U} A_m^{(\bar o_u)_m} \text{ where } O_m = \begin{cases} \{(\bar o_u)_m\} & \text{if } m \in U \\ O & \text{otherwise.} \end{cases}$$
\end{proof}

\begin{Ex}
The Hardy model (Figure \ref{fig:hardy}) satisfies $QM$, as proved in Theorem \ref{th:bell_possibilistic}.
\end{Ex}

\begin{Rem}
In quantum mechanics, the \emph{quantification principle} says that, when performing the measurement associated with the observable $A$, the only possible outcomes are the eigenvalues of $A$. However, it is easy to check that those eigenvalues must be real numbers. Indeed, assume $A|\psi\rangle = a |\psi\rangle$. On the one hand, composing with $\langle \psi|$ on the left, one gets $\langle \psi| A |\psi\rangle = a ||\psi||^2$. On the other hand, if one transconjugates and evaluates in $|\psi\rangle$, one gets $\langle \psi | A |\psi\rangle = a^*||\psi||^2$, thus $a = a^*$, and $a \in \mathbb R$ \textit{in fine}. 
\end{Rem}

\subsubsection{Bell's theorem}

\begin{The}[Bell, logical possibilistic version] \label{th:bell_possibilistic}
$QM \not\models GS$. In words, there is a non-local quantum system. 
\end{The}

\begin{proof}
Let $\mathcal L$ denote the Hardy model (see Figure \ref{fig:hardy}). We have already proved that $\mathcal L \not\models GS$ (Proposition \ref{pro:not_sheaf}). Remains to show $\mathcal L \models QM$. Let $\mathbb H = \mathbb C^2 \otimes \mathbb C^2$. Let $A_a^0$ (resp. $A_a^1$) be the orthogonal projector on $(1,0)$ (resp. $(0,1)$) tensored with $\text{id}_{\mathbb C^2}$. Let $A_{a'}^0$ (resp. $A_{a'}^1$) be the orthogonal projector on $(1/\sqrt 2, 1/\sqrt 2)$ (resp. $(-1/\sqrt 2, 1/\sqrt 2)$) tensored with $\text{id}_{\mathbb C^2}$, \textit{i.e.} the same as $a$ on the $\pi/4$-rotated basis. The $A_b$'s and $A_{b'}$'s are defined similarly, with $\text{id}_{\mathbb C^2} \otimes$ on the left this time. Finally, the density operator is $\rho := \langle\ v\ |\ \cdot\ \rangle v$ where $v = \frac{1}{\sqrt{12}}\left((1,0) \otimes (1,0) - (1,0) \otimes (0,1) - (0,1) \otimes (1,0) - 3 (0,1) \otimes (0,1)\right)$. In the canonical tensor-product basis, we have : 
$$\rho = \frac{1}{12}\left(\begin{array}{cccc}
1 & -1 & -1 & -3 \\ -1 & 1 & 1 & 3 \\ -1 & 1 & 1 & 3 \\ -3 & 3 & 3 & 9
\end{array}\right)$$
One can check this defines a finite dimensional quantum system which yields the Hardy model. For example : 
$$\text{Tr}(A_{a'}^0A_b^0\rho) = \text{Tr } \frac{1}{2} \left(\begin{array}{cc} 1 & 1 \\ 1 & 1 \end{array}\right) \otimes \left(\begin{array}{cc} 1 & 0 \\ 0 & 0 \end{array}\right) \rho = \left(\begin{array}{cccc} 0 & & \text{}_{(\ast)} & \\ & 0 & & \\ & & 0 & \\ & \text{}^{(\ast)} & & 0 \end{array}\right) = 0$$
\end{proof}

Let us now define another family of counter examples, which we shall later compare to the former. 

\begin{Def}
Assume $M \supseteq \{a_0,a_1,b_0,b_1,c_0,c_1\}$, $O \subseteq \{0,1\}$ and $\mathcal C \cup \{\{a_0, b_0, c_0\}\} \subseteq \mathcal M$, where $\mathcal C = \{\{a_0, b_1, c_1\}, \{a_1, b_0, c_1\}, \{a_1,b_1,c_1\}\}$. We say $\mathcal L$ is a \emph{GHZ polyteam} whenever : 
\begin{align*}
\forall C \in \mathcal C, L_C & = \{(0,0,0), (0,1,1), (1,0,1), (1,1,0)\} \\ 
\text{and } L_{\{a_0,b_0,c_0\}} & = \{(0,0,1), (0,1,0), (1,0,0), (1,1,1)\} 
\end{align*}
\end{Def}

\begin{Le}
All the GHZ polyteams satisfying $QM$ are counter examples for Bell's theorem. 
\end{Le}

\begin{proof}
See \cite{abramsky_relational_2013}, Section 6.2.
\end{proof}

\subsubsection{The Kochen-Specker theorem}

\begin{The}[Kochen-Specker, logical version]
$QM \not\models NSC$. In words, there is a strongly contextual quantum system.
\end{The}

In fact, we will prove the following stronger result : there exists a maximal context set $\mathcal M$ such that any no-signalling empirical polyteam (with contexts from $\mathcal M$) that satisfies $QM$ fails $NSC$. For this, we need several lemmas. 

\begin{Def} \label{def:MKS}
Let $M = \{m_1, \cdots, m_{18}\}$, $O = \{0,1\}$, and let $\mathcal M_{KS}$ denote the maximal context set given by the columns of the following table. 
$$\begin{array}{c|c|c|c|c|c|c|c|c}
m_1 & m_1 & m_8 & m_8 & m_2 & m_9 & m_{16} & m_{16} & m_{17} \\ 
m_2 & m_5 & m_9 & m_{11} & m_5 & m_{11} & m_{17} & m_{18} & m_{18} \\ 
m_3 & m_6 & m_3 & m_7 & m_{13} & m_{14} & m_4 & m_6 & m_{13} \\ 
m_4 & m_7 & m_{10} & m_{12} & m_{14} & m_{15} & m_{10} & m_{12} & m_{15}
\end{array}$$
\end{Def}

\begin{Fac}
$\mathcal M_{KS}$ has 9 columns, and each $m_i$ appears exactly twice. 
\end{Fac}

In the remaining of this paragraph, let us assume that $\mathcal M = \mathcal M_{KS}$, and that the languages $\sigma, \tau$ have a relation symbol $R$ which is interpreted in the underlying structure as : 

$$R^{\mathcal A} = \{(x_1,x_2,x_3,x_4) \in A^4 :\ \text{exactly one } x_i \text{ is 1}\}$$

\begin{Def}
Consider the following formula. 
$$ONE(\mathcal M_{KS}) := \bigwedge_{C = \{m_i, m_j, m_k, m_\ell\} \in \mathcal M_{KS}} R(m_i^C, m_j^C, m_k^C, m_\ell^C)$$
\end{Def}

\begin{Le} \label{le:ks_lemma_1}
$\mathcal L \models ONE(\mathcal M_{KS}) \implies \mathcal L \not\models NSC$. In words, $ONE(\mathcal M_{KS})$ entails strong contextuality. 
\end{Le}

\begin{proof}
Assume $\mathcal L \models ONE(\mathcal M_{KS})$, and suppose by contradiction that $s \in \mathcal E(M)$ is consistent with $\mathcal L$. Re-draw previous table, and replace each $m_i$ with $s(m_i) \in \{0,1\}$. Let $S \in \mathbb N$ denote the sum of the whole table. By assumption, each column sums to 1, and there are 9 columns, whence $S = 9$. However, each $m_i$ appears exactly twice, thus $S = 2 \times \sum_i s(m_i)$ is even. A contradiction.  
\end{proof}

\begin{Le} \label{le:18_vectors}
There is a set of 18 vectors $\mathcal V = \{v_1, \cdots, v_{18}\} \subseteq \mathbb R^4$ such that replacing each $m_i$ by $v_i$ in the table of Definition \ref{def:MKS} yields 9 orthonormal bases of $\mathbb C^4$ (each column being one). 
\end{Le}

\begin{proof}
One may explicit them. That is done in \cite{cabello_bell-kochen-specker_1996}. 
\end{proof}

\begin{Le} \label{le:ks_lemma_2}
$ONE(\mathcal M_{KS})\ \wedge\ NS \models QM$. In words, all the no-signalling polyteams with $ONE(\mathcal M_{KS})$ represent some quantum system.  
\end{Le}

\begin{proof}
Assume $\mathcal L$ is such. Let $\mathbb H = \mathbb C^4$ and consider for each $j \in \{1, \cdots, 18\}$ the orthogonal projectors $A_j^1$ on $v_j$ and $A_j^0$ on $v_j^\perp$, acting on $\mathbb H$. It is easily shown that each $A_j^\varepsilon$ is self-adjoint, and thus is an observable. Define : 
\begin{align*}
E & := \{j \in \{1, \cdots, 18\} : 1 \in \mathcal L[m_j]\}, \\
p_j & := \frac{1}{|E|} \text{ if } j \in E \text{, but } 0 \text{ otherwise,} \\
\rho & := \sum_{j \in E} p_j \langle\ v_j\ |\ \cdot\ \rangle v_j  \text{ , and} \\
\mathcal S & := (\mathbb H, (A_j^\varepsilon)_{\underset{\varepsilon \in \{0,1\}}{1 \leq j \leq 18}}, \rho)
\end{align*}
Let us show that $\mathcal S$ is a finite dimensional quantum mechanical system. Let $C \in \mathcal M_{KS}$, without loss of generality $C = \{m_1,m_2,m_3,m_4\}$, and $\bar \varepsilon \in \{0,1\}^4$. Given $j \in \{1,2,3,4\}$, $A_j^{\varepsilon_j}$ can be written in the basis $(v_1, v_2, v_3, v_4)$ as the diagonal matrix with value $\varepsilon_j$ in position $j$, and $1 - \varepsilon_j$ elsewhere. For example : 
$$A_2^{\varepsilon_2} = \left(\begin{array}{cccc} 1 - \varepsilon_2 & & & \\ & \varepsilon_2 & & \\ & & 1 - \varepsilon_2 & \\ & & & 1 - \varepsilon_2 \end{array}\right)$$ 
It follows that $A_j^0 + A_j^1 = \text{ id}$ for each $j$, and $\mathcal S$ is a finite dimensional quantum mechanical system. Let us now check that $\mathcal L$ represents $\mathcal S$. $A_C^{\bar \varepsilon} = A_1^{\varepsilon_1}A_2^{\varepsilon_2}A_3^{\varepsilon_3}A_4^{\varepsilon_4}$ can be written in the same basis as : 
\fontsize{7}{7}{$$
\left(\begin{array}{cccc}
\varepsilon_1(1-\varepsilon_2)(1-\varepsilon_3)(1-\varepsilon_4) & & & \\
& (1-\varepsilon_1) \varepsilon_2 (1-\varepsilon_3)(1-\varepsilon_4) & & \\ 
& & (1-\varepsilon_1)(1 - \varepsilon_2) \varepsilon_3 (1 - \varepsilon_4) \\ 
& & & (1 - \varepsilon_1)(1-\varepsilon_2)(1-\varepsilon_3)\varepsilon_4
\end{array}\right)
$$} \normalsize
It is easy to verify that $A_C^{\bar \varepsilon} \neq 0 \iff \bar \varepsilon$ has exactly one 1, in which case the diagonal $A_C^{\bar \varepsilon}$ has exactly one 1, and 0's elsewhere. Moreover, in the same basis : 
$$\rho = \left(\begin{array}{cccc} p_1 & & & \\ & p_2 & & \\ & & p_3 & \\ & & & p_4 \end{array}\right) \text{ and thus Tr}(A_C^{\bar \varepsilon}\rho) = \sum_{j = 1}^4 p_j \varepsilon_j \prod_{k \neq j} (1 - \varepsilon_k) \text{, so that}$$
\begin{align*}
\text{Tr}(A_C^{\bar \varepsilon}\rho) > 0 & \iff \exists ! j, \varepsilon_j = 1 \text{, and for that } j, p_j > 0 \\ 
& \iff \exists ! j, \varepsilon_j = 1 \text{, and for that } j, 1 \in \mathcal L[m_j] \\
& \iff \exists ! j, \varepsilon_j = 1 \text{, and for that } j, 1 \in L_C[m_j] \hspace{50pt} \text{(by No-Signalling)} \\ 
& \iff (m_1m_2m_3m_4 \longmapsto \varepsilon_1\varepsilon_2\varepsilon_3\varepsilon_4) \in L_C \hspace{50pt} \text{(because } \mathcal L \models ONE(\mathcal M_{KS})\text{)}
\end{align*}
This for all $C \in \mathcal M_{KS}$. Therefore $\mathcal L$ represents $\mathcal S$. All in all, $\mathcal L \models QM$.
\end{proof}

\begin{proof}[Proof of the Theorem]
It suffices to show that $ONE(\mathcal M_{KS}) \wedge NS$ is consistent, since then for any polyteam $\mathcal L$ satisfying this formula (in some structure), we would have $\mathcal L \models QM$ by Lemma \ref{le:ks_lemma_2} but $\mathcal L \not\models NSC$ by Lemma \ref{le:ks_lemma_1}. The following table describes a polyteam satisfying the desired formula (see next Section) : the assignments in $L_C$ are exactly those Dirac measures selecting one of the squared variables. 
$$\begin{array}{c|c|c|c|c|c|c|c|c}
\boxed{m_1} & \boxed{m_1} & \boxed{m_8} & \boxed{m_8} & \boxed{m_2} & m_9 & \boxed{m_{16}} & \boxed{m_{16}} & m_{17} \\ 
\boxed{m_2} & \boxed{m_5} & m_9 & m_{11} & \boxed{m_5} & m_{11} & m_{17} & m_{18} & m_{18} \\ 
m_3 & m_6 & m_3 & m_7 & m_{13} & m_{14} & m_4 & m_6 & m_{13} \\ 
m_4 & m_7 & m_{10} & m_{12} & m_{14} & \boxed{m_{15}} & m_{10} & m_{12} & \boxed{m_{15}}
\end{array}$$
\end{proof}

\subsubsection{The one-player game of $ONE(\mathcal M_{KS}) \wedge NS$}
Consider the following one-player game, which could be stated more precisely, if one wishes. 

\begin{itemize}
\item The table of Definition \ref{def:MKS} is the "board". It is played in 9 rounds (one per column). 
\item At round $i$, the player must circle a non-empty subset $S_i$ of the $i$-th column such that : 
\begin{enumerate}
\item for all $m_k$ in column $i$ appearing in a column $j < i$, $m_k \in S_i \iff m_k \in S_j$
\item for all $m_k$ in column $i$ appearing in a column $j < i$, $S_i = \{m_k\} \iff S_j = \{m_k\}$
\end{enumerate}
\item The player loses if they have no legal move at some point, and wins otherwise. 
\end{itemize}

\begin{Le}
The final positions of this game are in bijection with the empirical polyteams satisfying $ONE(\mathcal M_{KS}) \wedge NS$.
\end{Le}

\begin{proof}
If we forget about the game's conditions on $S_i$, it is clear that the set of final positions is in bijection with the polyteams satisfying $ONE(\mathcal M_{KS})$. Then, those conditions are equivalent to the no-signalling character of the resulting polyteam. Indeed, the first condition forces 1 to appear for $m_k$ if it has already happened in another context, and the second condition forces 0 to appear for $m_k$ if it has already happened before. 
\end{proof}

\begin{Ex}
The final position exhibited in the proof of previous theorem is a minimal winning strategy, meaning that each move played was $\subseteq$-minimal among the available moves. In particular, there is a winning strategy. 
\end{Ex}

\begin{Fac}
Lemma \ref{le:ks_lemma_1} reformulates as : there is no singleton-only winning strategy for the game. 
\end{Fac}

For example, in this position after the eighth move, there is no legal ninth move that is a singleton. 
$$\begin{array}{c|c|c|c|c|c|c|c|c}
\boxed{m_1} & \boxed{m_1} & \boxed{m_8} & \boxed{m_8} & m_2 & m_9 & \boxed{m_{16}} & m_{16} & m_{17} \\ 
m_2 & m_5 & m_9 & m_{11} & m_5 & m_{11} & m_{17} & \boxed{m_{18}} & m_{18} \\ 
m_3 & m_6 & m_3 & m_7 & \boxed{m_{13}} & m_{14} & m_4 & m_6 & m_{13} \\ 
m_4 & m_7 & m_{10} & m_{12} & m_{14} & \boxed{m_{15}} & m_{10} & m_{12} & m_{15}
\end{array}$$

The game has a pairs-only winning strategy (just circle the first two rows), but does not seem to have a triplets-only winning strategy.

\subsection{Comparisons with prior (team-based) literature} \label{sec:previous_work}
There is much to say, but let us be rather brief here. 

\paragraph*{} Firstly, in order to evaluate the benefits of this work, it is important to note that Albert and Grädel \cite{albert_unifying_2022} have already provided a team-based account of both Bell's and the Kochen-Specker theorem. We modestly believe that our polyteam account allows two slight improvements compared to \cite{albert_unifying_2022}. The first one is that the authors of \cite{albert_unifying_2022} have been led to introduce three new team semantical atoms to account for the Kochen-Specker theorem, apart from $QM$. One of these is a team version of the polydependence atom. On the contrary, we have not introduced any, that is why we are led to view our polyteam approach as being somewhat more natural. The second improvement is that our approach supports a wider diversity of maximal context sets. Indeed, the maximal context sets supported in \cite{albert_unifying_2022} are exactly those whose contexts have a uniform cardinality. This formalism may be interpreted as accounting for experiments in which there is a fixed number of agents (either spatially separated or not), each of whom performs a measurement, and in which "who performs measurement $m$" may matter. Even though our counter examples involve maximal context sets that are such, some more generality might turn out to be useful in future works. 

\paragraph*{} Secondly, this work is to be thought in the lineage of \cite{abramsky_team_2026} (Abramsky, Väänänen, Puljujärvi, 2026). The authors have provided a team semantical interpretation of the relational models found in \cite{abramsky_relational_2013}. They have succesfully proved a logical version of Bell's theorem in that framework, without introducing alien atoms. However, they have not been able to do the same with the Kochen-Specker theorem, they have given the foundational idea of this work (using polyteams) instead. Moreover, they have been investigating the implication problem, and listing sound rules from which they were able to derive syntactically several semantical implications. We were able to show there is a translation of empirical teams, hidden variable teams, physical properties and most team formulas appearing in \cite{abramsky_team_2026} as polyteams and polyteam formulas, respecting the satisfaction relation. We omit the details here. Simply put, the team framework only accounts for the so called \emph{Bell-type} polyteam models.

\paragraph*{} Thirdly, the authors of \cite{barlag_locally_2026} have defined and studied team-syntax dependency atoms in what they call \emph{contextual $K$-families}, where $K$ is a semiring. Let us focus on $K = \mathbb B = \{0,1\}$ seen as a boolean algebra, so that contextual $\mathbb B$-families are exactly locally consistent (possibilistic) polyteams. In short, when $\mathcal L$ is locally consistent, $\mathcal L \models\ \dep(\bar u, \bar v)$ is unambiguous because $L_C \models\ \dep(\bar u^C, \bar v^C)$ has the same truth value regardless of $C \supseteq \bar u\bar v$. It is important to note that, although this syntax is simpler to deal with, transitivity now fails. In line with \cite{barlag_locally_2026}, we have made sense of all team-syntax independence logic formulas in locally consistent polyteams and shown that every such formula has an equivalent formula in regular polyteam syntax. How much of the converse holds is an interesting question. Notably, all the physical properties from $\mathcal P$, as well as $GS$ and $NSC$, are expressible in locally consistent polyteam logics. We believe it would be interesting to tackle syntactical deduction in this simplified syntax.

\paragraph*{} Let us end this section with a table summarizing the formulas involved.
\newpage
\begin{table}[h]
\centering
\renewcommand{\arraystretch}{2}
\begin{tabular}{ccc}
\textbf{Physical properties} & \underline{Polyteam semantics} & \underline{LC polyteam semantics} \\ 
Weak Determinism & $\displaystyle \bigwedge_{m \in C \in \mathcal M} \dep(\emptyset^C, m^C / \emptyset^C, m^C)$ & $\displaystyle \bigwedge_{m \in M} \dep(m)$\\ 
Strong Determinism & $\displaystyle \bigwedge_{m \in C,C' \in \mathcal M} \dep(\emptyset^C, m^C / \emptyset^{C'}, m^{C'})$ & $\displaystyle \bigwedge_{m \in M} \dep(m)$ \\ 
No-Signalling & $\displaystyle \bigwedge_{U \subseteq C,C' \in \mathcal M} \bar u^C \subseteq u^{C'}$ & $\top$ \\ 
Global Section & $\texists z\ HSD \wedge zI$ & $\texists z\ HSD$ \\ 
Not Strongly Contextual & $(WD \wedge NS) \vee \top$ & $WD \vee \top$ \\ 
Hidden Weak Determinism & $\displaystyle \bigwedge_{C \in \mathcal M} \dep(z^C, \bar c / z^C, \bar c)$ & $\displaystyle \bigwedge_{m \in M} \dep(z,m)$ \\ 
Hidden Strong Determinism & $\displaystyle \bigwedge_{U \subseteq C, C' \in \mathcal M} \dep(z^C, \bar u^C / z^{C'}, \bar u^{C'})$ & $\displaystyle \bigwedge_{m \in M} \dep(z, m)$ \\ 
Single-Valuedness & $\displaystyle \bigwedge_{C, C' \in \mathcal M} \dep(\emptyset^C, z^C / \emptyset^{C'}, z^{C'})$ & $\dep(z)$ \\ 
Parameter Independence & $\displaystyle \bigwedge_{U \subseteq C,C' \in \mathcal M} \bar u^C \subseteq_{z^C, z^{C'}} \bar u^{C'}$ & $\top$ \\ 
Outcome Independence & $\displaystyle \bigwedge_{m \in C \in \mathcal M} m^C / m^C \perp_{z^C, z^C / z^C} \overline{c \setminus m}^C / \overline{c \setminus m}^C$ & $\displaystyle \bigwedge_{m \in C \in \mathcal M} m \perp_z \overline{c \setminus m}$ \\ 
$z$-Independence & $\displaystyle \bigwedge_{C,C' \in \mathcal M} z^C \subseteq z^{C'}$ & $\top$ \\ 
Locality & $PI \wedge OI$ & $OI$ \\ 
Local Realism & $PI \wedge OI \wedge zI$ & $OI$ \\ 
(Hidden) Local Consistency & $PI \wedge zI$ & $\top$
\end{tabular}
\caption{Formulas for physical properties in polyteam semantics}
\label{tab:physical_properties}
\end{table}

\newpage
\section{Semiring- and semifield-weighted polyteams} \label{sec:semirings}
This section aims at introducing logics of dependence and independence in semiring-weighted polyteams, and studying their general properties. As far as we know, this has not been done before in polyteams. There is already some work about semiring-weighted teams though, in particular probabilistic teams, and we will get inspiration from there. 

\subsection{$R$-weighted teams and polyteams}
\begin{Def}
A (non-trivial) \emph{semiring} is a tuple $(R,+,\times,0,1)$ where $R$ is a set, $+$ and $\times$ are binary operations on $R$ such that both $(R,+,0)$ and $(R,\times,1)$ are commutative monoids, $\times$ distributes on $+$, and $0 \neq 1$.
\end{Def}

\begin{Fac}
Let $R$ be a semiring. $R$ has a natural partial ordering given by : 
$$a \leq b \iff \exists c \in R,\ a + c = b$$
\end{Fac}

\begin{proof}
Reflexivity is seen with $c = 0$, and transitivity stems from the associativity of $+$.
\end{proof}

\begin{Ex}
Rings are semirings. 
\end{Ex}

\begin{Def}
A semiring $R$ is called \emph{positive} whenever $a + b = 0 \implies a = b = 0$. 
\end{Def}

\begin{Ex}
Typical positive semirings are the boolean algebra $\mathbb B = \{0,1\}$, $\mathbb N$, $\mathbb Q_+$, $\mathbb R_+$, and the tropical semiring $\mathbb T = (\mathbb R \cup \{+\infty\}, \min, +, +\infty, 0)$.
\end{Ex}

\begin{Def}
Let $R$ be a semiring, and $X$ be a finite set. An $R$-distribution on $X$ is a map $d : X \to R$ such that 
$$\sum_{x \in X} d(x) = 1$$
The set $\supp(d) := \{x \in X : d(x) \neq 0\}$ is called the \emph{support} of $d$. The set of $R$-distributions on $X$ is denoted by $\mathcal D_RX$.
\end{Def}

\begin{Rem}
It is very natural to extend $d$ to $\mathcal P(X)$ with $d(U) = \sum_{x \in U} d(x)$, in order to speak about the total weight of an aggregate of points. In various scenarios, and indeed for our purposes, it must hold that any aggregate of non-zero-weight points remains of non-zero weight. That is precisely what positivity asserts, and that is why, unless otherwise specified, all the semirings from Section \ref{sec:semirings} will be assumed positive. Intuitively, positivity prompts us to think of weights as degrees of existence. 
\end{Rem}

\begin{Ex}
The sets $\mathcal D_{\mathbb B}X, \mathcal D_{\mathbb N}X, \mathcal D_{\mathbb R_+}X$ are respectively in bijection with the set $\mathcal P^+(X)$ of non-empty subsets of $X$, the set of multisubsets of $X$, and the set of probability distributions on $X$.  
\end{Ex}

\begin{Def}
Let $\mathcal A$ be a structure, $\bar x$ be variables with arbitrary sorts, and $R$ be a semiring. Denote by $\bar x \to A$ the set of assignments on $\mathcal A$ with domain $\bar x$. An \emph{$R$-weighted team} (or simply \emph{$R$-team}) is any $X \in \mathcal D_R(\bar x \to A)$. 
\end{Def}

\begin{Def}
An \emph{$R$-weighted polyteam} (or simply $R$-polyteam) on $\mathcal A$ with domain $\bar x$ labeled as $\chi : \bar x \to \mathcal C$ is a set of $R$-teams $\mathcal L = \{L_C, C \in \mathcal M\}$ such that $\dom(L_C) = \{v \in \bar x : \chi(v) = C\}$ for each $C \in \mathcal M$.
\end{Def}

\begin{Rem}
$X$ or $\mathcal L$ is called possibilistic or probabilistic respectively when $R = \mathbb B$ and $R = \mathbb R_+$.
\end{Rem}

\begin{Ex}
Figure \ref{fig:bell_model} shows a probabilistic polyteam. The missing assignments are implicitly given probability 0. Since, in this example, all probabilities are rational, it may be seen as a $\mathbb Q_+$-polyteam. Note that $\mathbb Q_+$-(poly)teams are essentially the same as the so called \emph{multi(poly)teams}, \textit{i.e.} $\mathbb N$-weighted (poly)teams without the normalisation condition, once all weights are set to the same denominator.  
\begin{figure} 
\centering
$\begin{array}{ccc}
a & b & \mathbb R_+ \\ \hline
0 & 0 & 1/2 \\ 
1 & 1 & 1/2
\end{array}$
\hspace{20pt}
$\begin{array}{ccc}
a' & b & \mathbb R_+ \\ \hline
0 & 0 & 3/8 \\ 
0 & 1 & 1/8 \\ 
1 & 0 & 1/8 \\ 
1 & 1 & 3/8
\end{array}$
\hspace{20pt}
$\begin{array}{ccc}
a & b' & \mathbb R_+ \\ \hline
0 & 0 & 3/8 \\ 
0 & 1 & 1/8 \\ 
1 & 0 & 1/8 \\ 
1 & 1 & 3/8
\end{array}$
\hspace{20pt}
$\begin{array}{ccc}
a' & b' & \mathbb R_+ \\ \hline
0 & 0 & 1/8 \\ 
0 & 1 & 3/8 \\ 
1 & 0 & 3/8 \\ 
1 & 1 & 1/8
\end{array}$
\caption{Bell's model}
\label{fig:bell_model}
\end{figure}
\end{Ex}

\begin{Def}
Let $X$ be an $R$-team on domain $\bar v$, and let $\bar u \subseteq \bar v$. The weight map $X \restriction \bar u$ defined for $s : \bar u \to A$ by : 
$$(X \restriction \bar u)(s) := \sum_{t \in X \text{ extending } s} X(t)$$
is an $R$-team, called the \emph{restriction} of $X$ to $\bar u$. 
\end{Def}

\begin{Not}
We denote $X(\bar u = \bar o)$ or $X(\bar u \mapsto \bar o)$ for $(X \restriction \bar u)(\bar u \mapsto \bar o)$.
\end{Not}

\begin{Def}
Assume we are working with the language $\sigma$, the maximal context set $\mathcal M$, and the sheaf of events $\mathcal E$ as in previous sections. Then, if $U \subseteq C \in \mathcal M$, elements of $\mathcal D_R\mathcal E(U)$ may be identified with $R$-teams with domain $\bar u^C$. An \emph{empirical $R$-polyteam} is an $R$-polyteam $\mathcal L = \{L_C, C \in \mathcal M\}$ such that $L_C \in \mathcal D_R\mathcal E(C)$ for each $C$. If $\sigma, \mathcal E$ are replaced with $\tau, \mathcal E^h$, one gets the definition of \emph{hidden variable $R$-polyteam}.
\end{Def}

\begin{Fac}
$\mathcal D_R\mathcal E$ and $\mathcal D_R\mathcal E^h$ are presheaves. 
\end{Fac}

\begin{proof}
Restriction maps are given by res$^U_{U'}(X) := X \restriction \bar u'$ and res$^U_{U'}(Y) := Y \restriction \bar u'z$. The proof then consists in elementary manipulations of sums. 
\end{proof}

\begin{Fac}
Bell's model (Figure \ref{fig:bell_model}) is locally consistent. 
\end{Fac}

\begin{proof}
Any measurement in any context takes value 0 or 1 both with probability $1/2$.
\end{proof}

\subsection{Dependence and independence atoms in $R$-weighted polyteams} \label{sec:R_polyatoms}
It is unclear at first sight how dependence atoms should be defined in $R$-(poly)teams. That is why we will define them from independence atoms, so that they still satisfy Lemma \ref{le:polydep_atoms_as_polyindep_atoms}. As for independence atoms in $R$-polyteams, there are some constraints coming from prior literature we would like to fit in, so that the theory of team and polyteam logics remain as consistent as possible. 

\begin{Def}
Let $\bar u, \bar v, \bar w$ be variables. The syntactical expression $\bar v \independent[R][\bar u] \bar w$ is called an \emph{$R$-independence atom} for teams. It is given (for example in \cite{abramsky_team_2026}), the following semantics. For any structure $\mathcal A$, and any $X$ $R$-team on $\mathcal A$ whose domain contains the variables $\bar u, \bar v, \bar w$ : 
$$\mathcal A \models_X \bar v \independent[R][\bar u] \bar w \iffdef \forall \bar o_u\bar o_v\bar o_w \subseteq O, X(\bar u\bar v\bar w = \bar o_u\bar o_v\bar o_w)X(\bar u = \bar o_u) = X(\bar u\bar v = \bar o_u\bar o_v)X(\bar u\bar w = \bar o_u\bar o_w)$$
\end{Def}

\begin{Rem}
Let us see why this is a natural way to define an independence notion. Suppose $R = \mathbb R_+$. Note that if $X(\bar u = \bar o_u) = 0$, the right hand side equality holds, by a crucial use of $R$'s positivity. Otherwise, the equality is equivalent to : 
\begin{align*}
\frac{X(\bar u\bar v\bar w = \bar o_u\bar o_v\bar o_w)}{X(\bar u = \bar o_u)} & = \frac{X(\bar u\bar v = \bar o_u\bar o_v)}{X(\bar u = \bar o_u)} \times \frac{X(\bar u\bar w = \bar o_u\bar o_w)}{X(\bar u = \bar o_u)} \\
i.e. \text{ to } X(\bar v\bar w = \bar o_v\bar o_w\ |\ \bar u = \bar o_u) & = X(\bar v = \bar o_v\ |\ \bar u = \bar o_u) \times X(\bar w = \bar o_w\ |\ \bar u = \bar o_u)
\end{align*}
where $X(E\ |\ F)$ denotes the $X$-probability of $E$ conditioned to $F$. All in all : 
\begin{multline*}
\mathcal A \models_X \bar v \independent[R][\bar u] \bar w \iff \text{The random variables } \bar v \text{ and } \bar w \text{ are independent when conditioned to } \bar u \text{,} \\ \text{with respect to the probability measure } X
\end{multline*}
The equivalent formulation of the Definition simply allows to define this notion in semirings that are not semifields, \textit{i.e.} the Definition does not need a multiplicative inverse for $X(\bar u = \bar o_u)$. 
\end{Rem}

\begin{Rec}
The semantics of independence atoms in possibilistic polyteams are so defined : 
$$\hspace{-35pt}\mathcal A \models_{\mathcal L} \bar y^C / \bar v^{C''} \perp_{\bar x^C, \bar a^{C'} / \bar u^{C''}} \bar b^{C'} / \bar w^{C''} \iffdef \forall s \in L_C, \forall s' \in L_{C'}, s(\bar x^C) = s'(\bar a^{C'}) \to \exists s'' \in L_{C''}\ \begin{array}[t]{|l} s''(\bar u^{C''}) = s'(\bar a^{C'}) = s(\bar x^C),\\ s''(\bar v^{C''}) = s(\bar y^C), \text{ and} \\ s''(\bar w^{C''}) = s'(\bar b^{C'}) \end{array}$$
\end{Rec}

Ideally, we would like to define an $R$-independence atom $\bar y^C / \bar v^{C''} \independent[R][\bar x^C, \bar a^{C'} / \bar u^{C''}] \bar b^{C'} / \bar w^{C''}$ such that : 
\begin{enumerate}[label=(\alph*)]
\item If $\bar x^C = \bar a^{C'} = \bar u^{C''}$, $\bar y^C = \bar v^{C''}$ and $\bar b^{C'} = \bar w^{C''}$, then : 
$$\mathcal L \models \bar v^C / \bar v^C \independent[R][\bar u^C, \bar u^C / \bar u^C] \bar w^C / \bar w^C \iff L_C \models \bar v^C \independent[R][\bar u^C] \bar w^C$$
\item If $R = \mathbb B$, then : 
$$\mathcal L \models \bar y^C / \bar v^{C''} \independent[\mathbb B][\bar x^C, \bar a^{C'} / \bar u^{C''}] \bar b^{C'} / \bar w^{C''} \iff \mathcal L \models \bar y^C / \bar v^{C''} \perp_{\bar x^C, \bar a^{C'} / \bar u^{C''}} \bar b^{C'} / \bar w^{C''}$$
\end{enumerate}

\begin{Fac}
(a) and (b) are consistent, meaning $L_C \models \bar v^C \independent[\mathbb B][\bar u^C] \bar w^C \iff \mathcal L \models \bar v^C / \bar v^C \perp_{\bar u^C, \bar u^C / \bar u^C} \bar w^C / \bar w^C$.
\end{Fac}

\begin{proof}
\begin{itemize}
\item[$\implies$.] Assume $L_C \models \bar v \independent[\mathbb B][\bar u] \bar w$, \textit{i.e.} $\forall \bar o_u\bar o_v\bar o_w, L_C(\bar u\bar v\bar w = \bar o_u\bar o_v\bar o_w) L_C(\bar u = \bar o_u) = L_C(\bar u\bar v = \bar o_u\bar o_v) L_C(\bar u\bar w = \bar o_u\bar o_w) \in \mathbb B$. Suppose $s,s' \in L_C$ are such that $s(\bar u) = s'(\bar u)$. The hypothesis for $\bar o_u = s(\bar u) = s'(\bar u)$, $\bar o_v = s(\bar v)$ and $\bar o_w = s'(\bar w)$ says that $L_C(\bar u\bar v\bar w = \bar o_u\bar o_v\bar o_w) = 1$, \textit{i.e.} that there is $s''$ as desired. 
\item[$\impliedby$.] Assume $L_C \models \bar v \perp_{\bar u} \bar w$, let $\bar o_u\bar o_v\bar o_w \subseteq O$ and let us show that $L_C(\bar u\bar v\bar w = \bar o_u\bar o_v\bar o_w) L_C(\bar u = \bar o_u) = L_C(\bar u\bar v = \bar o_u\bar o_v) L_C(\bar u\bar w = \bar o_u\bar o_w)$. The $\leq$ direction is obvious, and the $\geq$ direction is a reformulation of the hypothesis. 
\end{itemize}
\end{proof}

Thus, in principle, it is possible to define $\independent[R]$ such that (a) and (b) hold. However, both the natural generalization of (a) to arbitrary sorts and the natural generalization of (b) to any $R$ fail to satisfy (a) and (b). Indeed, the natural generalization of (b) would be the following : 
\begin{multline}
\mathcal L \models \bar y^C / \bar v^{C''} \independent[R][\bar x^C, \bar a^{C'} / \bar u^{C''}] \bar b^{C'} / \bar w^{C''} \iffdef \forall \bar o_{xau}\bar o_{yv}\bar o_{bw} \subseteq O,\\ L_C(\bar x\bar y = \bar o_{xau}\bar o_{yv}) L_{C'}(\bar a\bar b = \bar o_{xau}\bar o_{bw}) \leq L_{C''}(\bar u\bar v\bar w = \bar o_{xau}\bar o_{yv}\bar o_{bw})
\end{multline}
In this case (a) asks that for all $\bar o_{xau}\bar o_{yv}\bar o_{bw} \subseteq O$ : 
\begin{multline*}
L_C(\bar u\bar v = \bar o_{xau}\bar o_{yv}) L_C(\bar u\bar w = \bar o_{xau}\bar o_{bw}) \leq L_C(\bar u\bar v\bar w = \bar o_{xau}\bar o_{yv}\bar o_{bw}) \\ \iff L_C(\bar u\bar v\bar w = \bar o_{xau}\bar o_{yv}\bar o_{bw})L_C(\bar u = \bar o_{xau}) = L_C(\bar u\bar v = \bar o_{xau}\bar o_{yv})L_C(\bar u\bar w = \bar o_{xau}\bar o_{bw})
\end{multline*}
which fails, because of the left-to-right direction. Now, consider the natural generalization of (a) : 
\begin{multline} \label{eqn:R_polyindependence_atom}
\mathcal L \models \bar y^C / \bar v^{C''} \independent[R][\bar x^C, \bar a^{C'} / \bar u^{C''}] \bar b^{C'} / \bar w^{C''} \iffdef \forall \bar o_{xau}\bar o_{yv}\bar o_{bw} \subseteq O, L_{C''}(\bar u\bar v\bar w = \bar o_{xau}\bar o_{yv}\bar o_{bw})L_C(\bar x = \bar o_{xau})L_{C'}(\bar a = \bar o_{xau}) \\ = L_C(\bar x\bar y = \bar o_{xau}\bar o_{yv}) L_{C'}(\bar a\bar b = \bar o_{xau}\bar o_{bw}) L_{C''}(\bar u = \bar o_{xau})
\end{multline}

\begin{Not}
To lighten notations, we omit $\bar o$'s starting from now, unless ambiguity shows up. 
\end{Not}
In the latter case, (b) asks that 
$$L_{C''}(\bar u\bar v\bar w)L_C(\bar x)L_{C'}(\bar a) = L_C(\bar x\bar y) L_{C'}(\bar a\bar b) L_{C''}(\bar u) \iff L_C(\bar x\bar y) L_{C'}(\bar a\bar b) \leq L_{C''}(\bar u\bar v\bar w)$$
which fails, in both directions. A detailed analysis (that we choose to silence here) shows that the direct sense would hold under the assumption that common values of $\bar x$ and $\bar a$ are values for $\bar u$. More generally, the equivalence would hold under the assumption that $L_C \restriction \bar x\bar y \simeq L_{C''} \restriction \bar u\bar v$ and $L_{C'} \restriction \bar a\bar b \simeq L_{C''} \restriction \bar u\bar w$, which should be construed as a local consistency hypothesis (see Lemma \ref{le:B_indep_atoms_coincidence}).

\paragraph*{} Finally, we are more enclined to accept a slight difference between $\independent[\mathbb B]$ and $\perp$ than a non-uniform definition of $\independent[R]$. Indeed, we aim in Section \ref{sec:collapses_realizations} at establishing connections between $R$-truth and $R'$- (typically $\mathbb B$-) truth. This task would be much less fruitful if $\independent[R]$ was defined separating the case $R = \mathbb B$. 

\begin{Def}
Consider the labeled variables $|\bar x^C| = |\bar a^{C'}| = |\bar u^{C''}|$, $|\bar y^C| = |\bar v^{C''}|$ and $|\bar b^{C'}| = |\bar w^{C''}|$. The syntactical expression $\bar y^C / \bar v^{C''} \independent[R][\bar x^C, \bar a^{C'} / \bar u^{C''}] \bar b^{C'} / \bar w^{C''}$ is called an \emph{$R$-polyindependence atom}. Its semantics is defined by Equation (\ref{eqn:R_polyindependence_atom}), \textit{i.e.} in simplified notations : 
\begin{multline*}
\mathcal L \models \bar y^C / \bar v^{C''} \independent[R][\bar x^C, \bar a^{C'} / \bar u^{C''}] \bar b^{C'} / \bar w^{C''} \iffdef \forall \bar o_{xau}\bar o_{yv}\bar o_{bw} \subseteq O,\\ L_{C''}(\bar u\bar v\bar w)L_C(\bar x)L_{C'}(\bar a) = L_C(\bar x\bar y) L_{C'}(\bar a\bar b) L_{C''}(\bar u)
\end{multline*}
\end{Def}

\begin{Rem}
Again, let us suppose $R = \mathbb R_+$. If $L_C(\bar x)$, $L_{C'}(\bar a)$ or $L_{C''}(\bar u)$ is zero, the right-hand side equality holds (again, positivity is crucial). Otherwise, this equality is equivalent to : 
$$L_{C''}(\bar v\bar w\ |\ \bar u) = L_C(\bar y\ |\ \bar x) L_{C'}(\bar b\ |\ \bar a)$$
and the notion of independence thus obtained does generalize that of $R$-teams (see previous Remark). 
\end{Rem}

\begin{Ex}
Figure \ref{fig:R_polyindependence_atom} shows a probabilistic polyteam $\mathcal L$, and $\mathcal L \models y / v \independent[\mathbb R_+][x, a / u] b / w$. Indeed : 
\begin{align*}
L_{C''}(uvw = 001)L_C(x=0)L_{C'}(a=0) = L_C(xy=00)L_{C'}(ab=01)L_{C''}(u=0) & \iff \frac{1}{4} \times \frac{1}{3} \times \frac{1}{2} = \frac{1}{3} \times \frac{1}{2} \times \frac{1}{4} \\ 
L_{C''}(uvw = 110)L_C(x=1)L_{C'}(a=1) = L_C(xy=11)L_{C'}(ab=10)L_{C''}(u=1) & \iff \frac{1}{2} \times \frac{1}{3} \times \frac{1}{2} = \frac{1}{3} \times \frac{1}{3} \times \frac{3}{4}\\
L_{C''}(uvw = 112)L_C(x=1)L_{C'}(a=1) = L_C(xy=11)L_{C'}(ab=12)L_{C''}(u=1) & \iff \frac{1}{4} \times \frac{1}{3} \times \frac{1}{2} = \frac{1}{3} \times \frac{1}{6} \times \frac{3}{4} \\  
\end{align*}

\begin{figure}
\centering 
$\begin{array}{ccc}
x & y & \mathbb R_+ \\ \hline
0 & 0 & 1/3 \\ 
1 & 1 & 1/3 \\ 
2 & 0 & 1/3
\end{array}$
\hspace{20pt}
$\begin{array}{ccc}
a & b & \mathbb R_+ \\ \hline
0 & 1 & 1/2 \\ 
1 & 0 & 1/3 \\ 
1 & 2 & 1/6
\end{array}$
\hspace{20pt}
$\begin{array}{cccc}
u & v & w & \mathbb R_+ \\ \hline
0 & 0 & 1 & 1/4 \\
1 & 1 & 0 & 1/2 \\ 
1 & 1 & 2 & 1/4
\end{array}$
\caption{A probabilistic polyteam satisfying a polyindependence atom}
\label{fig:R_polyindependence_atom}
\end{figure}
\end{Ex}

\begin{Fac}[Recovery of (a)]
Further assume that $R$ is multiplicatively cancellative. Then : 
$$\mathcal L \models \bar v^C / \bar v^C \independent[R][\bar u^C, \bar u^C / \bar u^C] \bar w^C / \bar w^C \iff L_C \models \bar v^C \independent[R][\bar u^C] \bar w^C$$
\end{Fac}

\begin{proof}
Easy verifications. 
\end{proof}

From now on, unless otherwise specified, $R$ is taken to have multiplicative inverses for each non-zero element. $R$ is thus a positive semiring with multiplicative inverses, \textit{i.e.} a so called positive \emph{semifield}. Let us now define $R$-polydependence atoms from $R$-polyindependence atoms, and examine their meaning. 

\begin{Def}
Consider variables $|\bar x^C| = |\bar u^{C'}|$ and $|\bar y^C| = |\bar v^{C'}|$. The syntactical expression $\dep[R](\bar x^C, \bar y^C / \bar u^{C'}, \bar v^{C'})$ is called an \emph{$R$-polydependence atom}. It has the following semantics : 
$$\mathcal L \models\ \dep[R](\bar x^C, \bar y^C / \bar u^{C'}, \bar v^{C'}) \iff \mathcal L \models \bar y^C / \bar y^C \independent[R][\bar x^C, \bar u^{C'} / \bar x^C] \bar v^{C'} / \bar y^C$$
in line with Lemma \ref{le:polydep_atoms_as_polyindep_atoms}. 
\end{Def}

\begin{Le} \label{le:R_dep_atoms_weakly_flat}
For any $R$-polyteam $\mathcal L$ : $\mathcal L \models\ \dep[R](\bar x^C, \bar y^C / \bar u^{C'}, \bar v^{C'}) \iff \supp(\mathcal L) \models\ \dep(\bar x^C, \bar y^C / \bar u^{C'}, \bar v^{C'})$. In words, $\dep[R]$ is essentially a possibilistic notion : it only sees whether weights are zero or non-zero. 
\end{Le}

\begin{proof}
\begin{itemize}
\item[$\implies$.] Apply Proposition \ref{pro:pushing_semiring_truth}. Crucially, one may prove by hand that $=$ and $\dep[\mathbb B] $ coincide.  
\item[$\impliedby$.] Assume $\mathcal L$ is an $R$-polyteam such that $\supp(\mathcal L) \models\ \dep(\bar x^C, \bar y^C / \bar u^{C'}, \bar v^{C'})$. Let $\bar o_{xu}, \bar o_{yv} \subseteq O$ and let us show that $L_C(\bar x\bar y)L_C(\bar x)L_{C'}(\bar u) = L_C(\bar x\bar y)L_{C'}(\bar u\bar v)L_C(\bar x)$. Without loss of generality, $L_C(\bar x\bar y = \bar o_{xu}\bar o_{yv}) \neq 0$. By multiplicative cancellation, this is equivalent to showing $L_{C'}(\bar u) = L_{C'}(\bar u\bar v)$. By positivity, this is the same as saying that $L_{C'}(\bar u\bar v = \bar o_{xu}\bar o'_{yv}) = 0$ as soon as $\bar o'_{yv} \neq \bar o_{yv}$, which is true from the hypothesis and the fact that $L_C(\bar x\bar y = \bar o_{xu}\bar o_{yv}) \neq 0$.
\end{itemize}
\end{proof}

\begin{Le} \label{le:B_indep_atoms_coincidence}
The semantics of $\bar y^C / \bar v^{C''} \perp_{\bar x^C, \bar a^{C'} / \bar u^{C''}} \bar b^{C'} / \bar w^{C''}$ and $\bar y^C / \bar v^{C''} \independent[\mathbb B][\bar x^C, \bar a^{C'} / \bar u^{C''}] \bar b^{C'} / \bar w^{C''}$ coincide as soon as $L_C \restriction \bar x\bar y \simeq L_{C''} \restriction \bar u\bar v$ and $L_{C'} \restriction \bar a\bar b \simeq L_{C''} \restriction \bar u\bar w$, meaning $L_C(\bar x\bar y \mapsto \bar o_{xau}\bar o_{yv}) = L_{C''}(\bar u\bar v \mapsto \bar o_{xau}\bar o_{yv}) \in \mathbb B$ for all $\bar o_{xau}\bar o_{yv} \subseteq O$, and similarly for the other one. 
\end{Le}

\begin{proof}
\begin{itemize}
\item[$\models$.] Assume $\mathcal L \models \bar y^C / \bar v^{C''} \perp_{\bar x^C, \bar a^{C'} / \bar u^{C''}} \bar b^{C'} / \bar w^{C''}$, let $\bar o_{xau}\bar o_{yv}\bar o_{bw} \subseteq O$, and let us show that : 
$$L_{C''}(\bar u\bar v\bar w)L_C(\bar x)L_{C'}(\bar a) = L_C(\bar x\bar y)L_{C'}(\bar a\bar b)L_{C''}(\bar u) \in \mathbb B$$
By positivity, if either $L_C(\bar x)$, $L_{C'}(\bar a)$ or $L_{C''}(\bar u)$ is 0, then the result holds. Therefore, assume all three are 1. Remains to show that :
$$L_{C''}(\bar u\bar v\bar w) = L_C(\bar x\bar y)L_{C'}(\bar a\bar b)$$
$\geq$ is precisely what the $\perp$ atom says. $\leq$ follows easily from $L_C \restriction \bar x\bar y \simeq L_{C''} \restriction \bar u\bar v$ and $L_{C'} \restriction \bar a\bar b \simeq L_{C''} \restriction \bar u\bar w$.
\item[$\delsmo$.] Assume $\mathcal L \models \bar y^C / \bar v^{C''} \independent[\mathbb B][\bar x^C, \bar a^{C'} / \bar u^{C''}] \bar b^{C'} / \bar w^{C''}$, and suppose $s \in L_C$ and $s' \in L_{C'}$ are such that $s(\bar x^C) = s'(\bar a^{C'}) =: \bar o_{xau}$. $L_C \restriction \bar x\bar y \simeq L_{C''} \restriction \bar u\bar v$ (or the other one) implies that $L_{C''}(\bar u = \bar o_{xau}) = 1$, and we already have $L_C(\bar x) = L_{C'}(\bar a) = 1$. If we let $\bar o_{yv} := s(\bar y^C)$ and $\bar o_{bw} := s'(\bar b^{C'})$, then we have $L_C(\bar x\bar y) = L_{C'}(\bar a\bar b) = 1$ as well. Finally, the $\independent[\mathbb B] $ atom says that $L_{C''}(\bar u\bar v\bar w) = 1$, meaning there is $s'' \in L_{C''}$ as wanted. 
\end{itemize}
\end{proof}

\subsection{Logics of dependence and independence in $R$-weighted polyteams} \label{sec:R_logics}
Let us now define full logics of dependence and independence in $R$-weighted polyteams. 

\begin{Def}
Formulas for $R$-polyteam logics are build according to the folowing syntax. 
\begin{center}
$\varphi :=$ Single-labeled first order literal | Non-first order $R$-atom | $\psi \wedge \theta$ | $\psi \vee \theta$ | $\exists x\ \psi$ | $\forall x \psi$
\end{center}
We will also often use $\texists$ and $\tforall$, as before. Again, when working in the physical framework $(\tau, \mathcal M)$, we will use $\texists z$ as a shorthand for $\texists z^{C_1}\ \exists z^{C_2} \cdots \exists z^{C_n}$. If one only allows $R$-polydependence atoms, one gets \emph{$R$-polydependence logic} FO$(\dep[R]\ )$. Else, if one allows all $R$-polyindependence atoms, one gets \emph{$R$-polyindependence logic} FO$(\independent[R]\!)$.
\end{Def}

\begin{Def}
Let $\mathcal L, \mathcal L', \mathcal L_1, \mathcal L_2$ be $R$-polyteams. We say that : 
\begin{itemize}
\item $\mathcal L'$ is a \emph{weak $R$-subpolyteam} of $\mathcal L$ whenever $\supp(\mathcal L') \subseteq \supp(\mathcal L)$. We denote $\mathcal L' \subseteq \mathcal L$. 
\item $\mathcal L'$ is a \emph{strong $R$-subpolyteam} of $\mathcal L$ whenever there is $r \in R$ such that $r \times L'_C(s) = L_C(s)$ for each $C$ and $s$. In that case $r$ is the necessarily-common value of $L_C(\supp(L'_C))$, and we denote $\mathcal L' \subseteq_r \mathcal L$. In words, $\mathcal L'$ is a weak $R$-subpolyteam of $\mathcal L$ and proportions between non-zero weights are preserved from one to the other. 
\item $(\mathcal L_1, \mathcal L_2)$ is a \emph{$(p,q)$-split} of $\mathcal L$ whenever $\mathcal L_1, \mathcal L_2 \subseteq \mathcal L$ and $p,q \in R$ are such that $p + q = 1$ and $p \times L_{1,C}(s) + q \times L_{2,C}(s) = L_C(s)$ for each $C$ and $s$. We denote $p\mathcal L_1 + q\mathcal L_2 = \mathcal L$. 
\item $F : L_C \to \mathcal D_RA_i$ is a \emph{supplement function} for $\mathcal L$, and for any variable $x^{C,i} \notin \dom(L_C)$, we denote by $\mathcal L(F/x)$ the $R$-polyteam such that $\mathcal L(F/x)_C(s(a/x)) = L_C(s) \times F(s)(a)$. When $A_i$ is finite and $F$ is constantly equal to the uniform distribution ($1/|A_i|$ exists by positivity), $F$ is called a \emph{duplicate function}, and we write $\mathcal L(A/x)$ instead of $\mathcal L(F/x)$.  
\end{itemize}
\end{Def}

\begin{Def} $R$-polyteam formulas are given the following semantics. 
\begin{itemize}
\item[(First order literal)] $\mathcal A \models_{\mathcal L}\ \alpha \iffdef \forall s \in \supp(L_C), \mathcal A \models_s\ \alpha$, where $C$ is the uniform label of FV$(\alpha)$
\item[(Non-first order atom)] As detailed in section \ref{sec:R_polyatoms}.
\item[($\wedge$)] $\mathcal A \models_{\mathcal L}\ \varphi \wedge \psi \iffdef \mathcal A \models_{\mathcal L}\ \varphi$ and $\mathcal A \models_{\mathcal L}\ \psi$
\item[($\vee$)] $\mathcal A \models_{\mathcal L}\ \varphi \vee \psi \iffdef \exists \mathcal L_1, \mathcal L_2 \subseteq \mathcal L, \exists p,q \in R,\ p\mathcal L_1 + q\mathcal L_2 = \mathcal L$, $\mathcal A \models_{\mathcal L_1}\ \varphi$ and $\mathcal A \models_{\mathcal L_2}\ \psi$
\item[($\exists$)] $\mathcal A \models_{\mathcal L}\ \exists x^{C,i}\ \varphi \iffdef \exists F : L_C \to \mathcal D_RA_i,\ \mathcal A \models_{\mathcal L(F/x)}\ \varphi$
\item[($\forall$)] $\mathcal A \models_{\mathcal L}\ \forall x^{C,i}\ \varphi \iffdef \mathcal A \models_{\mathcal L(A/x)}\ \varphi$ 
\item[($\texists$)] $\mathcal A \models_{\mathcal L}\ \texists z^{C,h}\ \varphi \iffdef \mathcal B \models_{\mathcal L}\ \exists z^{C,h}\ \varphi$ for some $\mathcal B$ expanding $\mathcal A$ by the sort $h$
\item[($\tforall$)] $\mathcal A \models_{\mathcal L}\ \tforall z^{C,h}\ \varphi \iffdef \mathcal B \models_{\mathcal L}\ \forall z^{C,h}\ \varphi$ for all $\mathcal B$ expanding $\mathcal A$ by the sort $h$
\end{itemize}
\end{Def}

\begin{Rem}
In the semantics of $\vee$, one may wonder whether it is appropriate to impose that the $(p,q)$-split is made of parts that are weak $R$-subpolyteams of $\mathcal L$. It is in fact crucial to do so in order to recover the following equivalence : 
\begin{align*}
\mathcal L \models \varphi \vee \top & \iff \exists \mathcal L' \subseteq \mathcal L,\ \mathcal L' \models \varphi \\ 
& \niff \varphi \not\equiv \top \text{ (as would have happened without the assumption)}
\end{align*}
\end{Rem}

\begin{Ex}
Consider the probabilistic polyteam $\mathcal L$ from Figure \ref{fig:R_polyindependence_atom}, on some single-sorted structure $\mathcal A$ with domain $A = \{0,1,2\}$, interpreting the binary relation symbol $R$ as $R^{\mathcal A} := \{(m,n) \in A^2 : m + n \geq 2\}$, and the following formulas. 
\begin{center}
\begin{minipage}{0.5\textwidth}
$\begin{array}{ccc}
x & y & \mathbb R_+ \\ \hline
0 & 0 & 1/3 \\ 
1 & 1 & 1/3 \\ 
2 & 0 & 1/3
\end{array}$
\hspace{10pt}
$\begin{array}{ccc}
a & b & \mathbb R_+ \\ \hline
0 & 1 & 1/2 \\ 
1 & 0 & 1/3 \\ 
1 & 2 & 1/6
\end{array}$
\hspace{10pt}
$\begin{array}{cccc}
u & v & w & \mathbb R_+ \\ \hline
0 & 0 & 1 & 1/4 \\
1 & 1 & 0 & 1/2 \\ 
1 & 1 & 2 & 1/4
\end{array}$
\end{minipage}
\hspace{15pt}
\begin{minipage}{0.22\textwidth}
\begin{align*}
\alpha &\ :\ R(x,y) \\ 
\beta &\ :\ \dep[\mathbb R_+](x,y / u,v) \\ 
\gamma &\ :\ y/v \independent[\mathbb R_+][x,a/u] b/w
\end{align*}
\end{minipage}
\begin{minipage}{0.22\textwidth}
\begin{align*}
\varphi &\ :\ \alpha \vee \beta\\ 
\psi &\ :\ \gamma \wedge \varphi
\end{align*}
\end{minipage}
\end{center}

\begin{itemize}
\item $\alpha$ is false, because $xy \mapsto (0,0)$ is in $\supp(\mathcal L)$.
\item $\beta$ is true, it suffices to check that $\supp(\mathcal L) \models\ \dep(x,y / u,v)$.
\item $\gamma$ is true, as proved in previous Example.
\item $\varphi$ is true. Consider $\mathcal L_1$ (resp. $\mathcal L_2$) obtained from $\mathcal L$ by changing the $xy$-weights to $0, 1/3, 2/3$ (resp. $2/3, 1/3, 0$). Then $(\mathcal L_1, \mathcal L_2)$ is a suitable $(1/2,1/2)$-split. 
\item $\psi$ is true. 
\end{itemize}
\end{Ex}

\begin{Fac}
As for possibilistic logics, the distributivity of $\vee$ and $\wedge$ on one another fail. 
\end{Fac}

\begin{proof}
(Partial.) Consider $\mathcal L$ and $\alpha$ as in previous example. Let $\top$ denote any tautology, such as $x=x$. Then, $\mathcal L \models (\alpha \vee \top) \wedge (\neg \alpha \vee \top)$, but $\mathcal L \not\models (\alpha \wedge \neg\alpha) \vee \top$ (there is no such thing as an empty $\mathbb R_+$-polyteam). Note that the former is almost a tautology whereas the latter is a self-contradiction. 
\end{proof}

\begin{Le}
$R$-polyteam logics support domain locality, as in the Corollary of page \pageref{cor:domain_locality}.
\end{Le}

\begin{proof}
Replay the proof of Lemma \ref{le:locality}. The only major difference is that the final union should be changed to a sum. 
\end{proof}

\subsection{Collapses and realizations} \label{sec:collapses_realizations} 
Results such as Lemma \ref{le:R_dep_atoms_weakly_flat} incite us to investigate the links between satisfaction in different semirings (positive semifields, rather). That is the purpose of this section. 

\begin{Def}
Let $R,R'$ be positive semifields, and $\varphi$ be a formula in FO$(\independent[R]\!)$. Denote by $\varphi^{R'}$ the formula in FO($\independent[R']$\!) obtained from $\varphi$ by replacing each occurrence of $\independent[R]$ with $\independent[R']$ (and, each occurrence of $\dep[R]$ with $\dep[R']$). Denote $\varphi^\perp$ if one replaces with $\perp$.
\end{Def}

\begin{Def}
Let $f : R \to R'$ be a semifield homomorphism, and $\mathcal L$ be an $R$-polyteam. The \emph{$f$-collapse} of $\mathcal L$ is the $R'$-polyteam $f(\mathcal L) = \{L^f_C, C \in \mathcal M\}$ defined by : 
$$L^f_C(s) = f(L_C(s)) \in R' \text{ , for each } C \text{ and } s$$
\end{Def}

\begin{Ex}
$\supp(\mathcal L)$ can be thought of as the $\supp$-collapse of $\mathcal L$, where $\supp : \begin{array}[t]{rcl} R & \longrightarrow & \mathbb B \\ r & \longmapsto & \begin{cases} 1 & \text{if } r \neq 0 \\ 0 & \text{if } r = 0 \end{cases} \end{array}$
\end{Ex}

\begin{Rem}
Any semifield homomorphism $f : R \to R'$ is such that $f^{-1}(0_{R'}) = \{0_R\}$. Indeed, if $r \neq 0$, then $\exists r' \in R, rr' = 1_R$, and $f(r)f(r') = 1_{R'}$, whence $f(r) \neq 0_{R'}$. Therefore, $\supp\ f(\mathcal L) = \supp(\mathcal L)$. 
\end{Rem}

\begin{Pro} \label{pro:pushing_semiring_truth}
Let $f : R \to R'$ be a semifield homomorphism, $\mathcal L$ an $R$-polyteam, and $\varphi$ in FO$(\independent[R]\!)$. We have : 
$$\mathcal L \models \varphi \implies f(\mathcal L) \models \varphi^{R'}$$ 
\end{Pro}

\begin{proof}
The proof goes by induction on $\varphi$. 
\begin{itemize}
\item[FO.] If $\varphi$ is a first order atom, the result follows from the fact $\supp\ f(\mathcal L) = \supp(\mathcal L)$. 
\item[$\independent \!\!$.] If $\varphi$ is the atom $\bar y^C / \bar v^{C''} \independent[R][\bar x^C, \bar a^{C'} / \bar u^{C''}] \bar b^{C'} / \bar w^{C''}$, and $\mathcal L \models \varphi$, then we have : 
\begin{align*}
L_{C''}(\bar u\bar v\bar w)L_C(\bar x)L_{C'}(\bar a) & = L_C(\bar x\bar y)L_{C'}(\bar a\bar b)L_{C''}(\bar u) \hspace{20pt} \text{ for any } \bar o_{xau}\bar o_{yv}\bar o_{bw} \subseteq O \text{, and thus} \\ 
L^f_{C''}(\bar u\bar v\bar w)L^f_C(\bar x)L^f_{C'}(\bar a) & = L^f_C(\bar x\bar y)L^f_{C'}(\bar a\bar b)L^f_{C''}(\bar u) \hspace{20pt} \text{ after applying } f
\end{align*}
\item[$\wedge$.] Easy verification. 
\item[$\vee$.] Assume $\varphi = \psi \vee \theta$, and fix $\mathcal L_1, \mathcal L_2 \subseteq \mathcal L$, $p,q \in R$ such that : 
$$p + q = 1_R \hspace{20pt} p\mathcal L_1 + q\mathcal L_2 = \mathcal L \hspace{20pt} \mathcal L_1 \models \psi \hspace{10pt} \text{and} \hspace{10pt} \mathcal L_2 \models \theta$$
Applying $f$, its homomorphism regularities, and using the induction hypothesis, we get : 
$$f(p) + f(q) = 1_{R'} \hspace{20pt} f(p)f(\mathcal L_1) + f(q)f(\mathcal L_2) = f(\mathcal L) \hspace{20pt} f(\mathcal L_1) \models \psi^{R'} \hspace{10pt} \text{and} \hspace{10pt} f(\mathcal L_2) \models \theta^{R'}$$
In particular $f(\mathcal L) \models \varphi^{R'}$.
\item[$\exists$.] Assume $\varphi = \exists x^{C,i}\ \psi$, and $F : \supp(L_C) \to \mathcal D_RA_i$ is such that $\mathcal L(F/x) \models \psi$. By induction hypothesis : $f(\mathcal L(F/x)) \models \psi^{R'}$. One can check that $f(\mathcal L(F/x)) = f(\mathcal L)(f^*F/x)$, where :
$$f^*F(s)(a) := f(F(s)(a)) \text{ for each } s \in \supp\ f(\mathcal L) = \supp(\mathcal L) \text{ and } a \in A_i$$ 
In particular $f(\mathcal L) \models \varphi^{R'}$.
\item[$\forall$.] Similar. 
\end{itemize}
\end{proof}

\begin{Cor} \label{cor:weak_flatness_easy_direction}
Let $\mathcal L$ be an $R$-polyteam and $\varphi$ in FO$(\independent[R]\!)$.
$$\mathcal L \models \varphi \implies \supp(\mathcal L) \models \varphi^{\mathbb B} \overset{\tiny{\begin{array}{c}\text{under hyp. of} \\ \text{Lemma \ref{le:B_indep_atoms_coincidence}}\end{array}}}{\iff}\normalsize \supp(\mathcal L) \models \varphi^\perp$$
\end{Cor}

\begin{proof}
The result follows from the combination of previous Proposition and Lemma \ref{le:B_indep_atoms_coincidence}.
\end{proof}

Let us now investigate in which cases satisfaction in the possibilistic collapse suffices for truth to occur in $R$-polyteams logics, \textit{i.e.} as in Lemma \ref{le:R_dep_atoms_weakly_flat}. 

\begin{Def}
$\varphi$ in FO$(\independent[R]\!)$ is called \emph{weakly flat} whenever $\mathcal L \models \varphi \iff \supp(\mathcal L) \models \varphi^{\mathbb B}$. 
\end{Def}

\begin{Fac}
First order atoms and $R$-polydependence atoms are weakly flat. 
\end{Fac}

\begin{proof}
First order atoms are weakly flat by very definition. For $R$-polydependence atoms, it is the content of Lemma \ref{le:R_dep_atoms_weakly_flat}. 
\end{proof}

\begin{Fac}
Weak flatness is not preserved by $\vee$. 
\end{Fac}

\begin{proof}
Both $R(x^C,y^C)$ and $R(u^{C'},v^{C'})$ are weakly flat, but $\varphi = R(x,y) \vee R(u,v)$ is not. Indeed, let $\mathcal A$ be a single-sorted structure with domain $\{0,1\}$, interpreting $R$ as $R^{\mathcal A} := \{(m,n) \in A^2 : m + n \geq 2\}$, and $\mathcal L$ be the following $\mathbb R_+$ polyteam. 
$$\begin{array}{ccc}
x & y & \mathbb R_+ \\ \hline
0 & 0 & 3/4 \\ 
1 & 1 & 1/4
\end{array}
\hspace{20pt}
\begin{array}{ccc}
u & v & \mathbb R_+ \\ \hline
0 & 0 & 3/4 \\ 
1 & 1 & 1/4
\end{array}$$
Then $\supp(\mathcal L) \models \varphi^{\mathbb B}$. However, if $\mathcal L \models \varphi$ occurred, it would forcedly be witnessed by the $(p,q)$-split : 
$$\mathcal L_1 =\hspace{10pt}\begin{array}{ccc}
x & y & \mathbb R_+ \\ \hline
1 & 1 & 1
\end{array}
\hspace{20pt}
\begin{array}{ccc}
u & v & \mathbb R_+ \\ \hline
0 & 0 & 1
\end{array}
\hspace{20pt}\mathcal L_2 =\hspace{10pt}\begin{array}{ccc}
x & y & \mathbb R_+ \\ \hline
0 & 0 & 1
\end{array}
\hspace{20pt}
\begin{array}{ccc}
u & v & \mathbb R_+ \\ \hline
1 & 1 & 1
\end{array}$$
But then $p\mathcal L_1 + q\mathcal L_2 = \mathcal L$ yields $p = 1/4$ (evaluated in $xy \mapsto (1,1)$), and $p = 3/4$ (evaluated in $uv \mapsto (0,0)$), a contradiction. Note that this proof can be replayed in any semiring $R$ such that $\exists r, r' \in R, r,r' \leq 1_R$ and $r \neq r'$. 
\end{proof}

\begin{Rem}
On the contrary (Lemma 4.11 in \cite{abramsky_team_2026}), $\vee$ does preserve weak flatness in $R$-\emph{team} semantics. A sufficient condition for $\supp(\mathcal L) \models \psi^{\mathbb B} \vee \theta^{\mathbb B}$ to imply $\mathcal L \models \psi \vee \theta$ is that $\supp(\mathcal L) \models \psi^{\mathbb B} \vee \theta^{\mathbb B}$ occurs with respect to a split $(\mathcal S_1, \mathcal S_2)$ such that $L_C(\mathcal S_i)$ is a constant when $C$ varies, for both $i = 1$ and $i = 2$. In that case, a proof similar to that of \cite{abramsky_team_2026} can be settled. 
\end{Rem}

\begin{Le} \label{le:weak_flatness_regularities}
Weak flatness is preserved by $\wedge$ and quantifiers. 
\end{Le}

\begin{proof}
\begin{itemize}
\item[$\wedge$.] Easy verification. 
\item[$\exists$.] Assume $\supp(\mathcal L) \models \exists x^{C,i}\ \psi^{\mathbb B}$, with $\psi$ weakly flat. Fix $F : \supp(L_C) \to \mathcal P^+(A_i)$ such that $\supp(L)(F/x) \models \psi^{\mathbb B}$. For each $s(a/x) \in \supp(L_C)(F/x)$, define 
$$\tilde F(s)(a) := \frac{L_C(s)}{|F(s)|}$$
Then $\supp(\mathcal L)(F/x) = \supp(\mathcal L(\tilde F/x))$, by weak flatness of $\psi$ we have $\mathcal L(\tilde F/x) \models \psi$, and therefore $\mathcal L \models \exists x\ \psi$.
\item[$\forall$.] Similar. 
\end{itemize}
\end{proof}

\begin{Pro} \label{pro:R_dep_logic_weakly_flat}
The splitjunction-free fragment of FO$(\dep[R]\ )$ is weakly flat. 
\end{Pro}

\begin{proof}
Easy induction using Lemma \ref{le:R_dep_atoms_weakly_flat} and Lemma \ref{le:weak_flatness_regularities}. 
\end{proof}

\begin{Def}
Let $\varphi$ in FO$(\independent[R])$. $\varphi$ is called \emph{downwards closed} (or \emph{$\downarrow$-closed}) whenever $\varphi$ is stable by strong subpolyteams, in the following sense : if $\mathcal L \models \varphi$ and $\mathcal L'$ is a strong subpolyteam of $\mathcal L$, then $\mathcal L' \models \varphi$. $\varphi$ is called \emph{strongly downwards closed} (or \emph{strongly-$\downarrow$-closed}) whenever $\varphi$ is stable by weak subpolyteams. 
\end{Def}

\begin{Le} \label{le:weakly_flat_downwards_closed}
If $\varphi$ in FO$(\independent[R])$ is weakly flat and $\varphi^{\mathbb B}$ is $\downarrow$-closed, then $\varphi$ is strongly-$\downarrow$-closed. 
\end{Le}

\begin{proof}
Assume that $\varphi$ is as in the statement, $\mathcal L \models \varphi$ and $\mathcal L'$ is a weak subpolyteam of $\mathcal L$. Then $\supp(\mathcal L') \subseteq \supp(\mathcal L)$ and $\supp(\mathcal L) \models \varphi^{\mathbb B}$ by the Corollary at page \pageref{cor:weak_flatness_easy_direction}. By $\downarrow$-closedness, $\supp(\mathcal L') \models \varphi^{\mathbb B}$, and by weak flatness $\mathcal L' \models \varphi$. 
\end{proof}

\begin{Pro}
The splitjunction-free fragment of FO$(\dep[R]\ )$ is strongly-$\downarrow$-closed. 
\end{Pro}

\begin{proof}
The result follows from Lemma \ref{le:weakly_flat_downwards_closed}, Proposition \ref{pro:R_dep_logic_weakly_flat}, and Proposition \ref{pro:polyteam_logics_closures}.
\end{proof}

\subsection{Meaning of physical $R$-properties} \label{sec:R_properties_meaning} 
There is no \textit{a priori} reason that the $R$-transcriptions of physical properties from $\mathcal P$ keep a natural, somewhat faithful meaning. Still, they do. 

\begin{Pro}
Let $\mathcal L$ (resp. $\mathcal H$) be an empirical (resp. hidden variable) $R$-polyteam. 
\begin{align*}
\mathcal L \models NS^R & \iff \mathcal L \text{ is locally consistent} \\ 
\mathcal L \models WD^R & \iff \supp(\mathcal L) \models WD\\ 
\mathcal L \models SD^R & \iff \supp(\mathcal L) \models SD\\ 
\mathcal H \models HWD^R & \iff \supp(\mathcal H) \models HWD\\ 
\mathcal H \models HSD^R & \iff \supp(\mathcal H) \models HSD\\ 
\mathcal H \models SV^R & \iff \supp(\mathcal H) \models SV\\
\mathcal H \models PI^R \wedge zI^R & \iff \mathcal H \text{ is locally consistent}\\ 
\mathcal H \models Loc^R & \overset{\tiny{\begin{array}{c}\text{under}\\ zI^R\end{array}}\normalsize}{\iff} H_C(\bar uz^C = \bar o\lambda)H_C(z^C = \lambda)^{k-1} = \prod_{j<k} H_{C_j}(u_jz^{C_j} = o_j\lambda) \\
& \hspace{20pt} \text{for all } u_0 \cdots u_{k-1} \subseteq C \in \mathcal M, (C_i)_{i<k} \text{ such that } u_i \in C_i, \bar o \in O^k \text{ and } \lambda \in \Lambda 
\end{align*}
\end{Pro}

\begin{proof}
The second, third, fourth, fifth and sixth equivalences are direct applications of Proposition \ref{pro:R_dep_logic_weakly_flat}. 
\begin{itemize}
\item[$NS$.] $\begin{array}[t]{rcl}
\mathcal L \models NS^R & \iff & \forall U \subseteq C,C' \in \mathcal M, \mathcal L \models \bar u^C / \bar u^{C'} \independent[R] \emptyset^C / \emptyset^{C'}\\ 
& \iff & \forall U \subseteq C,C'\in \mathcal M, \forall \bar o_u \subseteq O, L_{C'}(\bar u) = L_C(\bar u)\\ 
& \iff & \mathcal L \text{ is locally consistent}
\end{array}$
\item[$PI \wedge zI$.] We omit it here. 
\item[$Loc$.] \begin{itemize}
\item[$\implies$.] Assume $\mathcal H \models Loc^R \wedge zI^R$, and let $\bar u, C, (C_i)_{i<k}, \bar o$ and $\lambda$ be as in the statement. If $H_C(z^C = \lambda) = 0$, then $H_{C_j}(z^{C_j} = \lambda) = 0$ as well by $z$-Independence, and the goal is true. Assume the contrary. If $H_{C_{j_0}}(u_{j_0}z^{C_{j_0}} = o_{j_0}\lambda) = 0$, then $o_{j_0} \notin \supp(\mathcal H)_{C_{j_0}}^\lambda[u_{j_0}]$, and since $\supp(\mathcal H) \models Loc$, that implies $\bar o \notin \supp(\mathcal H)_{C}^\lambda[\bar u^C]$, \textit{i.e.} $H_C(\bar uz^C = \bar o\lambda) = 0$ and the goal holds. We may thus assume without loss of generality that all factors in the desired equality are non-zero. Fix $i < k$. $PI^R$ says : 
$$H_C(u_iz = o_i\lambda)H_{C_i}(z = \lambda)H_C(z = \lambda) = H_{C_i}(u_iz = o_i\lambda)H_{C_i}(z = \lambda)H_C(z = \lambda)$$
\textit{i.e.} $H_C(u_iz = o_i\lambda) = H_{C_i}(u_iz = o_i\lambda) \neq 0$. Assume now that $H_C(u_0 \cdots u_iz = o_0 \cdots o_i\lambda) \neq 0$ and $i < k-1$. $OI^R$ says : 
$$\mathcal H \models \bar u^C_{\leq i} / \bar u^C_{\leq i} \independent[R][z^C,z^C/z^C] u_{i+1}\bar v^C / u_{i+1}\bar v^C \text{ where } \bar v = C \setminus \{u_0, \cdots, u_{i+1}\} \ \ ,i.e.$$
$$H_C(\bar u_{\leq i+1}\bar vz = \bar o_{\leq i+1}\bar p\lambda)H_C(z = \lambda)^2 = H_C(\bar u_{\leq i}z = \bar o_{\leq i}\lambda)H_C(u_{i+1}\bar vz = o_{i+1}\bar p\lambda)H_C(z = \lambda)$$ 
Summing over $\bar p \subseteq O$, one gets : 
\begin{align*}
H_C(\bar u_{\leq i+1}z = \bar o_{i+1}\lambda) & = \sum_{\bar p} \frac{H_C(\bar u_{\leq i}z = \bar o_{\leq i}\lambda)}{H_C(z = \lambda)} H_C(u_{i+1}\bar vz = o_{i+1}\bar p\lambda) \\ 
& = \frac{H_C(u_{i+1}z = o_{i+1}\lambda)}{H_C(z = \lambda)} H_C(\bar u_{\leq i}z = \bar o_{\leq i}\lambda) \\ 
& = \frac{H_{C_{i+1}}(u_{i+1}z = o_{i+1}\lambda)}{H_C(z = \lambda)} H_C(\bar u_{\leq i}z = \bar o_{\leq i}\lambda) \text{ by what we proved before}
\end{align*}
We easily conclude by induction on $0 \leq i \leq k-1$. 
\item[$\impliedby$.] Assume $\mathcal H \models zI^R$ and the right hand side of the equivalence holds. Let us show $\mathcal H \models PI^R$ first. Let $U \subseteq C, C' \in \mathcal M$, $\bar o \subseteq O$ and $\lambda \in \Lambda$. We need to show 
$$H_{C'}(\bar uz = \bar o\lambda)H_C(z = \lambda)H_{C'}(z = \lambda) = H_C(\bar uz = \bar o\lambda)H_{C'}(z = \lambda)^2$$
which trivially holds if $H_C(z = \lambda) \overset{zI^R}{=} H_{C'}(z = \lambda) = 0$. In the opposite case, we are left to show $H_{C'}(\bar uz = \bar o\lambda) = H_C(\bar uz = \bar o\lambda)$. For this, it suffices to develop them as products, with $C_j = C$ all the time, and to note that the denominators are equal by $z$-Independence. Whence $\mathcal H \models PI^R$. Now, let us show that $\mathcal H \models OI^R$. Let $U \subseteq C \in \mathcal M$, $V := C \setminus U$, $\bar o, \bar p \subseteq O$ and $\lambda \in \Lambda$. We need to show that : 
$$H_C(\bar u\bar vz = \bar o\bar p\lambda)H_C(z = \lambda)^2 = H_C(\bar uz = \bar o\lambda)H_C(\bar vz = \bar p\lambda)H_C(z = \lambda)$$
This follows directly from developping $H_C(\bar u\bar vz = \bar o\bar p\lambda)$, $H_C(\bar uz = \bar o\lambda)$ and $H_C(\bar vz = \bar p\lambda)$ in products, with $C_j = C$ all the time. 
\end{itemize}
\end{itemize}
\end{proof}

\begin{Rem}
We can reformulate the last equivalence as : 
$$\mathcal H \models Loc^R \overset{\tiny{\begin{array}{c}\text{under}\\ zI^R\end{array}}\normalsize}{\iff} H_C(\bar u\ |\ z) = \prod_{j<k} H_{C_j}(u_j\ |\ z)$$
from which we see that $Loc^R$ naturally generalizes its possibilistic version. We do not know whether it is possible to find a counter example if $\mathcal H \not\models zI^R$. 
\end{Rem}

\begin{Pro} \label{pro:R_global_properties_meaning}
Let $\mathcal L$ be an empirical $R$-polyteam. 
\begin{align*}
\mathcal L \models GS^R & \iff \mathcal L \text{ has a global section}\\ 
\mathcal L \models NSC^R & \iff \supp(\mathcal L) \models NSC
\end{align*}
\end{Pro}

\begin{proof}
\begin{itemize}
\item[$GS$.] Assume $\mathcal H \models Det^R$ realizes $\mathcal L$. Then $\supp(\mathcal H) \models HSD \wedge zI$. Define the hidden variable team $Y \in \mathcal D_R\mathcal E^h(M)$ by : 
$$Y(s) = \begin{cases} \text{the common value of } H_C(z = s(z)) & \text{if } s \text{ is consistent with } \mathcal H\\ 0 & \text{otherwise} \end{cases}$$
and let $X$ denote the empirical team realized by $Y$. One can easily check that $X$ is a global section for $\mathcal L$. Conversely, assume $\mathcal L$ has a global section $X$. Let $Y$ be a hidden variale realization of $X$ such that $Y \models \dep[R](z,\bar m)$ (for example $\Lambda = X$). Let $\mathcal H$ denote the contextualization of $Y$ to $\mathcal M$. It is easy to check that $\mathcal H$ realizes $\mathcal L$, and $\mathcal H \models Det^R$.
\item[$NSC$.] Recall that $NS^{\mathbb B} \equiv NS$, because both are equivalent to local consistency. Then : 
\begin{align*} 
\mathcal L \models NSC^R & \iff \exists \mathcal L' \subseteq \mathcal L, \mathcal L' \models WD^R \wedge NS^R\\ 
& \iff \supp(\mathcal L) \text{ has a subpolyteam with } WD \wedge NS\\ 
& \iff \supp(\mathcal L) \models NSC 
\end{align*}
\end{itemize}
\end{proof}

\subsection{Entailments, and Bell's original theorem} \label{sec:R_entailments}
There is no \textit{a priori} reason that the entailments of Section \ref{sec:entailments} remain in $R$-polyteam semantics. Still, they mostly do. 

\begin{Le} \label{le:R_entailments}
$HWD^R \models OI^R$ and $HSD^R \equiv PI^R \wedge HWD^R$.
\end{Le}

\begin{proof}
Omitted. 
\end{proof}

\begin{The} \label{th:rational_prob}
An empirical $\mathbb Q_+$-model has a deterministic hidden variable realization if and only if it has a locally realistic hidden variable realization. In formulas : 
$$\texists z\ Det^{\mathbb Q_+} \overset{\text{not}}{\equiv} \texists z\ HSD^{\mathbb Q_+} \wedge zI^{\mathbb Q_+} \equiv \texists z\ PI^{\mathbb Q_+} \wedge OI^{\mathbb Q_+} \wedge zI^{\mathbb Q_+} \overset{\text{not}}{\equiv} \texists z\ LocReal^{\mathbb Q_+}$$
\end{The}

\begin{proof}
Again, we get inspiration from \cite{albert_unifying_2022}. The $\models$ direction is easy with Lemma \ref{le:R_entailments}. Let us prove $\delsmo$. Assume $\mathcal H \models LocReal^{\mathbb Q_+}$ realizes $\mathcal L$. As in the possibilistic version, we have : 
$$\mathcal H^\lambda[\bar c] = \prod_{m \in C} \mathcal H^\lambda[m]$$
where $\mathcal H^\lambda[m] := \{o \in O : \exists C \in \mathcal M, H_C(mz = o\lambda) \neq 0\}$. Without loss of generality, assume $\Lambda = \mathcal H[z]$, and note that $\mathcal H^\lambda[\bar c] \neq \emptyset$ by $z$-Independence. Define $Z_m := \mathcal H[m] \times \Lambda$ for each $m \in M$. By Parameter Independence and $z$-Independence (\textit{i.e.} local consistency), $H_C(m = o\ |\ z = \lambda)$ does not depend on $C$. Thus, the following definition makes sense : 
$$H_C(m = o\ |\ z= \lambda) =: p_m(o\lambda) = \frac{a_m(o\lambda)}{b_m(o\lambda)} \in \mathbb Q_+ \text{ where } a_m(o\lambda), b_m(o\lambda) \in \mathbb N \text{ are coprime}$$
Let $N_m := \gcd\ \{b_m(o\lambda), o\lambda \in Z_m\}$ and $\Lambda_m$ be any set of size $N_m$. Define the following. 
\begin{align*}
\Lambda' & := \Lambda \times \prod_{m \in M} \Lambda_m\\ 
N_m(o\lambda) & := p_m(o\lambda)N_m \in \mathbb N \\ 
\Lambda_m & =: \bigsqcup_{o \in \mathcal H[m]} \Lambda_m(o\lambda) \text{ any partition such that } |\Lambda(o\lambda)| = N_m(o\lambda)
\end{align*}
Such a partition exists for $\displaystyle \sum_{o \in \mathcal H^\lambda[m]} N_m(o\lambda) = N_m \sum_{m \in \mathcal H^\lambda[m]} p_m(o\lambda) = N_m$. Finally, define $\mathcal H'$ as follows. 

\begin{multline*}
H'_C(\bar c z\bar \zeta = \bar o\lambda\bar \mu) = \frac{H_C(\bar cz = \bar o\lambda)}{\displaystyle \prod_{m \in C} p_m(o_m\lambda) \cdot \prod_{m \in M} N_m} \ \text{ if } \bar \mu \text{ is such that } \forall m \in C, \mu_m \in \Lambda(o_m\lambda), \\ \text{and } 0 \text{ otherwise}
\end{multline*}
where the hidden variable is conveniently denoted $z\bar \zeta$ instead of $z$, so that it matches the definition of $\Lambda'$. This does define a probability distribution for each $C$ because the denominator is precisely the number of $\bar \mu$'s satisfying the condition. Firstly, let us show $\mathcal H' \models HSD^{\mathbb Q_+}$. Let $m \in C, C' \in \mathcal M$, $\lambda\bar \mu \in \Lambda'$, $o \in O$, and suppose without loss of generality that both $H'_C(mz\bar \zeta = o\lambda\bar \mu)$ and $H'_{C'}(mz\bar \zeta = o'\lambda\bar \mu)$ are non-zero. By definition then, $\mu_m \in \Lambda_m(o\lambda)$ and $\mu_m \in \Lambda_m(o'\lambda)$, which forces $o = o'$ since otherwise these sets would be disjoint. Thus $\mathcal H' \models HSD^{\mathbb Q_+}$. Secondly, let us show $\mathcal H' \models zI^{\mathbb Q_+}$. Let $C \in \mathcal M$ and $\lambda\bar \mu \in \Lambda'$. By construction, there is a unique $\bar o \subseteq O$ such that $\forall m \in C, \mu_m \in \Lambda_m(o_m\lambda)$. Thus : 
$$H'_C(z\bar \zeta = \lambda\bar \mu) = H'_C(\bar cz\bar \zeta = \bar o\lambda\bar \mu) = \frac{H_C(\bar cz = \bar o\lambda)}{\displaystyle \prod_{m \in C} p_m(o_m\lambda) \cdot \prod_{m \in M} N_m}$$
Crucially, we have : 
$$\displaystyle \prod_{m \in C} p_m(o_m\lambda) = \prod_{m \in C} H_C(m = o_m\ |\ z = \lambda) = H_C(\bar c = \bar o\ |\ z = \lambda)$$ by Outcome Independence of $\mathcal H$. Thus, substituting in previous equation, we get : 
$$H'_C(z\bar \zeta = \lambda\bar \mu) = \frac{H_C(z = \lambda)}{\displaystyle \prod_{m \in M} N_m}$$
And that does \emph{not} depend on $C$ by $z$-Independence of $\mathcal H$. Hence $\mathcal H' \models zI^{\mathbb Q_+}$. Finally, $\mathcal H'$ realizes $\mathcal L$, essentially because it is a "$\bar \zeta$-realization" of $\mathcal H$.
\end{proof}

\begin{Rem}
This proof strongly relies on the combinatorial nature of $\mathbb Q_+$, and one hardly sees how it could be adapted into a proof for $\mathbb R_+$. We have extensively tried to find a proof for $\mathbb R_+$ by "$\mathbb Q_+$-approximation", without success.  
\end{Rem}

Since Bell's theorem is about quantum mechanics, and that $QM^{\mathbb R_+}$ does not make sense \textit{a priori}, let us define it before stating Bell's original theorem. 

\begin{Def}
Let $\mathcal L$ be an $\mathbb R_+$-polyteam. We say that $\mathcal L$ is a \emph{finite dimensional quantum mechanical polyteam} of type $(\mathcal M,O)$ whenever there exists $\mathcal S$ a finite dimensional quantum system of type $(\mathcal M, O)$ (as in Definition \ref{def:finite_dim_quantum_system}) such that $L_C(s) = p^{\mathcal S}_C(s(\bar c))$ for each $C$ and $s$. We denote it by $\mathcal L \models QM^{\mathbb R_+}$. 
\end{Def}

\begin{Fac}
A possibilistic polyteam satisfies $QM$ exactly when it has a probabilistic realization satisfying $QM^{\mathbb R_+}$. 
\end{Fac}

\begin{The}[Bell 1964, logical version] \label{th:bell_original}
$QM^{\mathbb R_+} \not\models GS^{\mathbb R^+}$. In words, there is a non-local quantum system.
\end{The}

\begin{proof}
The original proof may be found in \cite{bell_einstein_1964}. At this stage of our work, we have a very easy proof. If by contradiction $QM^{\mathbb R_+} \models GS^{\mathbb R_+}$, let $\mathcal L$ be an empirical possibilistic polyteam such that $\mathcal L \models QM$. $\mathcal L$ has a probabilistic realization $\mathcal L'$ such that $\mathcal L' \models QM^{\mathbb R_+}$. As such, $\mathcal L' \models GS^{\mathbb R_+}$ and we can fix $X'$ witnessing it. It is then easy to check that $X := \supp(X')$ is a global section for $\mathcal L$. All in all, $QM \models GS$, contradicting Theorem \ref{th:bell_possibilistic}. 
\end{proof}

\begin{Rem}
The proof confirms that possibilistic Bell's theorem is "stronger" than its probabilistic analogue. 
\end{Rem}

\subsection{Contextuality strength hierarchy}
We have now seen sufficiently many witnesses of contextuality to order them in "wildness", in line with \cite{abramsky_sheaf-theoretic_2011}. 

\begin{Le}
Let $\mathcal L$ be an empirical probabilistic model. Then : 
$$\mathcal L \models GS^{\mathbb R_+} \implies \supp(\mathcal L) \models GS \implies \mathcal L \models NSC^{\mathbb R_+}$$
In words, strong contextuality is stronger than possibilistic non-locality, which is stronger than probabilistic non-locality. 
\end{Le}

\begin{proof}
The first implication is a direct application of the Corollary from page \pageref{cor:weak_flatness_easy_direction}. The second implication follows from Proposition \ref{pro:R_global_properties_meaning} and the fact that any global section of $\supp(\mathcal L)$ is made of global assignments that are consistent with $\mathcal L$. 
\end{proof}

\begin{Fac}
Figure \ref{fig:contextuality_strength_hierarchy} shows where the various examples studied in this work are located in the contextuality strength hierarchy. 
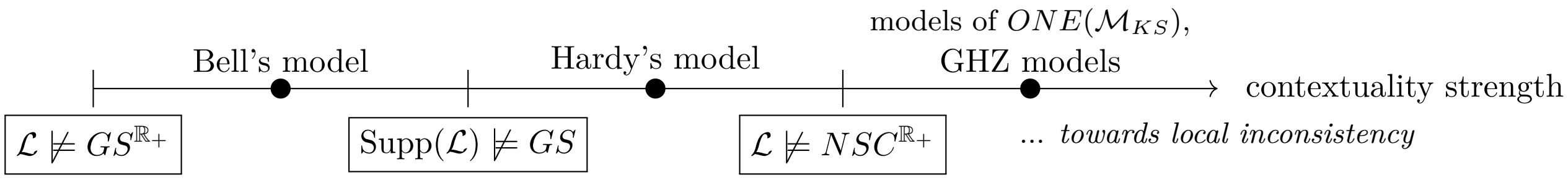
\begin{figure}
\centering
\begin{tikzpicture}
\draw[->] (0,0) -- (12,0);
\node at (14,0) {contextuality strength}; 
\draw (0,0.2) -- (0,-0.2);
\node at (0,-0.6) {$\boxed{\mathcal L \not\models GS^{\mathbb R_+}}$};  
\draw (4,0.2) -- (4,-0.2);
\node at (4,-0.6) {$\boxed{\supp(\mathcal L) \not\models GS}$};  
\draw (8,0.2) -- (8,-0.2);
\node at (8,-0.6) {$\boxed{\mathcal L \not\models NSC^{\mathbb R_+}}$};  
\filldraw (2,0) circle (0.1); 
\node at (2,0.3) {Bell's model}; 
\filldraw (6,0) circle (0.1); 
\node at (6,0.3) {Hardy's model}; 
\filldraw (10,0) circle (0.1); 
\node at (10,0.3) {GHZ models}; 
\node at (10,0.7) {\small{models of $ONE(\mathcal M_{KS})$,}};
\node at (12,-0.5) {\small{... \textit{towards local inconsistency}}};
\end{tikzpicture}
\caption{Contextuality strength hierarchy}
\label{fig:contextuality_strength_hierarchy}
\end{figure}
\end{Fac}

\begin{proof}
The position of Bell's model (Figure \ref{fig:bell_model}) is covered by Theorem \ref{th:bell_original} and analysis of page \pageref{fig:ex1}. That of Hardy's model (Figure \ref{fig:hardy}) is covered Proposition \ref{pro:not_sheaf} and the fact that the global assignment $aba'b' \mapsto 0010$ is consistent. The position of models of $ONE(\mathcal M_{KS})$ is explained by Lemma \ref{le:ks_lemma_1}. Finally, for the position of GHZ models, we refer the reader to \cite{abramsky_sheaf-theoretic_2011} Proposition 6.1.
\end{proof}

\subsection*{Questions}

\begin{enumerate}[label = (\arabic*)]
\item How much of Section \ref{sec:definitions} remains true of semiring-weighted polyteams? Are there normal forms, a sort of compacity? Can we compare expressive powers? 
\item How much of Section \ref{sec:semirings} remains true if we drop the positivity assumption? 
\item How much expressive power upon locally consistent polyteams does one lose if using locally consistent polyteam logics instead of regular polyteam logics? 
\item Does Theorem \ref{th:rational_prob} extend to $\mathbb R_+$? If not, what makes a counter example special? 
\item Can we find large definable subsets of polyteams satisfying $QM$ or $QM^{\mathbb R_+}$?
\item How much of polyteam logics can be simulated within other logics? And conversely? 
\item What to think of the cohabitation of $\perp$ and $\independent[\mathbb B] $? Is one of them more legitimate in certain contexts? 
\item In Section \ref{sec:semirings}, what remains or fails if $(p,q)$-splits are not assumed to be weak subpolyteams of the initial polyteam, and what if they are further assumed to be strong subpolyteams? 
\item Is there a procedure to computationally determine whether $\mathcal L \models QM$? 
\item Can we find partial axiomatizations for polyindependence logic? 
\end{enumerate}

We have gathered several results around the model checking problem and the implication problem of polyteam logics in another paper that we hope to publish soon. 

\printbibliography

\end{document}